# EXERCICES DE STYLE: A HOMOTOPY THEORY FOR SET THEORY, I

MISHA GAVRILOVICH AND ASSAF HASSON

Abstract. We construct a model category (in the sense of Quillen) for set theory, starting from two arbitrary, but natural, conventions. It is the simplest category satisfying our conventions and modelling the notions of finiteness, countability and infinite equi-cardinality. In a subsequent paper [GH10] we give a homotopy theoretic dictionary of set theoretic concepts, most notably Shelah's covering number $\mathrm{cov}(\lambda, \aleph_1, \aleph_1, 2)$, recovered from this model category.

We argue that from the homotopy theoretic point of view our construction is essentially automatic following basic existing methods, and so is (almost all) the verification that the construction works.

*the mischief of it is, nature will have to take its course:*

*Miguel de Cervantes Saavedra, Don Quixote*

## 1. Introduction

Arguably, homology represents one of the major developments of mathematics in the $20^{\text{th}}$ century. However, model theory and set theory are among the few fields of mathematics where homotopy theory has, essentially, never been applied. Indeed, with the exception of o-minimality, where homotopy/homology theories generalizing those arising in real geometry are used on a regular basis, we do not know of any applications of homotopy theory in either fields. In recent years, model theoretic questions arise to which, so it seems, homotopy theoretic tools should be applied.

The immediate motivation for carrying out the present work is a series of works by Zilber, Bays and the first author. These works are concerned with Zilber's program, launched in [Zil06], to apply model theoretic methods to the study of (non compact) complex analytic structures. In his D.Phil [Gav06], the first author, showed that Zilber's main (technical) result in [Zil06] can be naturally translated into familiar algebro-geometric terms. This translation also allowed Gavrilovich to generalize the statement of Zilber's result, casting them in the form of algebro-geometric conjectures. To prove some of these conjectures (modulo necessary corrections), Bays in his thesis ([Bay09], [BZ11]) translated them back into the language of model theory, where Shelah's machinery of *excellent classes* had to play a significant role. But in order to get Shelah's technology into work, non-trivial algebraic and algebraic-geometric information had to be obtained. To algebraic geometers and algebraists, information of the sort required in these proofs is usually given in the language of homology theory. Unfortunately, there is currently no dictionary translating modern homology theory into the language of model theory. Therefore, in order to be applied in this context, the algebraic tools first have to be "translated" into a more classical language. For Bays' work this translation was readily available, but it may turn out to be a serious obstacle for developing the theory further.

Motivated by these problems, the first author embarked on an attempt to develop at least some parts of the missing toolbox. As a baby version he started playing with the simplest, most naive, construction he could imagine of a homotopy theory akin to the one he was looking for: an almost degenerate homotopy theory for set theory. Quite surprisingly, playing around with this new toy, we - very naturally - reached, from a totally different angle, some set theoretic concepts playing a central role in Shelah's PCF theory: we recover the covering number of a cardinal $\lambda$ as the "correct" homotopy theoretic version of cardinality, and - by a slight variation - retrieve Shelah's revised power function. In the present note we describe the construction of a homotopy theory associated with set theory. In a subsequent paper [GH10] we will explain how to use this homotopy theory to recover, in a natural way (from the homotopy theory's point of view), ideas from Shelah's PCF theory.

The first author was partially supported by a MODNET (European Commission Research Training Network) grant and by the Skirball foundation as a post-doctoral fellow at Ben Gurion University.

The second author was partially supported by GIF grant No. 2266/2010 and by ISF grant No. 1156/10.

1.1. **A homotopy theoretical synopsis.** From the homotopy theoretic point of view the present note is, on the technical level, a triviality. To a reader with a basic familiarity with category theory and the first few pages of Quillen's book [Qui67] this synopsis should provide a fairly good idea of the contents of this note. Such readers may find it simpler to fill in the details themselves, rather than read our rendering of them. Readers less familiar with these concepts and ideas may find it more convenient to first read the background on model categories provided in the next section, before returning to this synopsis.

Imagine a simple minded homotopy theorist, or even better - a homotopy theory android - trying to understand the introductory, set-theoretic, chapter commonly preceding undergraduate texts in topology, say. What is the (simplest) category which could help our homotopy theorist understand at least some of the notions appearing in such a chapter. Sets are default candidates as objects. But what should the arrows be? The membership relation ($a \in A$) is not transitive, and therefore cannot serve as an arrow. Simple inclusions seem to be the next simplest idea. Knowing no better our homotopy android starts with this simple category, and tries to add information that will enable it to capture three basic notions of set theory: *finiteness*, *countability* and *equi-cardinality.*

This information is added following Quillen's formalism of model categories. In order to keep the construction as simple as possible, the android intends to construct a closed model category. Trying to construct a dictionary between Quillen's formalism and the unfamiliar set theoretic world, our homotopy theorist first tries to identify those arrows in our newly constructed category which could be naturally thought of as (weak) equivalences:
Set theoretically, adding a single element to a single set suggests itself as the least significant operation, and could therefore be declared a weak equivalence. As compositions of weak equivalences are also weak equivalences any arrow $A \longrightarrow B$ such that $B \setminus A$ is finite will be declared a weak equivalence. To keep track of the work done so far those arrows will be labelled (w). Of course, these arrows will not be the only weak equivalences in our category.

Now, our android - happy for dealing with the notion of finiteness - has to decide whether he should focus on defining fibrations or co-fibrations. He knows that on the category theoretic level he should have initial and terminal objects. Of those, the empty set, $\emptyset$ is indeed an initial object, whereas the terminal object, $V$ - the universe of set theory - is not a set, and has to be formally added to his category. This gives him reason to believe that co-fibrant objects will have a more natural set theoretic interpretation than fibrant objects. So he decides to focus his attention on co-fibrations. First, he declares that an arrow $A \longrightarrow B$ is a *cofibration* if $A$ and $B$ have the same cardinality. To keep track of this he labels such arrows $A \xrightarrow{(c)} B$. Checking the list of tasks to do, our hero notices that Axiom (M2) of model categories requires that - in particular - arrows of the form $A \xrightarrow{(w)} A \cup \{b\}$ (for $\{b\}$ a singleton) decompose as a (weak) cofibration followed by a fibration. For infinite $A$ this means (using induction) that for any finite set $\{\bar{b}\}$ he must declare $A \xrightarrow{(wc)} A \cup \{\bar{b}\}$. A little diagram chasing convinces him that, to avoid constructing a trivial model category the same must be true of any set $A$ (not necessarily infinite). Thus, till now $A \longrightarrow B$ would be labelled (c), if $\operatorname{card} A = \operatorname{card} B$ or both are finite.

So now he knows that if $A$ is finite $\emptyset \xrightarrow{(wc)} A$ (i.e. it is both a weak equivalence and a cofibration). He also knows that there must be sets $A$ such that $\emptyset \xrightarrow{(c)} A$ but not $\emptyset \xrightarrow{(wc)} A$. Such sets must be infinite. Since he still has to identify countability to achieve the task he set himself to do, he finds it natural to declare $\emptyset \xrightarrow{(c)} A$ if $A$ is countable. It is then easy for him to verify that $A \longrightarrow B$ has to be labelled (c) if either $\operatorname{card} A = \operatorname{card} B$ or $B$ is countable.

Here things become trickier, as our androids discovers that we are nearing a contradiction: On the one hand, Axiom (M2) assures that, up to weak equivalence, every arrow is a fibration, while - on the other hand - Axiom (M1) requires that any fibration has the left lifting property with respect to weak co-fibrations. But when considering arrows of the form $A \xrightarrow{(wc)} A \cup \{b\}$, it becomes obvious that our category simply does not have enough arrows (or, objects, which - since arrows between any two objects, when they exist, are unique - is an equivalent statement).

It now becomes clear that in order to solve this problem while keeping the category as simple as possible - i.e. arrows between objects are unique when they exist - we have to introduce as new objects in our category families of sets, but in order not to have two kinds of objects, and in order to keep the work already achieved, we identify a set $A$ with the singleton $\{A\}$. As a first approximation we allow all classes as objects. Viewing a class $\mathcal{A}$ as the direct limit of (the unique diagram of) all its members, there is no choice but to define $\mathcal{A} \longrightarrow \mathcal{B}$ if for all $A \in \mathcal{A}$ there exists $B \in \mathcal{B}$ such that $A \subseteq B$. This seems to be the simplest, and most natural definition our android could come up with compatible with everything considered up to this stage.

Now the process of producing a model category is almost automatic: First, take the structure (on our category) co-fibrantly generated by the two classes of arrows explicitly defined above. The resulting structure is not yet a model category - some obvious obstructions prevent our weak equivalences from satisfying the two-out-of-three axiom of model categories. Removing some objects and taking the resulting full subcategory we get an actual model category (this is, of course, not automatic, but simple enough for our android to verify).

The model category obtained in this way gives a homotopy theoretic interpretation to some basic set theoretic concepts. Most importantly for our homotopy android, a set $X$ (viewed as the *class* $\{X\}$) is countable if and only if $\emptyset \xrightarrow{(c)} \{X\}$, and if this arrow is also a weak equivalence then $X$ is finite; two infinite sets $A$ and $B$ have the same cardinality precisely when $\{A\} \xrightarrow{(c)} \{B\}$. Thus, the task it set for itself in the beginning is accomplished, he obtained a model category modelling the notions of finiteness, countability and equi-cardinality. But more interesting is the fact that the notion of a covering family of a set $A$ acquires a homotopical interpretation: $X$ is a covering family for $A$ if and only if $X \xrightarrow{(wf)} \{A\}$ is a weak equivalence and a fibration. It is now not hard to recover Shelah's covering numbers - a key notion in PCF theory. Though this is not very difficult, there are quite a few definitions to give, and some homotopy theory intuitions and technicalities to explain, so we postpone the discussion to a subsequent paper.

1.2. **The exposition.** There are three important guidelines to the exposition of this paper. The first is that homotopy theory is best written in homotopy theoretic language. Therefore the paper is written in the language of category theory. Combinatorial properties are, as a rule, transformed into diagrams and proofs are, quite often, translated into (simple) diagram chasing arguments. This is by no means an ideological choice. The standard set theoretic intuition is lost at the early stages of the construction, and we have to stick to Quillen's homotopy theory as a guide. Our choice of language allows us to keep track of this intuition.

The second guideline to the exposition is the realisation that potential readers (plural (!)) of this paper are, probably, set theorists, with little familiarity with category theory, and no familiarity with model categories. For that reason we spell out known definitions and our arguments are detailed enough for the paper to be - essentially - self contained.

The third and final guideline to the exposition is our intention to demonstrate that from a homotopy theoretic point of view our construction as well as its analysis are, essentially, automatic and could be achieved by an automaton akin (to the best of our understanding) to Gromov's Ergosystem ([Gro09]). In some respects the exposition of the present paper tries to address the following observation from [Gro09]:

> "The category/functor modulated structures can not be directly used by ergosystems, e.g., because the morphisms sets between even moderate objects are usually unlistable. But the ideas of the category theory show that there are certain (often non-obviuos) rules for generating proper concepts."

Our introduction of the *homotopy android* is a direct refernce to Gromov's text, and is used to to explicate, whenever possible, the algorithmic nature of our application of ideas coming from a (superficial) reading of Quillen's text on model categories. Indeed, our work could be considered as the exploration of a toy version of the ideas introduced in §6 of Gromov's paper. We do not claim that our construction captures any of the flavour of "real" homotopy theory, or that it would be considered natural by a seasoned homotopy theorist. Rather, we try to demonstrate that attempting to model some of the most basic notions of naive set theory using a purely syntactic analysis of Quillen's book [Qui67] one can construct a model category and that from the point of view of set theory, this model category may be of some interest. In the present paper we carry out the construction, demonstrating the first of these goals, in a second paper we will try to explain why set theorists might care about this construction.

The requirement that the construction and verification could be carried out by a "model category automaton" leads us (and see also our quotation from [Gro09]) to consider only categories where arrows between objects are unique (when they exist). Diagram chasing in such categories is considerably simpler than in a general category, and amounts to showing the existence of certain arrows (commutativity of diagrams being automatic). This allows us to keep all arguments (but one !) in this paper on a level accessible to our android. It is for that reason, for example, that when running into trouble with Quillen's Axiom (M2) the automaton chooses to add objects to his category, rather than consider the possibility of adding new arrows.

We start (Section 2) with a brief overview of what little we need from category theory and Quillen's homotopy theory. In our presentation we stick to the strict minimum required for this note to be self contained. We then proceed to define (Definition 3) a category, $StNaamen$, whose objects are arbitrary classes (of sets). The arrows of our category give it the structure of a quasi partially ordered class. We then label the arrows of $StNaamen$ (c),(f),(w) (Definition 4) in a first attempt to get a model category in the sense of Quillen - but fail; to overcome this failure we then define (Definition 1) a full subcategory, $QtNaamen$, of $StNaamen$ such that the induced labelling on $QtNaamen$ does give it the structure of a model category (Theorem 48, the main result of this paper). Proposition 22 gives an alternative, combinatorial, definition of our labelling. This characterisation is useful for some of the arguments and may also be helpful for readers not feeling as comfortable chasing commutative diagrams as our imaginary homotopy theorist.

**A set theoretic disclaimer.** Our use of set theory is naive, and we intentionally ignore set theoretic questions naturally arising in the context of the "category of sets". Since in the categories we are dealing with arrows are unique (when they exist) the only problem that may arise concerns the nature of the objects in our category. Any standard solution of such problems would address all questions of this nature which may arise in the construction described in this paper. For example, we could declare our objects to be classes (definable with parameters), and it is easy to describe a simple formalism in which our construction would fit. A simpler approach would be to fix a strongly inaccessible cardinal $\kappa$ and identify the objects of our category with $\mathbb{P}(V_\kappa)$. In such a setting the collection of objects of our category (the collection $\mathcal{Ob}$) can be identified with a subset of $\mathbb{P}(\mathbb{P}(V_\kappa))$. The arrows in our category ($\mathcal{Mor}$) can then be identified with a definable subset of $\mathcal{Ob} \times \mathcal{Ob}$, and the labelling associated with Quillen's model categories can be thought of as (definable) unary predicates on the set $\mathcal{Mor}$.

Having said that, we will from now on ignore all set theoretic questions of this nature, with the conviction that readers concerned with the possibility of set theoretic paradoxes arising as part of the construction can easily fill in all the details in either of above solutions, or any other standard solution they may find more attractive.

## 2. Model categories

Though we could give a set theoretic interpretation of our construction, we feel that such an exposition would entirely miss the point of our paper: to homotopy theorists it would be incomprehensible, and to set theorists it would seem arbitrary and unnatural. For that reason, and since the whole purpose of our work is to show how non-trivial set theoretic concepts arise naturally from a homotopy theoretic approach to the subject, our exposition is given wholly in the language of Quillen's model categories. On the technical level, our construction is a triviality and, being almost degenerate as a model category, the homotopy theoretic tools needed to comprehend it are very basic. In the present section we give a brief exposition of these tools. We focus more on ideas and intuitions that will guide us in our usage of these tools than on technical details. We first remind the reader of the basic notions of category theory, and in the subsequent subsection we proceed to describe model categories.

2.1. **Categories**[1]**.** Combinatorially, a category is a directed graph equipped with a notion of *composition* of some arrows, such that for every subgraph $X \xrightarrow{f} Y \xrightarrow{g} Z$, there is a *unique* arrow $h : X \longrightarrow Z$ called the *composition* of $f$ with $g$, denoted $h = f \circ g$ or simply $h = fg$. We allow multiple edges between two vertices as well as loops - edges leaving and entering the same vertex - called self-morphisms.

The philosophy underlying category theory is that, when viewed in the right category, the internal structural properties of any mathematical object are fully reflected in the combinatorics of the graph (or rather the 2-polyhedron) of morphisms-arrows around it. Rather amazingly, this language, if properly (often non-obviously) developed, allows a concise uniform description of mathematical structures in a vast variety of cases.[2]

An important axiom of category theory requires that every object, $\bullet$, in the category comes with a distinguished self-morphism into itself, called the *identity*, denoted usually $\mathrm{id}_\bullet$, such that $f \circ \mathrm{id}_\bullet = f$ and $\mathrm{id}_\bullet \circ f = f$ whenever the composition is defined. Thus, we can define a morphism $\bullet_1 \xrightarrow{f} \bullet_2$ to be an *isomorphism* if there exists a morphism $\bullet_2 \xrightarrow{g} \bullet_1$ such that $f \circ g = \mathrm{id}_{\bullet_1}$ and $g \circ f = \mathrm{id}_{\bullet_2}$. In the present note we will be interested in large categories, i.e. categories whose objects are classes (e.g. the class of all sets). Other than that our category will be extremely simple, as it will turn out that between any two objects there exists at most one morphism. In particular, in this context, a morphism $\bullet_1 \longrightarrow \bullet_2$ is an isomorphism precisely when there exists an arrow in the opposite direction.

The skeleton of a category is obtained by contracting to a point all arrows which are isomorphisms in the category. The objects of the resulting structure are isomorphism classes of objects of the original category. But the resulting arrows no longer give it the structure of a category: as a rule, morphisms cannot be *uniquely* composed unless we identify all arrows between isomorphic objects with the identiy morphism in the skeleton, sometimes causing more information loss than desired.

[1]This subsection borrows much from §3.5 of [Gro09]; in particular we cite a few sentences almost verbatim. We apologise for not daring to bother the author asking for permission to quote his work: our justification is that citing an in-text is, for all practical purposes, the same as inserting a link to the files which are available online.

[2]and Gromov adds: Some mathematicians believe that no branch of mathematics can claim maturity before it is set in a category theoretic or similar framework and some bitterly resent this idea.

Despite this deficiency, the skeleton encapsulates much of the significant mathematical information carried by the original category. By this we mean that most mathematical theorems and definitions concerning, say, rings, fields, Lie groups etc. are *isomorphism-invariant*, although many of the constructions involved are most certainly not (e.g. the value $f(x)$ of a function $f$ does not make sense when passing to morphisms between isomorphism classes of objects). In Quillen's words (originally applied to model categoties), such theorems or definitions concern the *skeleton* of our original category "with all extra structure which comes by performing constructions in the original category" with the collection of commutative diagrams.

Finally, it is hard to discuss category theory, even on the basic form we need it, without introducing the notion of *functors.* Those are simply functions between categories respecting arrows and the only operation such objects have - the composition of arrows. More precisely, if $\mathfrak{C} := (\mathcal{Ob}, \mathcal{Mor})$ and $\mathfrak{C}' := (\mathcal{Ob}', \mathcal{Mor}')$ are categories, a functor is a function $\gamma : \mathfrak{C} \to \mathfrak{C}'$ mapping the objects in $\mathfrak{C}$ to objects in $\mathfrak{C}'$ and the arrows in $\mathfrak{C}$ to arrows in $\mathfrak{C}'$ respecting the composition of morphisms.

2.2. **Commutative diagrams.** A basic notion in category theory is that of a commutative diagram. A diagram in a category $\mathfrak{C}$ is a (usually finite) subgraph $D$ of $\mathfrak{C}$, i.e. a set of objects and morphisms between them. A diagram with no loops is commutative if for every two objects $\bullet_1$ and $\bullet_2$ in $D$ and every two $D$-chains of composable morphisms, both starting at $\bullet_1$ and terminating at $\bullet_2$, the compositions of the morphisms in the two chains (which are certain morphisms from $\bullet_1$ to $\bullet_2$) are *equal.* Because in the categories we will be interested in arrows are unique, all diagrams are commutative. This is the key to the fact that most arguments in the present paper are essentially trivial.

In the present paper commutative diagrams play a crucial role in the exposition. We use them to convey in a concise and intuitive way the axioms of model categories, definitions and proofs. In order to incorporate into the standard schematic representation of commutative diagrams the additional information needed to express definitions (e.g. quantifiers, order of quantification etc.) and proofs (e.g. the order in which arguments appear and the justification for each argument) etc. we use different arrow styles and labelling. These will be introduced in more detail in Sub-section 2.3.1, and throughout the text.

2.3. **Quillen's Model category and homotopy theory. Commutative diagrams as $\forall\exists$-formulae.** As we have noted mathematicians tend (often implicitly) not to be interested in properties which are not isomorphism invariant. To a homotopy theorist the same is true for a wider class of functions: homotopy equivalences. For our purposes, Quillen's formalism of model categories can be viewed as a tool for extracting from a category (with some extra structure - a labelling of arrows) the "homotopy skeleton" in a way not dissimilar to what was described above. As Quillen points out, this formalism of labelled categories satisfying certain simple axioms has "sufficient generality to cover in a uniform way [...] different homotopy theories [allowing to carry out] a large number of arguments [...] formally similar to ones in algebraic topology." In a nutshell a syntactic rendering of Quillen's approach can be described as follows:

Take a category, and divide its morphisms into several (non-disjoint) families: simple morphisms (unlabelled), co-fibrations (labelled (c)), fibrations (labelled (f)) and weak homotopy equivalences (labelled (w)[3]). Make sure the labelling is rich enough (so that, e.g., every simple morphism is the composition of labelled ones and all isomorphisms are labelled (wcf)) and satisfies certain simple axioms (to be discussed below). "Localize" the resulting "coloured" category $\mathfrak{C}$ at (w) by formally inverting all arrows labelled (w) to obtain the *homotopy category* $S^{-1}\mathfrak{C}$. As pointed out by Quillen ([Qui67], §1.1, Definition 5), "except for set-theoretic difficulties, the category $S^{-1}\mathfrak{C}$ exists and may be constructed by mimicking the construction of a free group...", i.e. for each arrow $\mathcal{X} \xrightarrow{(w)} \mathcal{Y}$ add an arrow $\mathcal{X} \xleftarrow{(w)} \mathcal{Y}$ such that the compositions in both ways are $\mathrm{id}_{\mathcal{X}}$ and $\mathrm{id}_{\mathcal{Y}}$, respectively. Close the resulting class of arrows under compositions to obtain a category, $\mathrm{Ho}\mathfrak{C}$, the *homotopy category associated with* $\mathfrak{C}$.

Though, as mentioned above, $\mathrm{Ho}\mathfrak{C}$ is a category, it is usually considered with some extra structures. As with the passage from a category to its skeleton, some constructions, though homotopy invariant, have to be performed in the original category $\mathfrak{C}$, before passing to $\mathrm{Ho}\mathfrak{C}$ (see, e.g., [Qui67], §1.2). Thus, from Quillen's point of view "*the homotopy theory associated to a model category* $\mathfrak{C}$ *is defined to be the homotopy category* $\mathrm{Ho}\mathfrak{C}$ *with the extra structure [...] coming from performing constructions in* $\mathfrak{C}$" ([Qui67] §1.0 and §1.2). In the present context, however, most of this extra structure trivialises.

**Remark 1.** It is not hard to see that the homotopy category $\mathrm{Ho}\mathfrak{C}$ associated with a category $\mathfrak{C}$ depends only on the class of weak equivalences in $\mathfrak{C}$. The reader may wonder, therefore, why should we bother with fibrations and co-fibrations, if the extra structure induced on Ho$\mathfrak{C}$ from $\mathfrak{C}$ trivialises in our construction. The answer is twofold. First, not all the extra structure trivialises, and some of it will be of use to us. More importantly, as we will see shortly, we have to specify what are the fibrations and co-fibrations in our category in order to verify that our weak equivalences satify Quillen's axioms, i.e. that they are indeed the weak equivalences in *some* model category.

[3] This notation is adapted from [Dur07]; we use (wf), (wc) instead of his (af), (ac).

We will now give a brief overview of Quillen's axiomatisation of model categories. The topological (or homological) intuitions underlying this axiomatisation are not important for our purposes so we will not try to explain them.

**Axiom** ($M0$)**.** The category $\mathfrak{C}$ is closed under finite projective and injective limits.

Graphically, this axiom can be described by the following commutative diagram:

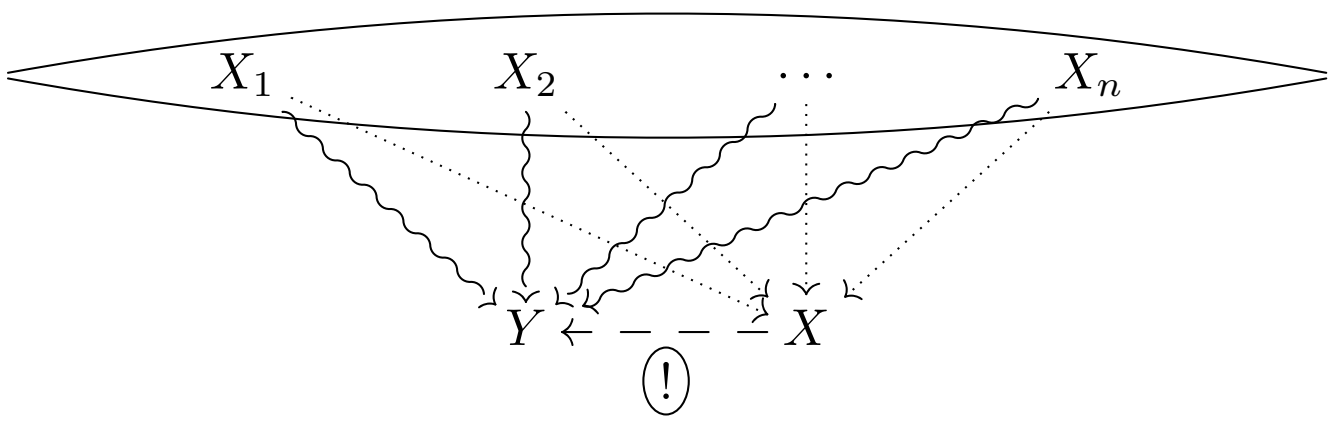

FIGURE 1. Inductive limits. The diagram reads (see, also, Remark 2 below): For every diagram with vertices $X_1, \ldots, X_n$ (the arrows in the diagram are not specified) there exists a system of dotted arrows, such that for any system of arrows ($\rightsquigarrow$) there exists a unique dashed arrow ($\dashrightarrow$).

The corresponding diagram for injective limits is obtained by inverting all arrows in the above diagram. Axiom (M0) is, essentially, a technical one, and deals solely with the pure category structure of $\mathfrak{C}$.

**Remark 2.** To be more precise, Figure 1 should be described as follows: "for every *commutative* diagram with vertices $X_1, \ldots, X_n$ there exists a system of dotted arrows *making the resulting diagram commute*, and such that for any system of arrows ($\rightsquigarrow$) *such that the resulting diagram commutes* there exists a unique dashed arrow ($\dashrightarrow$) *making the whole diagram commute.*". In order to keep things simpler and more comprehensible, we will assume from now on that all diagrams are drawn commutative, and given a diagram, whenever we refer to new arrows that can be added to it, we will mean "added such that the resulting diagram is commutative".

**Notation 3** (Commutative diagrams I)**.** Throughout this note *all diagrams we draw are assumed to be commutative.* Commutative diagrams will be used systematically throughout this note. Among others, diagrams will be used to introduce the axioms of model categories. The axioms we will describe using diagrams will be $\forall\exists$-statements (possibly, as in Axiom (M0), of deeper quantifier depth). In such cases we will keep the notation of Figure 1. Namely, solid arrows ($\longrightarrow$, suppressed in the previous figure) will denote those arrows associated with the first universal quantifier; dotted arrows ( $\cdots\!\!>$ ) will denote the arrows associated with the existential quantifier of depth 1, then ($\rightsquigarrow$) and ($\dashrightarrow$) for additional alternating quantifiers.

The next axiom lies at the heart of our construction:

2.3.1. *The lifting property.* Quillen identified the *lifting property* as a concise way for describing several notions and tools in algebraic topology. Moreover, it turned out that quite often various proofs could be reduced, by diagram chasing, to applying the lifting property to various arrows.

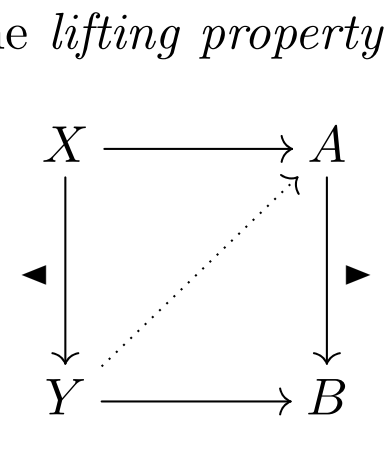

FIGURE 2. ( ⋌ ) The diagram reads: for every pair of horizontal arrows $X \longrightarrow A$ and $Y \longrightarrow B$ there exists a *lifting* arrow in the diagonal direction $Y \longrightarrow A$, or, equivalently, the arrow $Y \longrightarrow B$ *lifts to* $Y \longrightarrow A$ (recall the standing assumption that all diagrams commute).

**Definition 1.** An arrow $A \longrightarrow B$ has the (right) lifting property with respect to the arrow $X \longrightarrow Y$, denoted $X \longrightarrow Y \rightthreetimes A \longrightarrow B$, if Figure 2 is true for the arrows $A \longrightarrow B$ and $X \longrightarrow Y$.

More generally,

**Notation 4.** Given an arrow $f$ and a collection of arrows $F$ we let $f \rightthreetimes F$ denote $f \rightthreetimes g$ for all $g \in F$. Similarly for $F \rightthreetimes f$.

**Notation 5** (Commutative diagrams II)**.** Commutative diagrams will also be used to introduce new definitions. Identifying an object $\mathcal{X}$ with either $\emptyset \longrightarrow \mathcal{X}$ or $\mathcal{X} \longrightarrow \top$ (for the terminal object $\top$, see Claim 8) the properties defined will always be properties of arrows. We keep the notation introduced above. However, in diagrams encoding a definition we should be able to distinguish the object(s) being defined (the free variables of the formula being encoded). Those arrows we will denote by ◄ or ►.

Below we give a few examples of basic mathematical notions which can be defined through an appropriate lifting property:

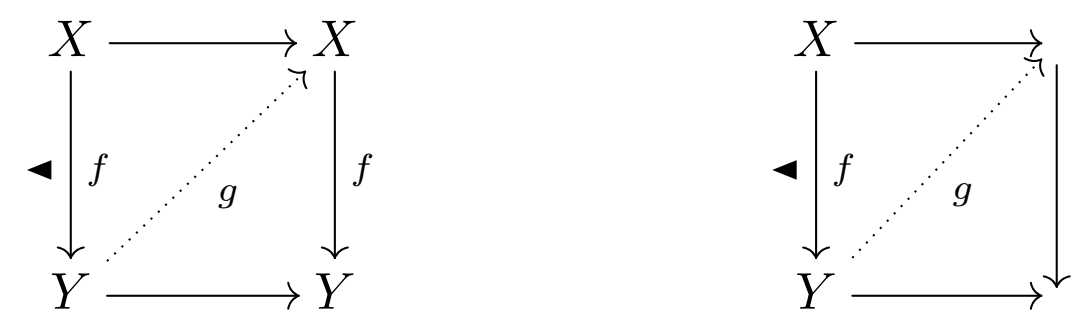

FIGURE 3. (Isomorphism) The arrow $f : X \longrightarrow Y$ is an isomorphism if and only if it has the (left) lifting property with respect to itself (equivalently, with respect to any arrow). For, taking $\mathrm{id}_X : X \to X$ and $\mathrm{id}_Y : Y \to Y$ as the horizontal arrows, $f \circ g = \mathrm{id}_X$ and $g \circ f = \mathrm{id}_Y$.

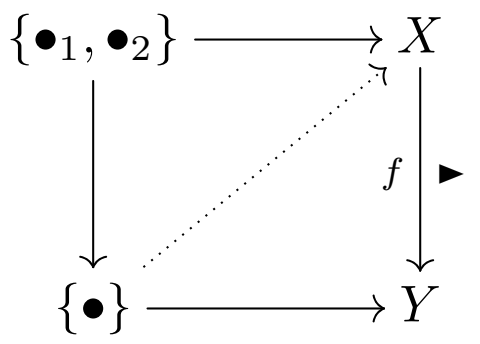

FIGURE 4. (Injectivity) The arrow $f : X \longrightarrow Y$ is injective if and only if it has the (right) lifting property with respect to any arrow of the form $\{\bullet_1, \bullet_2\} \to \{\bullet\}$ (in any category where such arrows make sense). For, if $f(x_1) = f(x_2)$ for some $x_1, x_2 \in X$ no commutative diagram as above has the lifting property, if $(\bullet_1, \bullet_2) \to (x_1, x_2)$ is taken for the upper horizontal arrow.

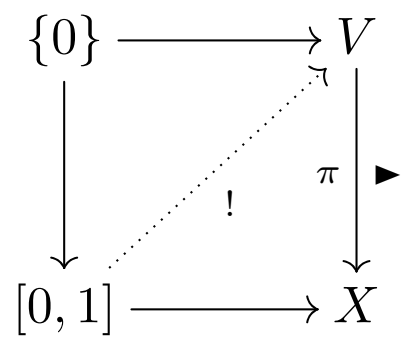

FIGURE 5. (Covering space) In the category of topological spaces, the pair $(V, \pi)$ is a covering space of $X$ if the lifting property described in the diagram holds.

In a similar way one can define surjective functions, injective and projective modules etc. Now that our reader, hopefully, has some intuition with respect to the lifting property, we can define Quillen's most important axiom (from our perspective).

**Axiom** (M1)**.** The following two lifting properties for labelled arrows hold:

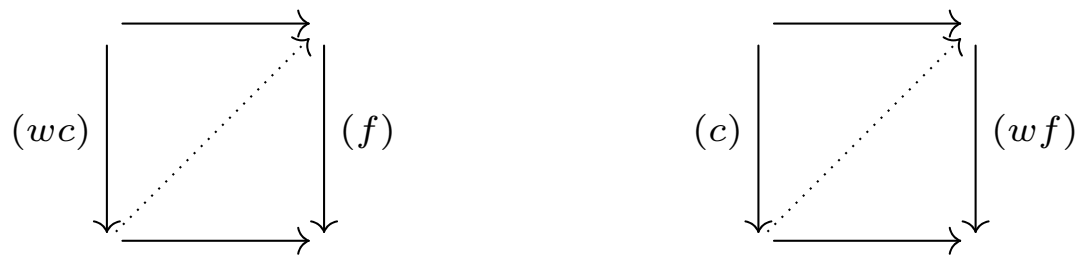

In our construction, once we identify those arrows we *require* to be weak equivalences and co-fibrations, these axioms will, essentially, determine the labelling of all other arrows. For that reason we consider Axiom (M1) as the key to the construction.

2.3.2. *The remaining axioms.* We now proceed to describe the remaining axioms in Quillen's formalism. The role played by these axioms in our construction is not as important as the role played by Axiom (M1). For that reason we will be brief. The next axiom assures that there are enough labelled arrows:

**Axiom** (M2)**.** The following two ∀∃-diagrams hold:

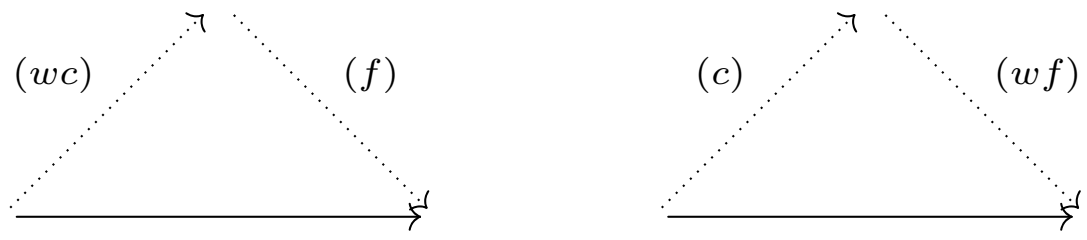

The next axiom assures that our notions of fibrations and co-fibrations are stable enough:

**Axiom** (M3(a))**.** Fibrations and co-fibrations are stable under compositions. Namely, the following two ∀∃-diagrams hold:

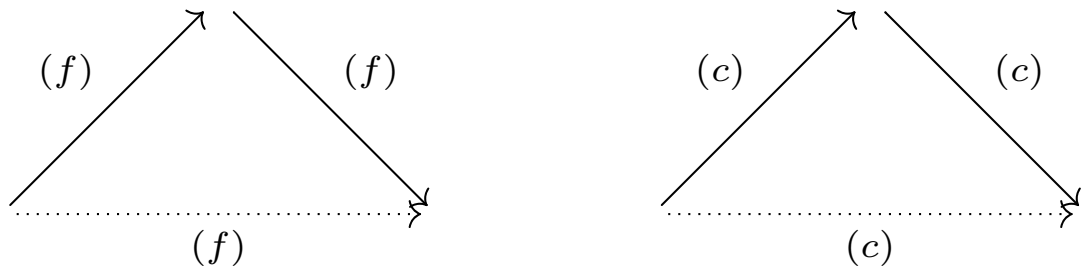

**Axiom** (M3(b))**.** Isomorphisms are fibrations, co-fibrations and weak equivalences:

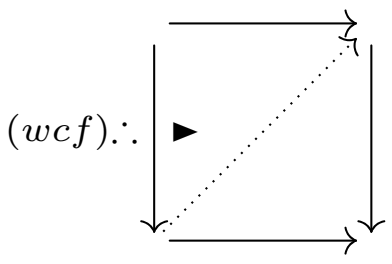

FIGURE 6. The figure reads: if the commutative ∀∃-diagram is true *then* the left arrow is labelled (wcf).

**Notation 6** (Commutative diagrams III)**.** As in Figure 6, we will use the notation $X \xrightarrow{\therefore(\cdot)} Y$ to mean "if the commutative diagram is true, then $X \longrightarrow Y$ is labelled $(\cdot)$".

**Axiom** (M3(c))**.** Fibrations and co-fibrations are stable under base change and co-base change respectively. I.e. the following diagrams are true (see Figure 1 for notation):

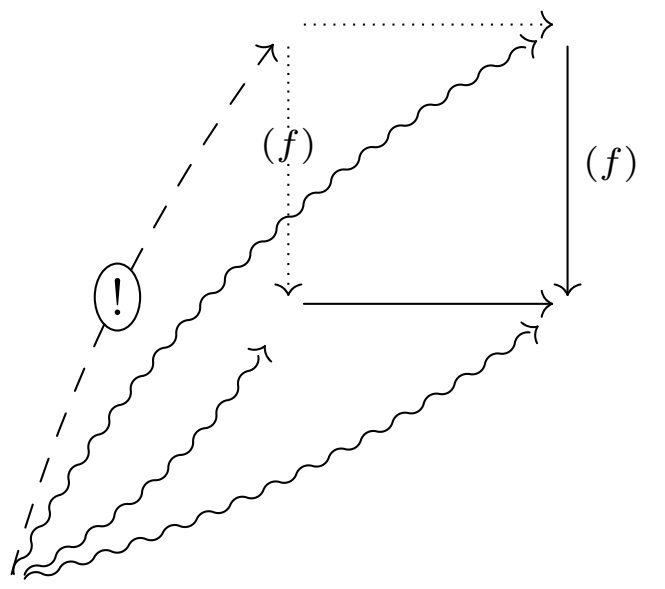

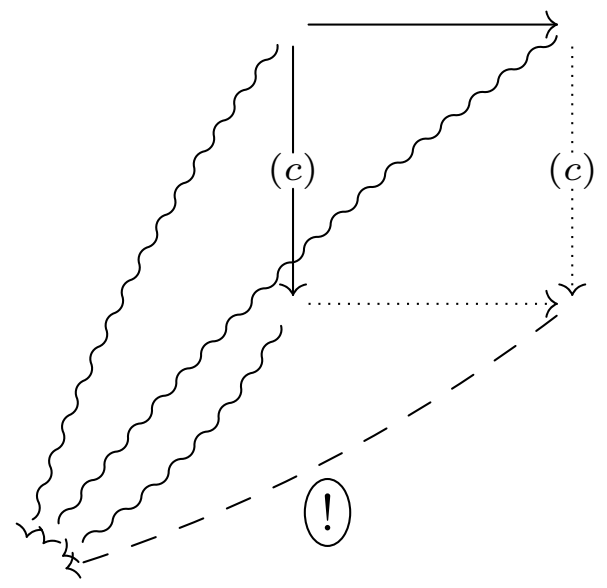

The next axiom extends Axiom M3(c) to weak equivalences, in this case a little more care is needed:

**Axiom** (M4)**.** The base extension of an arrow labelled (wc) and the co-base extension of an arrow labelled (wf) are both labelled (w) (again, see Figure 1 for the notation):

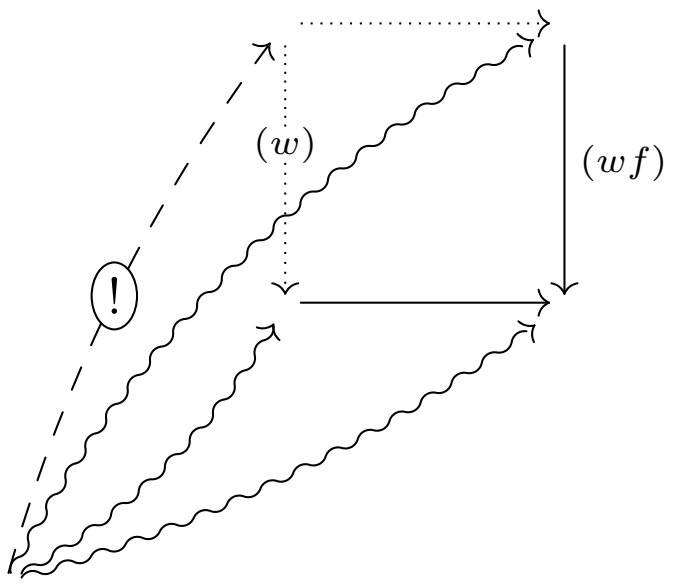

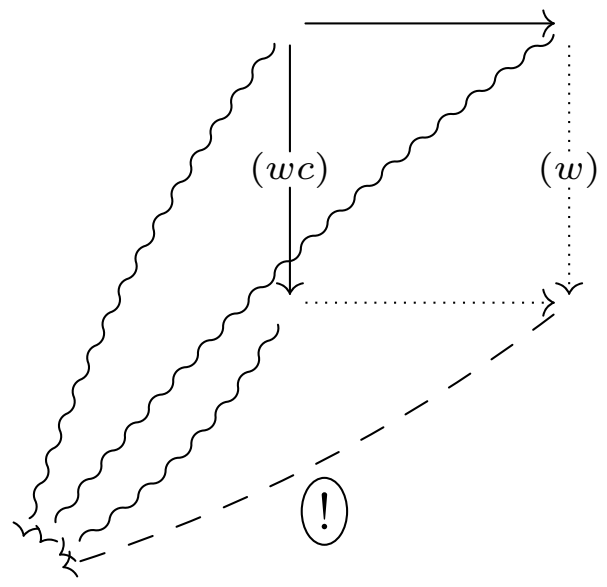

The last axiom assures that weak equivalence is close enough to being transitive:

**Axiom** (M5, Two out of three)**.** In a triangular diagram, if any two of the arrows are labelled (w) so is the third

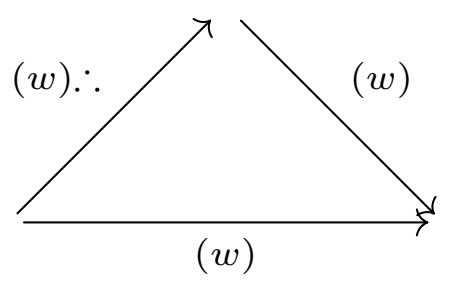

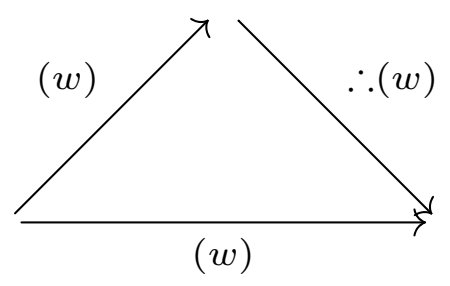

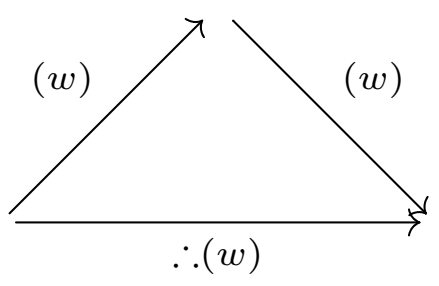

Axiom (M5) is a formal rendering of the informal ideology that (w)-arrows should be thought of as some sort of equivalence. The primary intention of this axiom is to ensure that Ho$\mathfrak{C}$ is a category for any model category $\mathfrak{C}$, i.e., that after localizing the category $\mathfrak{C}$ at the set of (w)-arrows we obtain a category.

Now that all the axioms of a model category have been spelled out, we observe that - somewhat trivially - every category $\mathfrak{C}$ can be labelled in a way making it a category: label every isomorphism (wcf), and (cf) all remaining arrows (observe that because of Axiom (M1) labelling, even more trivially, all arrows (wcf), will not necessarily work). For this labelling, the homotopy category Ho$\mathfrak{C}$ of $\mathfrak{C}$ is $\mathfrak{C}$ itself and the associated homotopy skeleton is precisely the skeleton of the category $\mathfrak{C}$. To be meaningful, our construction of the next section will have to be a little more careful.

**Remark 7.** A model category $\mathfrak{C}$ is called *closed* if any two of the following classes of arrows in $\mathfrak{C}$ - fibrations, co-fibrations and weak equivalences - determine the third by the rules:

(1) An arrow is labelled (f) if and only if it has the right lifting property with respect to all arrows labelled (wc).
(2) An arrow is labelled (c) if and only if it has the left lifting property with respect to all arrows labelled (wf).
(3) An arrow $g$ is labelled (w) if and only if $g = u \circ v$ where $v$ is an arrow having the left lifting property with respect to all arrows labelled (f) and $u$ has the right lifting property with respect to all arrows labelled (c).

We will later on refer to this additional axiom as (M6). Our construction of a model category for set theory will make it almost automatically a closed model category. This will become handy in the construction because once certain arrows have been labelled (c), (f) or (w) the axioms of a closed model category will guide us through the labelling of the remaining arrows.

Here are some simple facts concerning closed model categories:

**Claim 8.**

1. A model category has initial and terminal objects namely objects $\emptyset$ and $\top$ such that $\emptyset \longrightarrow \mathcal{C} \longrightarrow \top$ for every object $\mathcal{C}$.
2. A non-degenerate model category (i.e. a model category $\mathfrak{C}$ whose homotopy skeleton is not a singleton with a unique arrow) has non-trivial co-fibrant objects, namely there exists an object $\mathcal{C}$, not isomorphic to $\emptyset$ such that $\emptyset \xrightarrow{(c)} \mathcal{C}$.
3. Moreover, in a non-degenerate model category there exists non-trivial co-fibrant objects, namely there exists an object $\mathcal{C}$ such that $\emptyset \xrightarrow{(c)} \mathcal{C}$ is not a weak equivalence.
4. If all diagrams in the model category commute, then if $\mathcal{C} \longrightarrow \mathcal{D}$, $\mathcal{C} \xrightarrow{(c)} \mathcal{E}$ and $\mathcal{D} \longrightarrow \mathcal{E}$ then $\mathcal{D} \xrightarrow{(c)} \mathcal{E}$.
5. The composition of two (wc)-arrows is a (wc)-arrow.

*Proof.* The first item follows automatically from Axiom (M0): the initial and terminal objects of the category are the inverse and direct limits (respectively) of the empty diagram. For the second item observe that if $\mathcal{C} \longrightarrow \mathcal{D}$ for some objects $\mathcal{C}, \mathcal{D}$ and there are no non-trivial co-fibrant objects then by Axiom (M2) it must be that $\emptyset \xrightarrow{(wf)} \mathcal{C}$ and $\emptyset \xrightarrow{(wf)} \mathcal{D}$. By Axiom (M5) it follows that $\mathcal{C} \xrightarrow{(w)} \mathcal{D}$. This implies that all arrows in our model category are weak equivalences. This means that the associated homotopy category is trivial. A similar proof gives the third point as well.

For the fourth clause we only need to observe that:

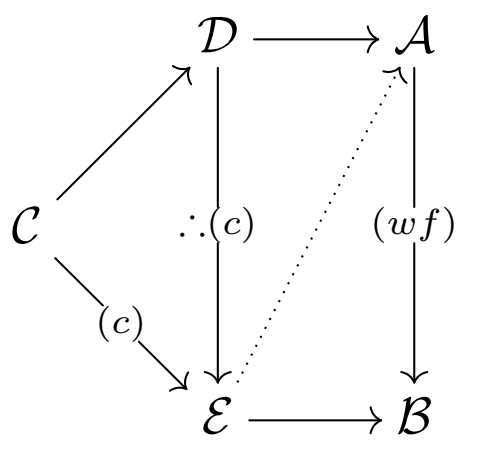

FIGURE 7. The dotted arrow exists by $\mathcal{C} \xrightarrow{(c)} \mathcal{E} \curlywedge \mathcal{A} \xrightarrow{(wf)} \mathcal{B}$. By assumption we know that the square commutes. By Axiom (M6) it follows that $\mathcal{D} \xrightarrow{(c)} \mathcal{E}$.

Finally, using Axiom (M6) the last clause reduces to simple diagram chasing:

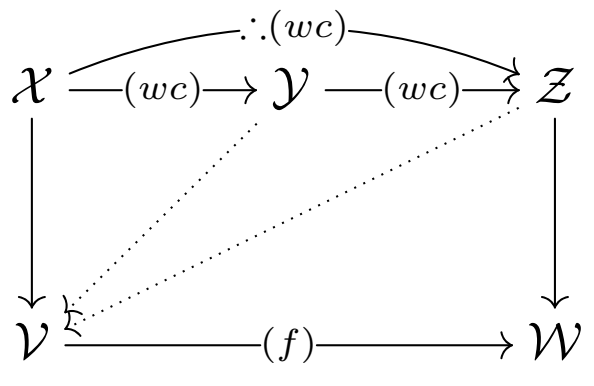

FIGURE 8. The arrow $\mathcal{Y} \longrightarrow \mathcal{V}$ is the lifting of $\mathcal{X} \xrightarrow{(wc)} \mathcal{Y} \curlywedge \mathcal{V} \xrightarrow{(f)} \mathcal{W}$. Now we know that the arrow $\mathcal{Z} \longrightarrow \mathcal{V}$ exists as the lifting $\mathcal{Y} \xrightarrow{(wc)} \mathcal{Z} \curlywedge \mathcal{V} \xrightarrow{(f)} \mathcal{W}$.

□

Readers not familiar with standard homotopy theoretic language may find it more convenient if we remind that:

**Definition 2.** Given a model category, $\mathfrak{C}$, an object $\mathcal{C}$ of $\mathfrak{C}$ is *co-fibrant* if $\emptyset \xrightarrow{(c)} \mathcal{C}$. A co-fibrant object is *trivial* if $\emptyset \xrightarrow{(wc)} \mathcal{C}$.

## 3. Our fairy tale

This section describes how our hero, who (unlike the authors of this paper) is a homotopy theorist (albeit an automated one) and who (like the authors of this paper) knows very little (if any) set theory could try to use his field of expertise to make some sense of the alien world of set theory. In the sequel we will describe how, in fact, our automaton could have discovered some of PCF theory (had Shelah not already discovered it). Whether our hero could have managed to do anything with his discovery remains unclear.

Our hero is not a seasoned homotopy theorist and has never seen homotopy theory in action. For him, Quillen's axioms of a model category do not carry any topological or homological intuition. Rather, they are procedures he can use to add new arrows to existing diagrams and figuring our their correct labelling. Faced with a problem our hero would not know "what would a homotopy theorist do?", rather, he will try to translate the problem into the language of homotopy theory, and look for clues in a textbook, e.g. [Qui67, Ch.I], he has mastered with great success. Our hero's lack of any knowledge of set theory (and having to start his exploration at some point), is guided by his love for simplicity and syntactic analysis. Or perhaps our hero is an ergosystem device (see [Gro09]) performing a *syntactic structure* analysis of a set theory book along the lines of §6 of [Gro09] with hard-wired *model category* diagram chasing.

**Remark 9.** Obviously, axioms (M1)-(M3) and (M5), (M6) fit into our description of the android's way of applying Quillen's machinery to produce "proofs". We should, however, explain how this works for axioms (M0) and (M4). For simplicity, we will focus on the part of Axiom (M0) stating that a direct limit exists. Diagramatically the existence of direct limits can be described as follows:

(1) Given a (finite) commutative diagram with vertices $\{X_i\}_{i\in I}$ add a new vertex $\mathcal{D}$ and arrows $\{X_i \longrightarrow \mathcal{D}\}$, making the whole diagram commute. Mark the new vertex and arrows by a special symbol $\circledless$.
(2) Given a (finite) commutative diagram with vertices $\{X_i\}_{i\in I}$, $\mathcal{D}$ and $D$, such that the diagram contains the arrows $\{X_i \longrightarrow \mathcal{D}\}_{i\in I}$ marked $\circledless$ and the (ordinary) arrows $\{X_i \longrightarrow D\}$, construct an arrow $\mathcal{D} \longrightarrow D$ making the diagram commute.

Note that applying (2) above in the case that $D \longrightarrow \mathcal{D}$ is an arrow in the diagram, we obtain the uniqueness of the direct limit.

His algorithm is, essentially, greedy. At every step he invents a simple task (indeed, the simplest he can think of) *visibly* bringing him closer to the goal of having a model category, e.g. by satisfying yet another axiom, sometimes looking back, correcting and readjusting his advancement as he goes. Essentially all arguments our hero has to verify amount to proving the existence of certain arrows. Though this is not the way we describe it, in order to do that, our robot can simply start adding arrows according to some pre-ordained algorithm (taking limits, using lifting properties of labelled arrows etc.) until he finds the desired arrow. Since no argument in the present paper splits into cases, our hero is very successful with this process. To shorten the search time of such an algorithm we allow our robot to use some simple heuristics: if the argument he has to prove concerns countable sets he start the construction of arrows by considering processes associated primarilly with countable sets (co-fibrant objects, in our hero's terminology); if he knows an argument not to be true in general, he first tries to construct arrows using the specific assumptions etc. Also, we assume that our automaton does not search proofs forever. He can look for proofs of up to a certain length. If he cannot find such proofs, he adds simplifying assumption, trying to look for special cases or possible counter examples.

Besides his diagram chasing capabilities he can apply basic inductive arguments (as long as they refer to diagram chasing and similar concepts), and he does have some understanding (that we do not explain precisely) of the notion of finiteness. He learns, and knows how to incorporate what he learns into his proof-searching algorithm (e.g., if he learns that certain arrows have more lifting properties than he originally was aware of, he will use this new knowledge later on, when trying to construct new arrows). He tends to believe in universal axioms: when faced with several examples sharing a common feature our android will believe that this common feature is universally true, unless there are obvious reasons suggesting otherwise. Finally, since closed model categories are easier to construct (one has to identify only two out of the three possible labellings), he tries to construct such a model category, though he does not really care whether he ends up with one or not, as long as the resulting labelling gives a model category.

3.1. **An alien world.** The first task our hero is facing is that of choosing the category that will best fit the task at hand. Sets are obvious candidates for objects. Arrows cannot designate the membership relation, as it is not transitive. Opting for the second simplest possibility, arrows are declared to be inclusion. To check himself he verifies that basic set theoretic operations can be expressed in this category, for example $Z = X \cup Y$ if and only if the following diagram is true

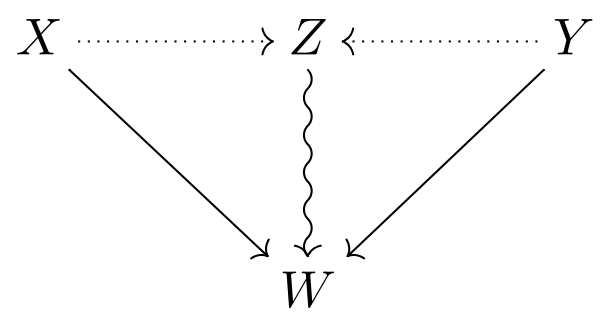

FIGURE 9. $Z$ is the union of $X$ and $Y$, the direct limit of the diagram whose vertices are $X, Y$.

and intersections can be defined by inverting all the arrows in this diagram. Similarly, stating that $A = \emptyset$ amounts to the diagram

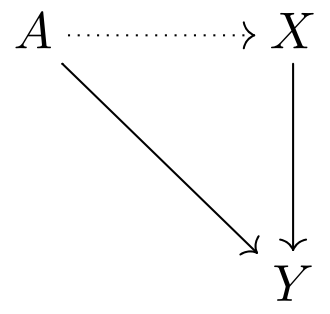

being true. Our hero, who is an expert in diagram chasing - and therefore not comfortable with definitions using negation - is quite happy with the definition he has found of the empty set. Confident, he now proceeds to tackle another basic definition using negation - the concept of the complement of a set (relative to a larger one). He finds that this, too, can be expressed in a simple diagram:

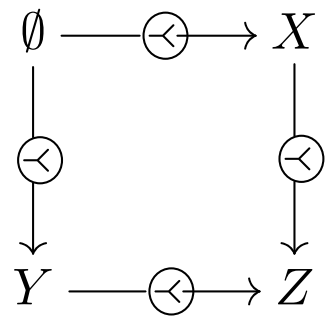

FIGURE 10. The set $Y$ is the complement of the set $X$ relative to the set $Z$ if $\emptyset$ is the inverse limit of $X$ and $Y$ and $Z$ is their direct limit.

He can now understand the notion of $A \setminus B$ as the complement (relative to $A$) of $A \cap B$. He also finds another definition of the complement, which he finds more appealing: if $X \longrightarrow Y$ then $Z$ is the complement of $X$ in $Y$ if the following diagram is true:

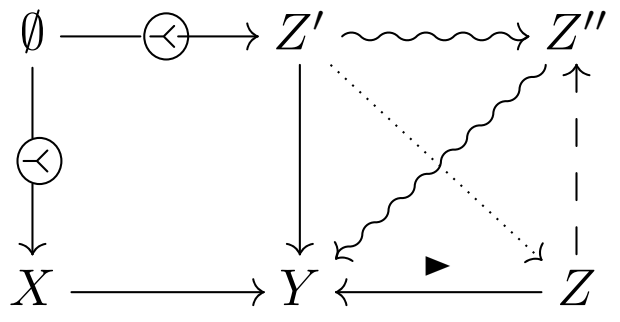

FIGURE 11. The complement of $X$ in $Y$ is the direct limit, $Z$, of all those objects $Z'$ such that the inverse limit of $Z'$ and $X$ is $\emptyset$. The diagram reads: for every $Z'$ such that $\emptyset = Z' \cap X$ there exists an arrow $Z' \longrightarrow Z$ such that any system $Z' \longrightarrow Z'' \longrightarrow Y$ factors through $Z$, $Z' \longrightarrow Z \longrightarrow Z'' \longrightarrow Y$.

He also notes that, the category being a lattice, the above observations imply that $\emptyset$ is the *initial object* of the category. Our hero decides that the model category he is about to construct should capture at least three basic notions of set theory: finiteness, countability and equi-cardinality. Before embarking on this adventure our hero notices, that regardless of how he chooses to label his arrows, (wcf)-arrows must be isomorphisms, as shows the following diagram:

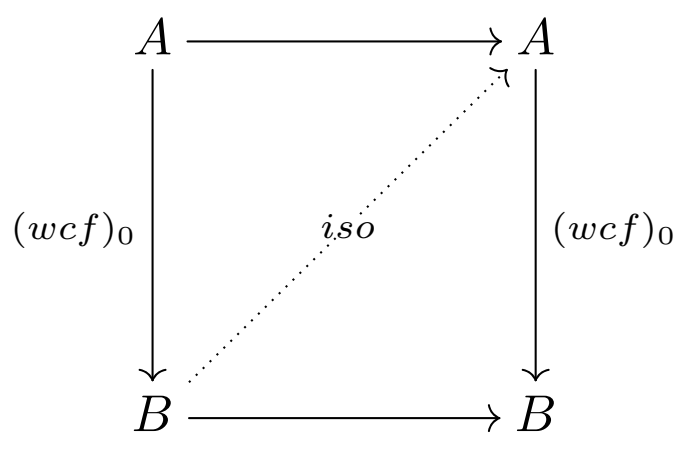

Thus, our automaton observes that whenever labelling an arrow it should not be labelled (wcf), unless it is an isomorphism. Moreover, his intention is to try to construct a closed model category. He knows that in such a model category once two of the labeles have been specified, the third is fully determined. So, in order to stay away from contradictions, he will focus on two of the labeles, letting them generate the third. Though the distinction between sets and classes is too subtle for him to understand, he knows that while $\emptyset$ is a natural object in his category and serves as the initial object, the set theoretic universe $V$, the terminal object, is not a set and has to be formally adjoined to the category in order for even Axiom (M0) to hold. Thus, it seems to him that though there is a natural homotopy theoretic duality between fibrant objects (i.e. objects $\mathcal{C}$ such that $\mathcal{C} \xrightarrow{(f)} \top$, where $\top$ is the terminal object) and co-fibrant objects, in the set theoretic context the latter should occur more naturally. Therefore he decides to try to define weak equivalences and co-fibrations (this seems more natural than focusing on co-fibrations and fibrations, since he has some hope that set theory could suggest a natural candidate for weak equivalences). But he has to be careful. His intention is to capture three set theoretic notions, and he only has two labels to achieve this with. Indeed, the situation is even more delicate: in order to understand these notions from a homotopy theoretic point of view, he has to make sure that they do not get mixed up in the passage to the associated homotopy category.

**Remark 10.** Though the construction of a co-fibrantly generated model category seems more natural in the present context, it is plausible that other approaches may also work. We do not know whether this is indeed the case, or whether the different possible approaches lead to important set theoretic concepts in the same way that the present construction leads (as will be seen in the next paper) naturally to the definition of the covering number.

So our hero tries to label the arrows to meet these goals. He first ponders which arrows to label (w). In algebraic topology (and in Quillen's model categories) these correspond to morphisms called *weak homotopy equivalences*, and indeed are thought of as some sort of equivalence, as is reflected by Axiom (M5). From the homotopy theorist point of view there is no interesting distinction between two weakly equivalent objects. When are two sets closest? When should we consider two sets "almost identical"? Two sets differing by one element are good candidates for that title. This idea is simple enough to be stumbled upon by our hero, and its connection to our hero's inner world of model categories is explicit enough to be recognised as such immediately. In our notation this corresponds to the arrow $\{a, b, ...\} \longrightarrow \{a, b, ..., \bullet\}$, so it is reasonable to label them (w).

Our hero does not know this yet, of course, but this first attempt of constructing a model category will fail. Nevertheless, it will lead him to the right track, and it will be useful to keep track of his achievements already at this early stage. In order to refer to this first approximation of a construction later on without confusing it with the actual construction, we will denote it with a 0 as subscript. Thus, for example, the arrow of the previous paragraph is denoted (without our hero knowing it) $\{a, b, ...\} \xrightarrow{(w)_0}_0 \{a, b, ..., \bullet\}$.

At this point our automaton is happy to note that by Axiom (M5) weak equivalences are stable under composition, implying, by induction, if $A \subseteq B$ and $B \setminus A$ is finite then $A \xrightarrow{(w)_0}_0 B$. This is a relation satisfying all the different cases of Axiom (M5). In particular, with this labelling, $\emptyset \xrightarrow{(w)_0}_0 B$ implies that $B$ is a finite set. Though he knows that in order to fulfil Axiom (M2) he'll have to add more weak equivalences, he hopes that in doing so he will not lose track of finite sets.

**Remark 11.** We will see later on that, in fact, if $B$ is a finite set then $\emptyset \xrightarrow{(wc)}_0 B$, and the trivial cofibrant objects will, indeed, capture the notion of finiteness. To be more precise, we will replace the set $B$ with the class $\{B\}$, and we will see that $\emptyset \xrightarrow{(wc)} \{B\}$ if and only if $B$ was a finite set.

Invigorated by this success our android decides to declare arrows co-fibrations if they connect sets of the same cardinality. In our notation this would mean that $A \xrightarrow{(c)_0}_0 B$ if $\text{card}\,(A) = \text{card}\,(B)$. Of course, our hero, who is aware of Claim 8 knows that his model category will have to have co-fibrant objects, implying that there will have to be more (c)-arrows than those he'd just declared. So he investigates what those objects might be.

He notices that Axiom (M2) requires that for any set $A$ and singleton $b$ the following diagram holds:

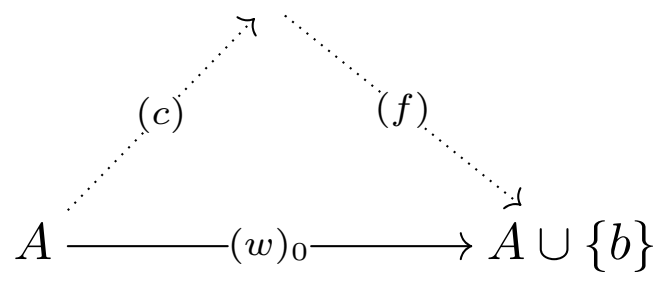

There are two possible solutions to this requirement. Either the target of the (c) arrow in the diagram is $A \cup \{b\}$, implying that $A \xrightarrow{(wc)_0}_0 A \cup \{b\}$ or the target of the (c) arrow is $A$ itself. Our android knows that for infinite sets the first option is dictated to him. So he wonders what to do with finite sets. Our hero's belief in universal axioms leads him to choose the first solution for all sets $A$ (not necessarily infinite). Since this argument goes through trivially by induction to any finite set $\{\bar{b}\}$ our android concludes that arrows should be declared co-fibrations if they connect either finite sets or sets of the same cardinality.

**Remark 12.** It is hard (or too easy) to argue rigorously in an axiomatic system we know to be inconsistent. But in fact, our android could come up with a more convincing argument for this choice.

By Claim 8 any non-trivial model category has non-trivial co-fibrant objects. Thus, there must be an infinite set $C$ such that $\emptyset \xrightarrow{(c)_0}_0 C$. Moreover, clause (4) of the same claim assures that if $\emptyset \xrightarrow{(c)} C$ then $\emptyset \xrightarrow{(wc)} \{c\}$ for all $c \in C$.

On the other hand, if $\emptyset \xrightarrow{(wf)} \mathcal{D}$ and $\mathcal{D}' \longrightarrow \mathcal{D}$ then $\emptyset \xrightarrow{(wf)} \mathcal{D}'$, as proves the following diagram:

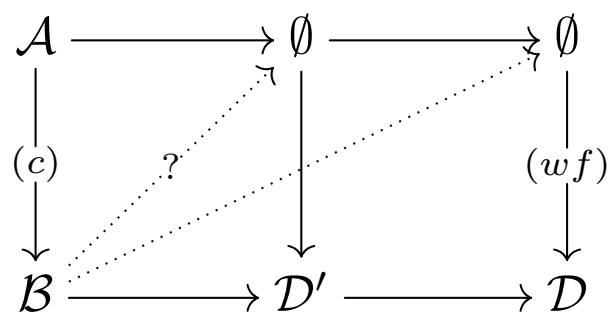

In the present context, these two observations imply that if $\emptyset \xrightarrow{(c)_0}_0 A$ and $\emptyset \xrightarrow{(f)_0}_0 B$ then $A \cap B = \emptyset$. For if $a \in A \cap B$ then, by definition $\emptyset \xrightarrow{(w)_0}_0 \{a\}$ and by what we have just seen $\emptyset \xrightarrow{(cf)_0}_0 \{a\}$. So $\emptyset \xrightarrow{(wcf)_0}_0 \{a\}$, contradicting the requirement that only isomorphisms are labelled (wcf).

Thus, if $\emptyset \xrightarrow{(f)_0}_0 A$ and $a \in A$ then $\emptyset \xrightarrow{(f)_0}_0 \{a\}$. So fix such a singleton $\{a\}$ and an infinite set $C$ such that $\emptyset \xrightarrow{(c)_0}_0 C$. So, on the one hand, by definition of the labelling, $C \xrightarrow{(wc)_0}_0 C \cup \{a\}$. By axiom (M3), $\emptyset \xrightarrow{(c)} C \cup \{a\}$. And as being a co-fibrant object is preserved when passing to subsets, $\emptyset \xrightarrow{(c)_0}_0 \{a\}$, a contradiction.

Thus, for a singleton $\{a\}$ we get $\emptyset \xrightarrow{(c)_0}_0 \{a\}$ and, by induction, the same is true if we replace $\{a\}$ with any finite set. Now let $A$ be a finite set. Then $\emptyset \xrightarrow{(wc)_0}_0 A$ and $\emptyset \xrightarrow{(wc)_0}_0 A \cup \{a\}$ for any singleton $\{a\}$. This is trivially seen to imply (in a closed model category) that $A \xrightarrow{(wc)} A \cup \{a\}$. □

Our android managed to capture the notions of infinite equi-cardinality and finiteness, but he still has to take care of the notion of countability. Finite sets are countable, so his definition of countability should be a proper generalisation of (wc)-arrows. The natural candidate is, therefore, to set $\emptyset \xrightarrow{(c)_0}_0 C$ if $C$ is countable.

But now our android runs into troubles. In his model category there should be arrows $A \longrightarrow B$ which are not weak co-fibrations. Axiom (M2) requires that there exists an object $\mathcal{A}_B$ such that $A \xrightarrow{(wc)_0}_0 \mathcal{A}_B \xrightarrow{(f)_0}_0 B$. But now the arrow $\mathcal{A}_B \xrightarrow{(f)_0}_0 B$ must have, by Axiom (M1), the right lifting property with respect to any arrow labelled $(wc)$. In particular, for any finite $b \subseteq B$ implying that if $\mathcal{A}_B$ were a set $A \cup b \subseteq \mathcal{A}_B$ for all finite $b \subseteq B$, so that $\mathcal{A}_B = B$, contradicting the assumption that $A \longrightarrow B$ is not a weak co-fibration. How can this be addressed? Our hero understands that the only possible solution (preserving the uniqueness of arrows between any two objects) is to add to his category, given $A$ and $B$ as above, a formal object - the collection of all sets $A \cup \{\bar{b}\}$ for $\bar{b} \subseteq B$ finite. Therefore, instead of working with sets he now has to work with families of sets. But which families to choose? Our simplistic hero decides to choose them all, oblivious of the fact that those will no longer be sets, but *classes*. Identifying a set $A$ with the family $\{A\}$, our android knows what the arrows between such objects should be. But how should arrows between an old object and a new one be defined, or even worse, between two new objects?

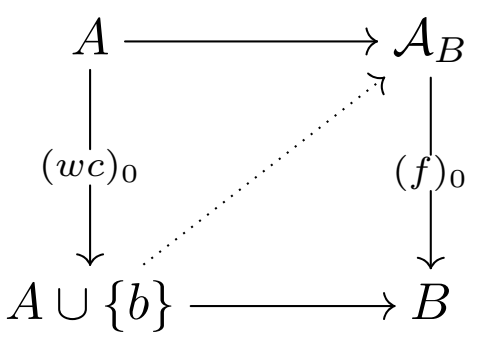

**Notation 13.** The objects of our category will be denoted by calligraphic letters ($\mathcal{A}, \mathcal{B}, \mathcal{X}, \mathcal{Y}, \mathcal{W}$ etc.) Sets we will be denoted by capital letters (e.g. $A \in \mathcal{A}$)).

Our android remembers that the only category theoretic restriction on the category he has to construct is that it should be closed under (finite) limits. So he aims to define the arrows in such a way that the resulting category will be closed under inductive and injective limits. He recalls that in the set theoretic universe he has already identified the notion of union with that of a direct limit (Figure 9). This means that for classes $\{A\}, \{B\}$ he gets an arrow $\{A\} \to \{A, B\}$. More generally, identifying a class $\mathcal{C}$ with $\bigcup_{C \in \mathcal{C}} \{C\}$, our hero immediately requires that $\mathcal{C}$ be the limit of the (possibly infinite) diagram whose nodes are $\{C\}$ for all $C \in \mathcal{C}$, concluding that

$$(\dagger) \qquad \{A\} \longrightarrow \mathcal{C} \iff A \longrightarrow_0 B \text{ for some } B \in \mathcal{C}.$$

Now, given two classes $\mathcal{C}$ and $\mathcal{D}$, our android wants to know when there should be an arrow $\mathcal{C} \longrightarrow \mathcal{D}$. Identifying $\mathcal{C}$ with the inductive limit of the diagram whose nodes are $\{C : C \in \mathcal{C}\}$ and referring to Axiom (M0), he sees that if for all $C \in \mathcal{C}$ there is an arrow $\{C\} \longrightarrow \mathcal{D}$, there should be an arrow $\mathcal{C} \to \mathcal{D}$. He notices that this relation is transitive, and therefore defines a category:

**Definition 3** (StNaamen)**.** Let *StNammen* be the category whose objects are classes of sets and whose morphisms are given by $\mathcal{X} \longrightarrow \mathcal{Y}$ if

$$(\forall X \in \mathcal{X} \exists Y \in \mathcal{Y})(X \subseteq Y).$$

**Remark 14.** Since in *StNaamen* arrows are unique, any collection of objects of the category describes a unique (necessarily commutative) diagram. Therefore one can talk of (direct or inverse) limits of any (small) collection of objects in *StNaamen* rather than on the limits of diagrams. We will use this convention freely in what follows.

Following his vocation, our hero tries to investigate lifting properties of some simple arrows in the new category, and joyfully observes, that combining simple set theoretic properties, with lifting diagrams quite often produces other set theoretic properties he has already run into.

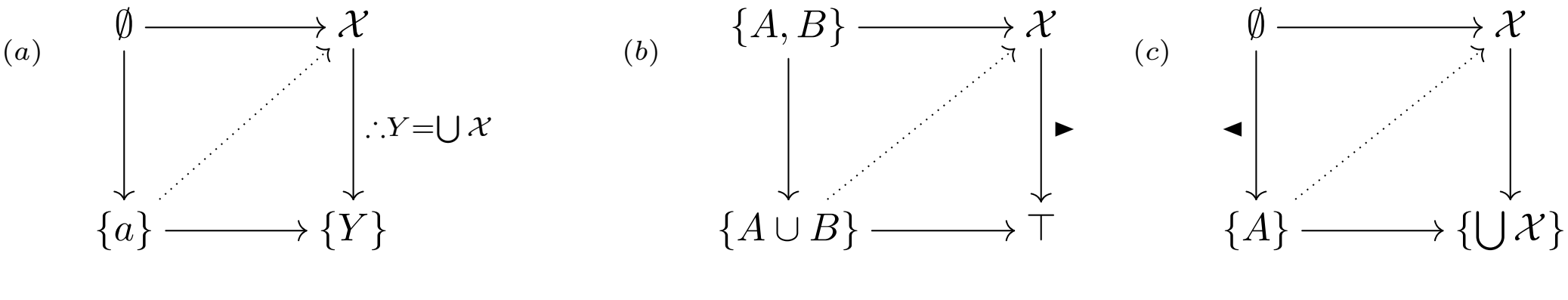

FIGURE 12. (a) Identifying a set $X$ with the class $\mathcal{X}$ whose elements are precisely the members of $X$ we have $\emptyset \longrightarrow \{a\} \rightthreetimes \mathcal{X} \longrightarrow \{Y\}$ for all $a$ if and only if $Y = \bigcup \mathcal{X}$
(b) The set (class) $\mathcal{X}$ is directed, if for all $A, B$ the lifting arrow exists. Note that the object $\mathcal{X}$ is identified with the arrow $\mathcal{X} \longrightarrow \top$ to the terminal object, $\top$ of our category.
(c) The set $A$ is finite if $\emptyset \to A \rightthreetimes \mathcal{X} \to \{\bigcup \mathcal{X}\}$ for all *linearly ordered* $\mathcal{X}$.

Of course, our hero has to use set theoretic concepts (such as linearly ordered sets, or - given a set $A$ - form the set $\{A\}$), but other than that he is quite satisfied. The world of set theory is no longer totally alien to him as he realises that in the new world objects still form a category and people apply the lifting property to define new concepts from old ones.

3.2. **Constructing a model category.** With these tools in hand, the labelling of arrows can now proceed almost axiomatically, a task on which our homotopy theorist now embarks.

3.2.1. *The labels (wc) and (f). Dual's dual argument.* In the category *StNaamen* some of the arrows are already labelled. Arrows of the form $\{A\} \longrightarrow \{B\}$ are labelled (w) if $A \xrightarrow{(w)_0}_0 B$, and similarly for the label (c). Following the guidelines of the axioms of a closed model category, the android has little choice but to adhere to Axiom (M1) and label (f) all those arrows which have the right lifting property with respect to all arrows labelled (wc). But now the same reasoning requires that we label (wc) all those arrows which have the left lifting property with respect to all the arrows we have labelled (f). This process could go on forever, but our hero does not worry. He knows from experience that:

**Claim 15.** The labelling defined satisfies the first part of Axiom 2.3.1. Labelling any other arrow (wc) or (f) will violate this axiom. All isomorphisms are labelled (wcf).

*Proof.* The argument is the same as the one showing that, for a vector space $V$, the duality $V^* = V^{***}$ always holds. Let $(wc)_0$ be the class of arrows labelled (wc) at the first step, i.e.

$$\{A\} \xrightarrow{(wc)_0} \{B\} \iff A \xrightarrow{(wc)_0}_0 B$$

Let $(f) := ((wc)_0)_{\rightthreetimes}$ be the class of all arrows $g$ such that $(wc)_0 \rightthreetimes g$, and $(wc) :={}_{\rightthreetimes}(f)$, the class of arrows $g$ such that $g \rightthreetimes (f)$. But now $(f)_2 := (wc)_{\rightthreetimes} = (f)$. To see this we apply the duality argument: as $(wc)_0 \rightthreetimes (f)$, we get $(wc)_0 \subseteq (wc)$ and infer by duality $(f) \supseteq (f)_2$. By construction $(wc) \rightthreetimes (f)$, i.e. $(f)_2 = (f)$, as required. For the second claim, observe that, by definition, if we label any arrow (wc) or (f) while preserving the first part of Axiom 2.3.1, this arrow must lie in ${}_{\rightthreetimes}(f)$ or $(wc)_{\rightthreetimes}$, respectively. But those arrows are already labelled. The last claim holds as isomorphisms have both the left and the right lifting property with respect to any arrow, as we observed in Figure 3. □

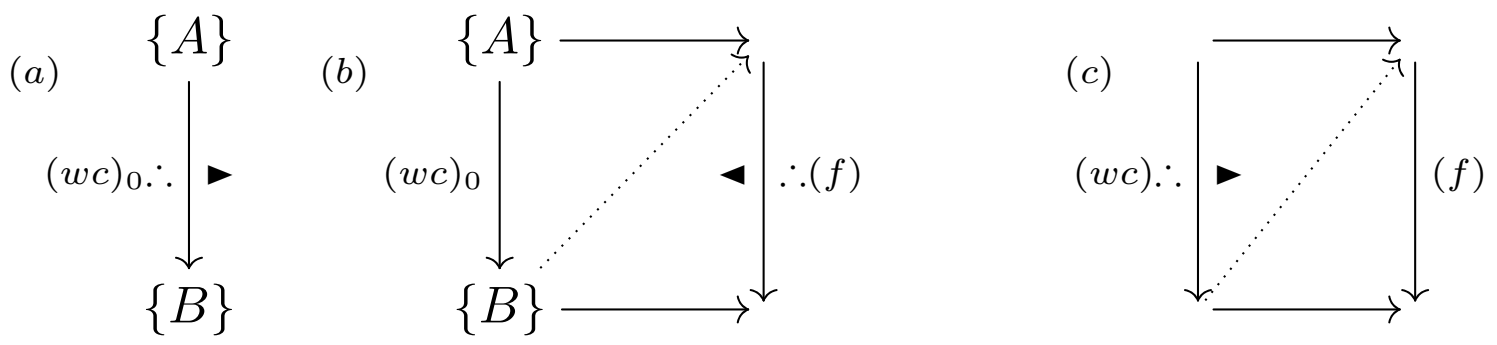

FIGURE 13. The three steps of the labelling:
(a) The basic $(wc)_0$ arrows: $B \setminus A$ is finite.
(b) The left hand side arrow is labelled (f) if it has the lifting property with respect to all arrows labelled $(wc)_0$ in step (a).
(c) the right hand side arrow is labelled (wc) if it has the lifting property with respect to all arrows labelled (f) in step (b).

Our hero is glad to find that already at this stage, a simplified definition of finiteness presents itself, one that does not have to use complicated set theoretic notions such as linearly ordered sets: a set $A$ is finite if and only if $\emptyset \xrightarrow{(wc)} \{A\}$ (see Proposition 22 for details).

**Remark 16.** This last definition of finiteness is still not a purely category theoretic one. While we have seen that the object $\emptyset$ can be defined category theoretically (being the initial object of the category), the operation $A \longmapsto \{A\}$ still has to be invoked from the world of set theory. To our android, this functor from the category of sets under inclusion to the category *StNaamen* is understood on a purely syntactic level. In all other respects, this functor is a mystery to him.

Content, the android goes on in the same vein of labelling.

3.2.2. *The labels (c) and (wf).* As before, we already have some arrows labelled (c): those labelled (wc) in the previous step (which we may ignore, see Remark 17 below) and

$$\{A\} \xrightarrow{(c)_0} \{B\} \iff A \xrightarrow{(c)_0}_0 B.$$

Using the second part of Axiom 2.3.1 our android proceeds precisely as in the previous step. More precisely, he sets $(c)_0 := \{\, \{A\} \longrightarrow \{B\} : \operatorname{card} B \leqslant \operatorname{card} A + \aleph_0 \,\}$, $(wf) := ((c)_0)_{\measuredangle}$, and finally $(c) :=_{\measuredangle} (wf)$. Evidently $(wc)_0 \subset (c)_0$ and by a similar $\measuredangle$-duality argument as above $(f) \supset (wf)$ and $(wc) \subset (c)$.

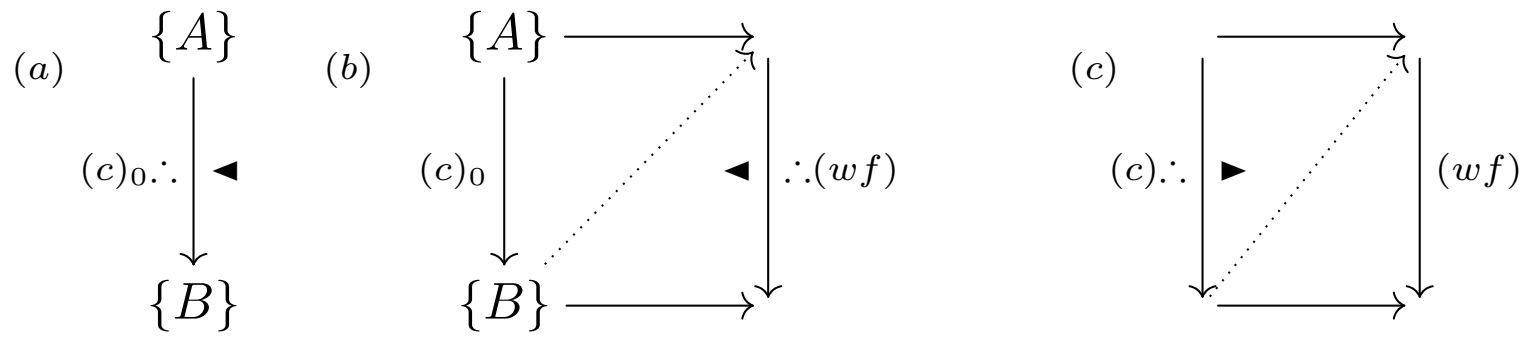

Figure 14. The three steps of the labelling:
(a) The basic $(c)_0$ arrows, $\operatorname{card} B \leqslant \operatorname{card} A + \aleph_0$.
(b) The right hand side arrow is labelled (wf) if it has the lifting property with respect to all arrows labelled $(c)_0$ in step (a).
(c) the left hand side arrow is labelled (c) if it has the lifting property with respect to all arrows labelled (wf) in step (b).

**Remark 17.** We have completed the first two steps of the construction (one described in Figure 13 and one in Figure 14). We verify that these two steps are compatible, and define, indeed, a labelling. If an arrow $h$ is labelled (wf) at stage (b) of Figure 14 this means that $(c)_0 \measuredangle h$, so in particular $(wc)_0 \measuredangle h$. Thus, the arrow $h$ has already been labelled (f) in Figure 13(b). It follows that if an arrow $g$ was labelled (wc) in Figure 13(c) then $g \measuredangle (f)$, so in particular $g \measuredangle (wf)$ and therefore $g$ is also labelled (c) in Figure 14(c).

3.2.3. *The label (w) and the Axiom M2: (w)=(wc)(wf).* Axiom (M2) requires that every arrow decomposes into the composition of a (c)-arrow and an (f)-arrow. By Axiom (M5) this implies that every (w)-arrow can be written as the composition of a (wc)-arrow and a (wf)-arrow. To satisfy this requirement, we label (w) all arrows which are the composition of a (wc)-arrow and a (wf)-arrow.

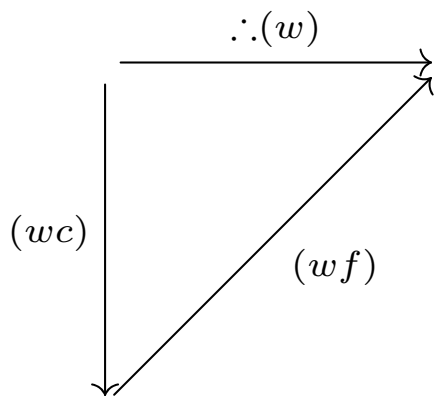

Before we proceed we need a simple observation:

**Claim 18.** The following ∀∃-diagrams are always true:

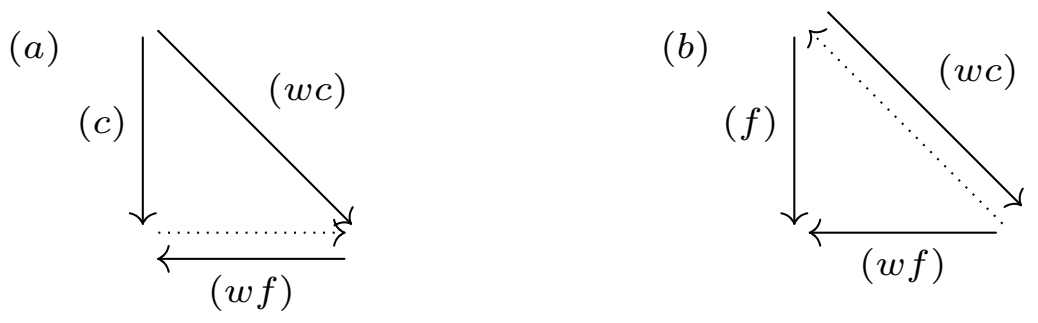

Figure 15. The dotted arrows exists whenever the solid ones are labelled as in the figures.

*Proof.* This follows trivially from our construction:

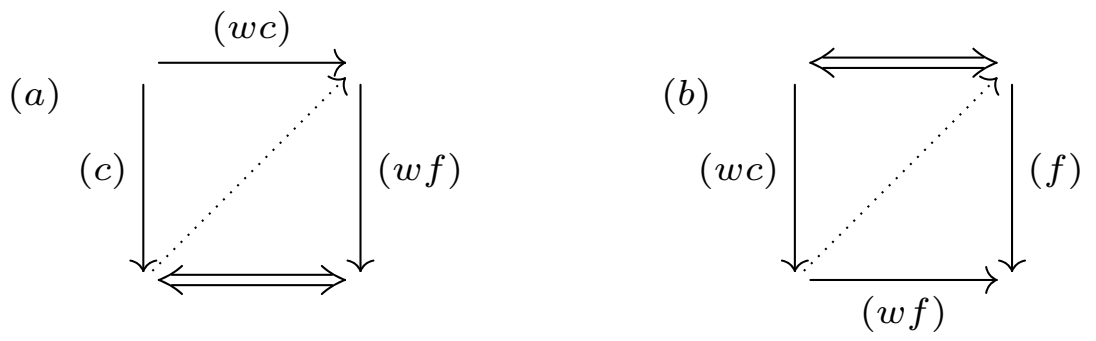

FIGURE 16. The lifting property holds in both figures, where $\Longleftrightarrow$ is the identity morphism.

□

Now, a simple diagram chasing, shows that this step of the construction is consistent:

**Claim 19.** If $\mathcal{X} \xrightarrow{(c)} \mathcal{Y}$ (resp. $\mathcal{X} \xrightarrow{(f)} \mathcal{Y}$) and if $\mathcal{X} \longrightarrow \mathcal{Y} = \mathcal{X} \xrightarrow{(wc)} \mathcal{Z} \xrightarrow{(wf)} \mathcal{Y}$ then $\mathcal{X} \longrightarrow \mathcal{Y}$ was labelled (wc) (resp. (wf)) before the last step of the construction.

*Proof.* We have to show that, under the above assumption, if $\mathcal{X} \xrightarrow{(c)} \mathcal{Y}$ then $\mathcal{X} \longrightarrow \mathcal{Y} \rightthreetimes \mathcal{V} \xrightarrow{(wf)} \mathcal{W}$ for all arrows $\mathcal{V} \xrightarrow{(wf)} \mathcal{W}$.

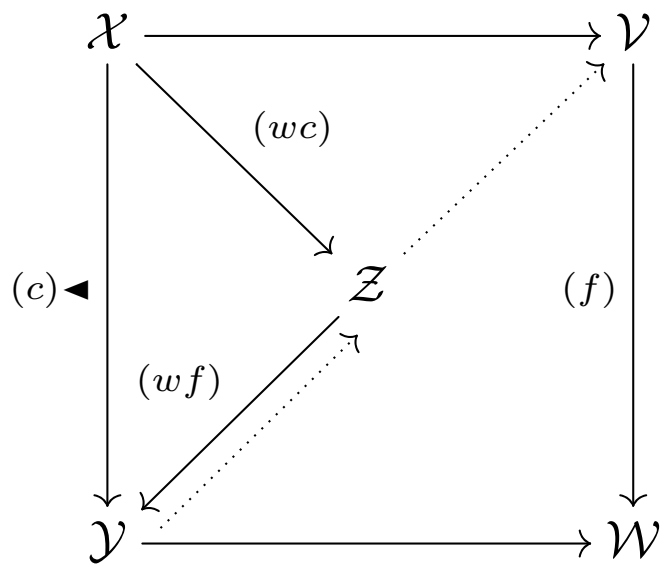

FIGURE 17. Applying the lifting $\mathcal{X} \xrightarrow{(wc)} \mathcal{Z} \rightthreetimes \mathcal{V} \xrightarrow{(f)} \mathcal{W}$ with respect to the arrows $\mathcal{X} \longrightarrow \mathcal{V}$ and the composition $\mathcal{Z} \longrightarrow \mathcal{Y} \longrightarrow \mathcal{W}$ produces the arrow $\mathcal{Z} \longrightarrow \mathcal{V}$. The arrow $\mathcal{Y} \longrightarrow \mathcal{Z}$ is given by Claim 3.2.3.

Since in our category arrows are unique, the above diagram is commutative. Thus, the composition $\mathcal{Y} \longrightarrow \mathcal{Z} \longrightarrow \mathcal{V}$ is the desired lifting. The proof of the analogous claim for $\mathcal{X} \xrightarrow{(f)} \mathcal{Y}$ is similar (dual, to be more precise). □

3.3. **Almost there.** Our hero counts his trophies. He observes that he now has some insight into the notions of *finiteness* ($\emptyset \xrightarrow{(wc)} \{A\}$), *countability* ($\emptyset \xrightarrow{(c)} \{B\}$), and *equi-cardinality* ($\{A\} \xrightarrow{(c)} \{B\}$ for infinite sets $A \subset B$). He has some understanding of the notion of *regularity* for an ordinal ($\lambda \xrightarrow{(wf)} \{\lambda\}$) (though the concept of ordinals still eludes him). He also mastered unions ($Y = \bigcup X$, see Figure 12(c)) and directed systems of sets (Figure 12(b)). He also has a feeling that he may have already constructed a model category. He knows he only has to check axioms (M2) and (M5). While our hero ponders this question, we summarise the new definitions and spell out their set theoretic meaning.

**Definition 4** (The labelling of *StNaamen*).

- **(wc)$_0$:** The arrow $\mathcal{A} \longrightarrow \mathcal{B}$ is labelled $(wc)_0$ if and only if $\mathcal{A} = \{A\}$ and $\mathcal{B} = \{B\}$ are sets and $B \setminus A$ is finite.
- **(f):** The arrow $\mathcal{X} \longrightarrow \mathcal{Y}$ is labelled (f) if and only if $(wc)_0 \rightthreetimes X \longrightarrow Y$.
- **(wf):** An arrow $\mathcal{X} \longrightarrow \mathcal{Y}$ is labelled (wf) if and only if $\{A\} \longrightarrow \{B\} \rightthreetimes \mathcal{X} \longrightarrow \mathcal{Y}$ for all sets $A, B$ such that $\operatorname{card} B \leqslant \operatorname{card} A + \aleph_0$.
- **(wc):** An arrow $\mathcal{X} \longrightarrow \mathcal{Y}$ is labelled (wc) if and only if $\mathcal{X} \longrightarrow \mathcal{Y} \rightthreetimes (f)$.
- **(c):** An arrow $\mathcal{X} \longrightarrow \mathcal{Y}$ is labelled (c) if and only if $\mathcal{X} \longrightarrow \mathcal{Y} \rightthreetimes (wf)$.
- **(w):** An arrow $\mathcal{X} \longrightarrow \mathcal{Y}$ is labelled (w) if and only if it can be written as the composition $\mathcal{X} \xrightarrow{(wc)} \mathcal{Z} \xrightarrow{(wf)} \mathcal{Y}$.

**Proposition 20.** *The labelling of the category StNammen given in Definition 4 satisfies Axioms (M0)-(M4).*

*Proof.* We have already seen (Figure 9) that *StNammen* has direct limits, which are given by unions. By duality, inverse limits are given by componentwise intersections. More precisely, given any finite set of objects $\mathcal{A}_1, \ldots, \mathcal{A}_n$ in *StNaamen* their inverse limit is

$$\mathcal{A}_1 \cap \cdots \cap \mathcal{A}_n := \bigcup \{A_1 \cap \cdots \cap A_n : A_i \in \mathcal{A}_i\}$$

as can be readily checked. Axiom (M1) holds by construction. By construction, also Axiom (M6) is satisfied. Thus, as already mentioned, Axioms (M3) and (M4) are automatic. We give the simple proofs for the sake of readers not as familiar with model categories as our hero. To verify Axiom (M3)(a), we have to check that if $\mathcal{X} \longrightarrow \mathcal{Y} = \mathcal{X} \xrightarrow{(c)} \mathcal{Z} \xrightarrow{(c)} \mathcal{Y}$ then $\mathcal{X} \xrightarrow{(c)} \mathcal{Y}$, i.e., that $\mathcal{X} \longrightarrow \mathcal{Y} \,\rangle\!\!\!\times (wf)$:

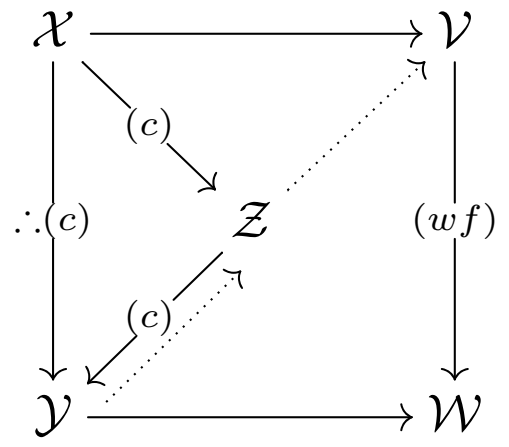

FIGURE 18. Apply the lifting $\mathcal{X} \xrightarrow{(c)} \mathcal{Z} \,\rangle\!\!\!\times \mathcal{V} \xrightarrow{(f)} \mathcal{W}$ with respect to the arrows $\mathcal{X} \longrightarrow \mathcal{V}$ and the composition $\mathcal{Z} \longrightarrow \mathcal{Y} \longrightarrow \mathcal{W}$ to produce the arrow $\mathcal{Z} \longrightarrow \mathcal{V}$. The arrow $\mathcal{Y} \longrightarrow \mathcal{Z}$ is obtained in a similar way from $\mathcal{Z} \xrightarrow{(c)} \mathcal{Y} \,\rangle\!\!\!\times \mathcal{V} \xrightarrow{(wf)} \mathcal{W}$.

A dual arguments shows that the composition of two (f)-labelled arrows is labelled (f). We have already verified that isomorphisms are labelled (wcf), so to verify Axiom (M3) it remains to check that (co)-fibrations are stable under (co)-base change. I.e. that if $\mathcal{X}' \longrightarrow \mathcal{Y}'$ satisfies the following diagram then $\mathcal{V} \xrightarrow{(wc)} \mathcal{W} \,\rangle\!\!\!\times \circledless \longrightarrow \mathcal{Y}'$:

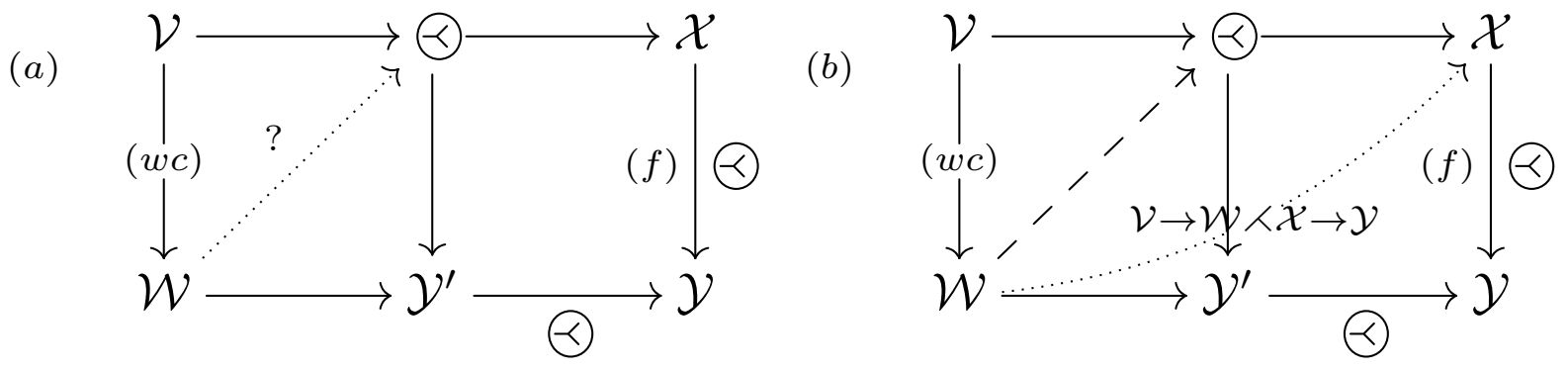

FIGURE 19. (a) We have to show that if $\circledless$ is the limit of the diagram of $\circledless$-labelled arrows then for any diagram of solid arrows there exists a lifting arrow $\mathcal{W} \longrightarrow \circledless$. (b) The arrow $\mathcal{W} \longrightarrow \mathcal{X}$ is the lifting of $\mathcal{V} \longrightarrow \mathcal{W} \,\rangle\!\!\!\times \mathcal{X} \longrightarrow \mathcal{Y}$. The arrow $\mathcal{W} \longrightarrow \circledless$ exists by the definition of limit.

By duality, the same argument works to show that co-fibrations are stable under co-base extensions. Axiom (M4) is proved in exactly the same way (the proof will show that, in fact, if an arrow $\mathcal{V} \longrightarrow \mathcal{W}$ is a base extension of a weak-fibration then $\mathcal{V} \xrightarrow{(wc)} \mathcal{W}$ and not only that $\mathcal{V} \longrightarrow \mathcal{W}$ decomposes into a weak co-fibration followed by a weak fibration).

So it only remains to verify that *StNaamen* satisfies Axiom (M2). Let $\mathcal{X} \longrightarrow \mathcal{Y}$ be any arrow. Let $\mathcal{L}$ be the (direct) limit of all the objects $\mathcal{V}$ such that $\mathcal{X} \xrightarrow{(wc)} \mathcal{V} \longrightarrow \mathcal{Y}$ (though our automaton is unaware of such problems, the next claim will show that indeed such a limit $\mathcal{L}$ exists in *StNaamen*). The corresponding diagram is:

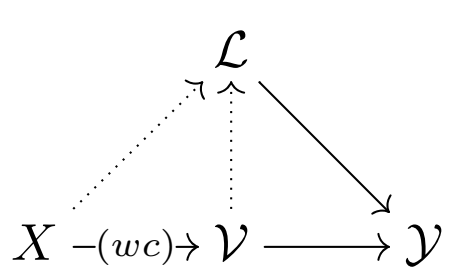

FIGURE 20

Then $\mathcal{X} \longrightarrow \mathcal{Y} = \mathcal{X} \longrightarrow \mathcal{L} \longrightarrow \mathcal{Y}$. So it remains to verify that $\mathcal{X} \xrightarrow{(wc)} \mathcal{L}$ and $\mathcal{L} \xrightarrow{(f)} \mathcal{Y}$. In order to prove the first assertion our automaton uses a simple claim:

**Claim:** Let $\{\mathcal{V}_i\}_{i\in I}$ be a class of objects in $StNaamen$ such that $\mathcal{X} \xrightarrow{(wc)} \mathcal{V}_i$ for all $i \in I$. Let $\mathcal{L}$ be the direct limit of $\{\mathcal{V}_i\}_{i\in I}$ (in particular, we assume that $\mathcal{L}$ is an object of $StNaamen$). Then $\mathcal{X} \xrightarrow{(wc)} \mathcal{L}$.

*Proof.* By the universal property of direct limits the following diagram is true:

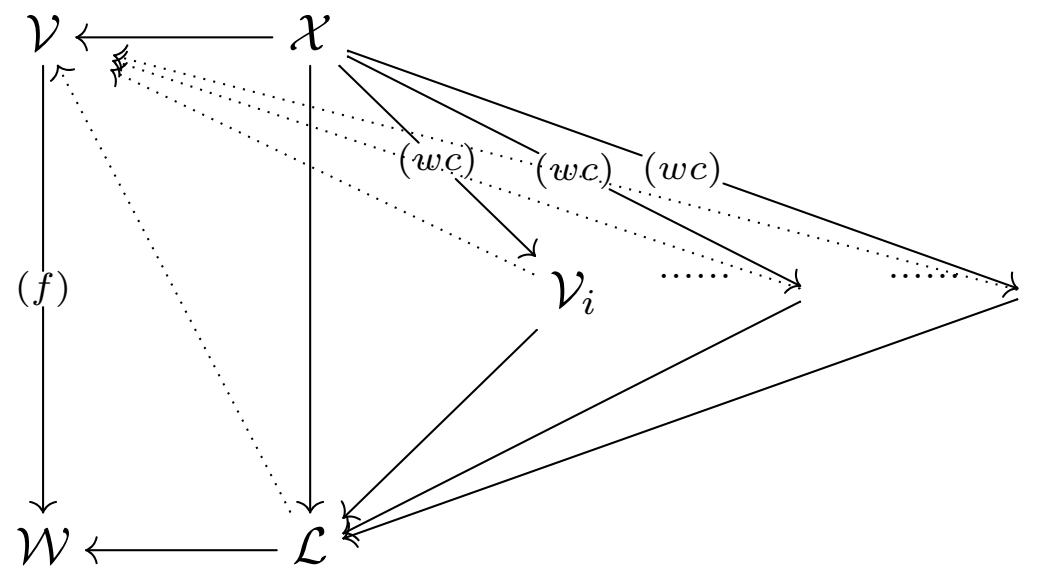

FIGURE 21. The arrows $\mathcal{V}_i \longrightarrow \mathcal{V}$ exists by $\mathcal{X} \xrightarrow{(wc)} \mathcal{V}_i \rightthreetimes \mathcal{V} \xrightarrow{(f)} \mathcal{W}$. So the arrow $\mathcal{L} \longrightarrow \mathcal{V}$ exists by virtue of it being the limit of the $V_i$.

$\square_{\text{Claim}}$

We now have to show that $\mathcal{L} \xrightarrow{(f)} \mathcal{Y}$ which, by Axiom (M6), amounts to $(wc) \rightthreetimes \mathcal{L} \longrightarrow \mathcal{Y}$. Consider the diagram:

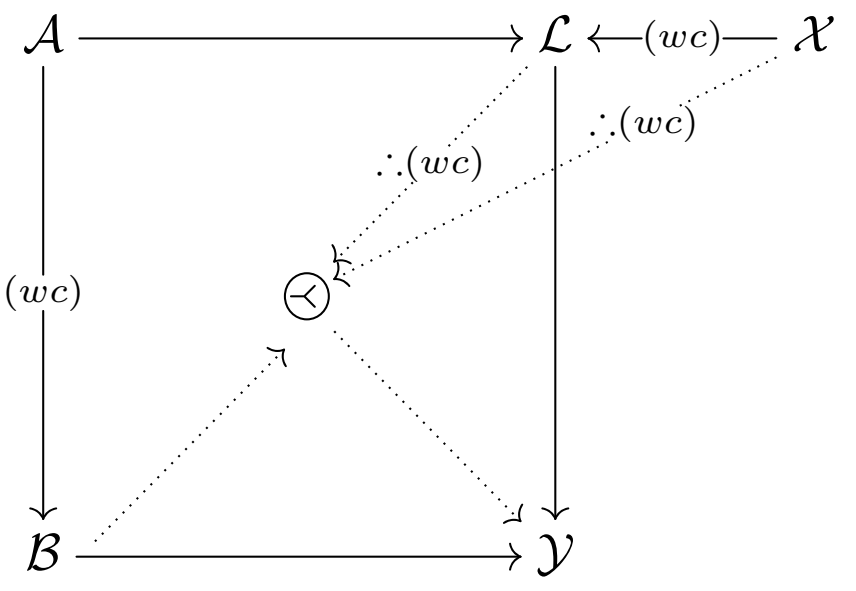

FIGURE 22. Consider the direct limit of $\mathcal{L}$ and $\mathcal{B}$. The dotted arrows exist by the universal property. The arrow connecting $\mathcal{L}$ with the limit is labelled $(wc)$ by Axiom (M3) and Axiom (M4). The arrow connecting $\mathcal{X}$ to the limit is labelled (wc) as the composition of two (wc)-arrows (see Claim 8).

It follows that $\mathcal{X} \xrightarrow{(wc)} \circledless \longrightarrow \mathcal{Y}$. By the definition of $\mathcal{L}$ this implies that $\circledless \longrightarrow \mathcal{L}$. Thus $\circledless \cong \mathcal{L}$ and as $\mathcal{B} \longrightarrow \circledless$ we also have $\mathcal{B} \longrightarrow \mathcal{L}$, showing that $\mathcal{A} \xrightarrow{(wc)} \mathcal{B} \rightthreetimes \mathcal{L} \longrightarrow \mathcal{Y}$. Since $\mathcal{A} \xrightarrow{(wc)} \mathcal{B}$ was arbitrary we are done. The other part of Axiom (M2) is proved in a similar way. $\square$

**Notation 21** (Commutative diagrams IV)**.** In a commutative diagrame a node labelled $\circledless$ denotes the limit (direct or inverse, as can be read from the direction of the arrows) of the sub-diagram whose *arrows* are labelled $\circledless$. Sometimes, in order not to over burden diagrams, we omit the $\circledless$-labelling of arrows. In such cases either no confusion can arise as to the nature of the limit, or we defer such clarifications to the caption describing the diagram.

3.3.1. *A set theoretic interlude.* Anticipating the future once again, we can now tell that our android has failed to prove in a similar way Axiom (M5). After failing long enough, he suspects that the axiom need not follow from the abstract homotopy theoretic nonsense he's been tinkering with. He still hopes, though, that the additional combinatorial information - that his objects are classes and that the arrows' labelling should have a precise set theoretic interpretation - could be applied in order to cover the gap in his argument. So he starts a project of finding the set theoretic interpretation of his labelling.

**Proposition 22.** *The set theoretic interpretation of the last definition is:*

(f) *an arrow* $\mathcal{A} \longrightarrow \mathcal{B}$ *is labelled (f) if and only if for every* $A \in \mathcal{A} \cup \{\emptyset\}$, $B \in \mathcal{B}$ *and a finite subset* $\{b_1, \ldots, b_n\} \subseteq B$ *there exists* $A' \in \mathcal{A} \cup \{\emptyset\}$ *such that* $(A \cap B) \cup \{b_1, ..., b_n\} \subseteq A'$.

(wf) *an arrow* $\mathcal{A} \longrightarrow \mathcal{B}$ *is labelled (wf) if and only if for every* $A \in \mathcal{A} \cup \{\emptyset\}$, $B \in \mathcal{B}$ *and subset* $B' \subseteq B$ *such that* $\operatorname{card} B' \leqslant \operatorname{card}(A \cap B) + \aleph_0$, *there exists* $A' \in \mathcal{A} \cup \{\emptyset\}$ *such that* $B' \subseteq A'$.

(wc) *an arrow* $\mathcal{A} \longrightarrow \mathcal{B}$ *is labelled (wc) if and only if every* $B \in \mathcal{B}$ *is contained, up to finitely many elements, in some* $A \in \mathcal{A} \cup \{\emptyset\}$ *(i.e.* $B \setminus A$ *is finite for some* $A \in \mathcal{A} \cup \{\emptyset\}$*).*

(c) *an arrow* $\mathcal{A} \longrightarrow \mathcal{B}$ *is labelled (c) if and only if for every* $\{B\} \longrightarrow \mathcal{B}$ *there exists* $A \in \mathcal{A} \cup \{\emptyset\}$ *such that* $A \xrightarrow{\mathcal{B}} B$, *where we define* $\mathcal{A} \xrightarrow{\mathcal{B}} B$ *if there exist* $n \in \mathbb{N}$ *and* $\{B_0, \ldots B_n\} \longrightarrow \mathcal{B}$ *such that:*
  (a) $\operatorname{card}(A \cap B_0) + \aleph_0 = \operatorname{card} B_0 + \aleph_0$ ,
  (b) $\operatorname{card}(B_i \cap B_{i+1}) + \aleph_0 = \operatorname{card} B_{i+1} + \aleph_0$ *for all* $0 \leq i < n$, *and*
  (c) $B = B_n$.

(w) *an arrow* $\mathcal{A} \longrightarrow \mathcal{B}$ *is labelled (w) if and only if for every* $A \in \mathcal{A} \cup \{\emptyset\}$, $B \in \mathcal{B}$ *and subset* $B' \subseteq B$ *such that* $\operatorname{card} B' \leqslant \operatorname{card}(A \cap B) + \aleph_0$, *there exists* $A' \in \mathcal{A} \cup \{\emptyset\}$ *such that* $B'$ *is contained in* $A'$ *up to finitely many elements.*

**Remark 23.** We observe that the set theoretic analysis of the last proposition is well within the abilities of our hero. He has accepted the notions of finiteness, countability and equi-cardinality as basic notions: notions he is trying to better *understand* homotopy theoretically, but he is willing to use syntactically. The remaining set theoretic concepts appearing in the proposition (unions, intersections, complements, finite complements) he mastered while constructing the category $StNaamen$.

Moreover, as we will see below, the "intuitive" understanding of the arrows $\emptyset \xrightarrow{(wc)} \{C\}$ ("$C$ is a finite set") and $\{A\} \xrightarrow{(c)_0} \{B\}$ ("$A$ is a 'large' subset of $\{B\}$") and their basic set theoretic properties suffice for our hero to produce diagram-chasing arguments to prove the proposition.

*Proof.* **(f):** By definition (see also Remark 17), an arrow $\mathcal{A} \longrightarrow \mathcal{B}$ is labelled (f) if it has the right lifting property relative to all arrows $\{C\} \xrightarrow{(wc)_0} \{D\}$, i.e. all arrows of the form $\{C\} \longrightarrow \{C, \bar{d}\}$ where $\bar{d}$ is a finite set. Thus, if $A \in \mathcal{A}$, $B \in \mathcal{B}$ and $\{b_1, \ldots, b_n\} \in B$ then

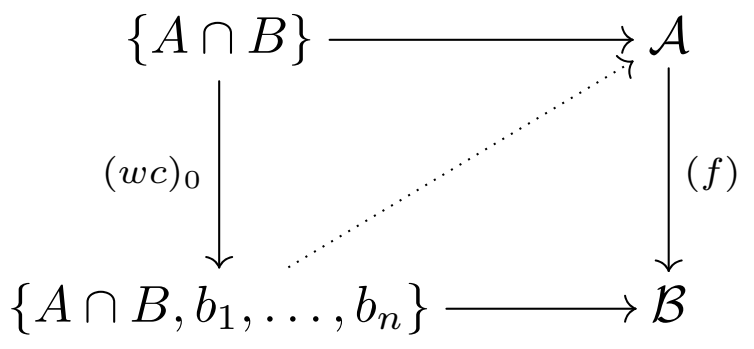

and by definition this means that there exists $B' \in \mathcal{A}$ such that $(A \cap B) \cup \{b_1, \ldots, b_n\} \subseteq A'$.
In the other direction:

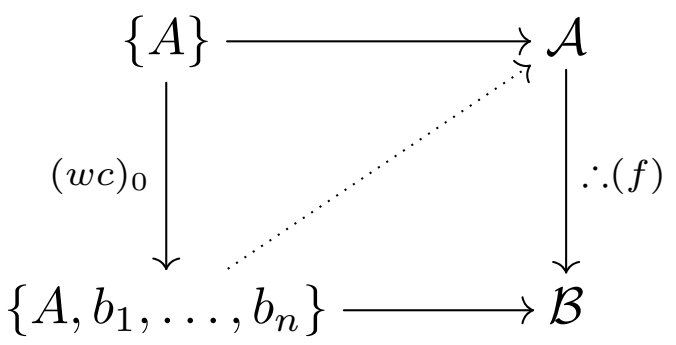

FIGURE 23. The lifting arrow exists by applying the assumption any $A \in \mathcal{A}$ and $\{\{b_1, \ldots, b_n\}\} \longrightarrow \mathcal{B}$. The conclusion $\mathcal{A} \xrightarrow{(f)} \mathcal{B}$ follows by the construction of Figutre 13.

**(wf):** Exactly the same proof works.
**(wc):** Given $\mathcal{C} \longrightarrow \mathcal{D}$ we let

$$\mathcal{C}_{\mathcal{D}} := \{(C \cap D) \cup \bar{d} : A \in \mathcal{A}, D \in \mathcal{D}, \ \bar{d} \subseteq D \text{ finite}\}.$$

Then obviously $\mathcal{C} \longrightarrow \mathcal{C}_{\mathcal{D}} \longrightarrow \mathcal{D}$, and by what we have just shown, $\mathcal{C}_{\mathcal{D}} \xrightarrow{(f)} \mathcal{D}$. Now assume that $\mathcal{C} \longrightarrow \mathcal{D} \times \mathcal{A} \xrightarrow{(f)} B$ for all $A \xrightarrow{(f)} B$. Then, in particular, $\mathcal{C} \longrightarrow \mathcal{D} \times \mathcal{C}_{\mathcal{D}} \longrightarrow \mathcal{D}$. Hence, as $StNaamen$ has unique arrows $\mathcal{D} \longrightarrow \mathcal{C}_{\mathcal{D}} \longrightarrow \mathcal{D}$ implies that $\mathcal{C}_{\mathcal{D}} \cong \mathcal{D}$. Now, let $D \in \mathcal{D}$ be any element. Then $D \subseteq C_D$ for some $C_D \in \mathcal{C}_{\mathcal{D}}$. But, by definition, $C_D$ is contained, up to finitely many elements, in some $C \in \mathcal{C}$. This proves that $\mathcal{C} \xrightarrow{(wc)} \mathcal{D}$ satisfies the combinatorial condition of the claim.

In the other direction, if $\mathcal{C} \longrightarrow \mathcal{D}$ satisfies the combinatorial condition, then

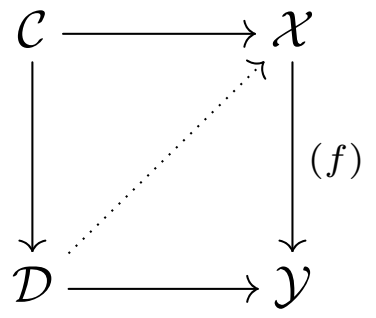

FIGURE 24. Given $D \in \mathcal{D}$ let $D \subseteq Y \in \mathcal{Y}$ and $C \in \mathcal{C}$ such that $C$ contains $D$ up to the finite set $\{d_1, \ldots, d_n\}$. Let $C \subseteq X \in \mathcal{X}$, and $X' \in \mathcal{X}$ such that $(X \cap Y) \cup \{d_1, \ldots, d_n\} \subseteq X'$, as assures $\mathcal{X} \xrightarrow{(f)} \mathcal{Y}$. So $D \subseteq X'$. Since $D$ was arbitrary, the lifting arrow exists.

By the construction of Figure 13, $\mathcal{C} \xrightarrow{(wc)} \mathcal{D}$, as required.

**(c):** The proof is quite similar to the classification of (wc)-arrows. Given $\mathcal{A} \longrightarrow \mathcal{B}$ let $\bar{\mathcal{B}} := \{B' \subseteq B : B \in \mathcal{B}\}$. Observe that $\mathcal{B} \cong \bar{\mathcal{B}}$. Define

$$\mathcal{A}^{\mathcal{B}} := \{B \in \bar{\mathcal{B}} : A \xrightarrow{\bar{\mathcal{B}}} B, \text{ some } A \in \mathcal{A}\}.$$

As before, $\mathcal{A} \longrightarrow \mathcal{A}^{\mathcal{B}} \longrightarrow \mathcal{B}$. We claim that $\mathcal{A}^{\mathcal{B}} \xrightarrow{(wf)} \mathcal{B}$. Indeed, let $A_B \in \mathcal{A}^{\mathcal{B}}$, $B \in \mathcal{B}$ and $B' \subseteq B$ be such that $\operatorname{card}(A_B \cap B) + \aleph_0 \geq \operatorname{card} B'$. We know that $A \xrightarrow{\bar{\mathcal{B}}} A_B$ for some $A \in \mathcal{A}$, so there are $B_0, \ldots, B_n \in \bar{\mathcal{B}}$ witnessing this, and $B_n = A_B$. Denote $B_{n+1} = (B \cap A_B) \cup B'$. Since, $B_{n+1} \in \bar{\mathcal{B}}$ and $\operatorname{card}(B_{n+1} \cap B_n) + \aleph_0 = \operatorname{card} B_{n+1} + \aleph_0$, we get that $A, B_0, \ldots, B_{n+1}$ witness that $A \xrightarrow{\bar{\mathcal{B}}} B_{n+1}$, so $B_{n+1} \in \mathcal{A}^{\mathcal{B}}$. Since $B' \subseteq B_{n+1}$ the combinatorial classification of (wf)-arrows yields $\mathcal{A}^{\mathcal{B}} \xrightarrow{(wf)} \mathcal{B}$.

Now assume that $\mathcal{A} \xrightarrow{(c)} \mathcal{B}$. We have to show that for $\{B\} \longrightarrow \mathcal{B}$ there exists $A \in \mathcal{A}$ such that $A \xrightarrow{\mathcal{B}} B$. By the previous paragraph, $\mathcal{A} \xrightarrow{(c)} \mathcal{B} \rightthreetimes \mathcal{A}^{\mathcal{B}} \xrightarrow{(wf)} \mathcal{B}$ and $\mathcal{B} \longrightarrow \mathcal{A}^{\mathcal{B}} \longrightarrow \mathcal{B}$. Namely, for every $\{B\} \longrightarrow \mathcal{B}$ there is $A_B \in \mathcal{A}^{\mathcal{B}}$ such that $B \subseteq A_B$, and by definition $A \xrightarrow{\bar{\mathcal{B}}} A_B$ for some $A \in \mathcal{A}$ witnessed by $B_0, \ldots, B_n \in \bar{\mathcal{B}}$. Note that $B_n = A_B \supseteq B$, so setting $B_{n+1} := B$ we get that $\operatorname{card}(B_{n+1} \cap B_n) = \operatorname{card} B_{n+1}$, implying that $A \xrightarrow{\mathcal{B}} B$.

For the other direction, we have to show that if $\mathcal{A} \longrightarrow \mathcal{B}$ satisfies the combinatorial condition, then the arrow is labelled $(c)$. By (M6) it will suffice to show that if $\mathcal{C} \xrightarrow{(wf)} \mathcal{D}$ and $\mathcal{A} \longrightarrow \mathcal{C}$, $\mathcal{B} \longrightarrow \mathcal{D}$ then for every element $B \in \mathcal{B}$ there exists $C \in \mathcal{C}$ such that $B \subseteq C$. We do this for each such $B$ separately by induction on the length of the (shortest possible) chain witnessing that $A \xrightarrow{\mathcal{B}} B$ for some $A \in \mathcal{A}$. For chains of length 0 this is obvious by construction (consider the lifting $\{A\} \xrightarrow{(c)} \{B\} \rightthreetimes \mathcal{C} \xrightarrow{(wf)} \mathcal{D}$ for any $A$ as above). Now assume that $A \xrightarrow{\mathcal{B}} B$ is witnessed by $B_0, \ldots, B_{n+1}$. Note that $\{B_0, \ldots, B_n\}$ witness that $A \xrightarrow{\mathcal{B}} B_n$. Therefore, by induction, there is some $C_n \in \mathcal{C}$ such that $B_n \subseteq C_n$. By definition, there is some $D \in \mathcal{D}$ such that $B_{n+1} \subseteq D$. Applying the combinatorial classification of (wf)-arrows to $D, C_n$ and $B_{n+1} \subseteq D$, the result follows.

**(w):** Assume, first, that $\mathcal{A} \xrightarrow{(w)} \mathcal{B}$. By construction, $\mathcal{A} \xrightarrow{(wc)} \mathcal{C} \xrightarrow{(wf)} \mathcal{B}$. Now let $A \in \mathcal{A}$, $B \in \mathcal{B}$ be any elements and $B' \subseteq B$ such that $\operatorname{card}(A \cap B) + \aleph_0 \geq \operatorname{card} B'$. Let $C \in \mathcal{C}$ such that $A \subseteq C$. So $\operatorname{card}(C \cap B) + \aleph_0 \geq \operatorname{card}(A \cap B) + \aleph_0 \geq \operatorname{card} B'$. By the classification of (wf)-arrows there is $C' \in \mathcal{C}$ such that $B' \subseteq C'$. By the classification of (wc)-arrows there is $A' \in \mathcal{A}$ such that $C'$ is contained in $A'$ up to finitely many elements. So $B'$ is contained in $A'$, up to finitely many elements as required.

Now assume that $\mathcal{A} \longrightarrow \mathcal{B}$ satisfies the combinatorial condition of the claim. We will show that there exists an object $\mathcal{C}$ such that $\mathcal{A} \xrightarrow{(wc)} \mathcal{C} \xrightarrow{(wf)} \mathcal{B}$. By Axiom (M2) we know that $\mathcal{A} \xrightarrow{(wc)} \mathcal{C} \xrightarrow{(f)} \mathcal{B}$. So our goal is to show that $\mathcal{C} \xrightarrow{(wf)} \mathcal{B}$, i.e. - using Axiom (M6) - that it has the lifting property with respect to any (c)-arrow. It will suffice to show (see Figure 14) that $(c)_0 \rightthreetimes \mathcal{C} \xrightarrow{(f)} \mathcal{B}$. Consider the following diagram:

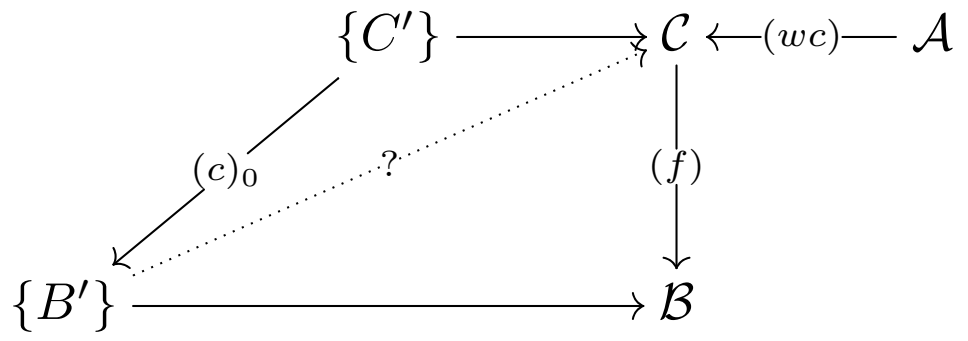

Let $B \in \mathcal{B}$ be such that $B' \subseteq B$, let $C \in \mathcal{C}$ be such that $C' \subseteq C$, and let $A \in \mathcal{A}$ be such that $C \subseteq A$ up to a finite set (as provided by the combinatorial classification of (wc)-arrows). So $\operatorname{card}(A \cap B) + \aleph_0 = \operatorname{card}(C \cap B) + \aleph_0 \geq \operatorname{card} C' + \aleph_0 = \operatorname{card} B' + \aleph_0$. Applying the assumption we get a set $A' \in \mathcal{A}$ such that $B'$ is contained in $A'$ up to a finite set $\bar{b}$. So Let $C'' \in \mathcal{C}$ be such that $A' \subseteq C''$. By the classification of (f)-arrows there exists $C''' \in \mathcal{C}$ such that $C''' \supseteq (C'' \cap B) \cup \bar{b}$. But $C'' \cap B \supseteq A' \cap B$, so $(C'' \cap B) \cup \bar{b} \supseteq (A' \cap B) \cup \bar{b} \supseteq B'$. Thus, $C''' \supseteq B'$, as required

□

As explained in Remark 23, an older model of our android - one not as capable of mastering set theoretic concepts - could have proved the above proposition using only his diagram chasing capabilities, and a few basic facts concerning finiteness, countability and cardinality. For example, the characterisation of (f)-arrows can be described as follows:

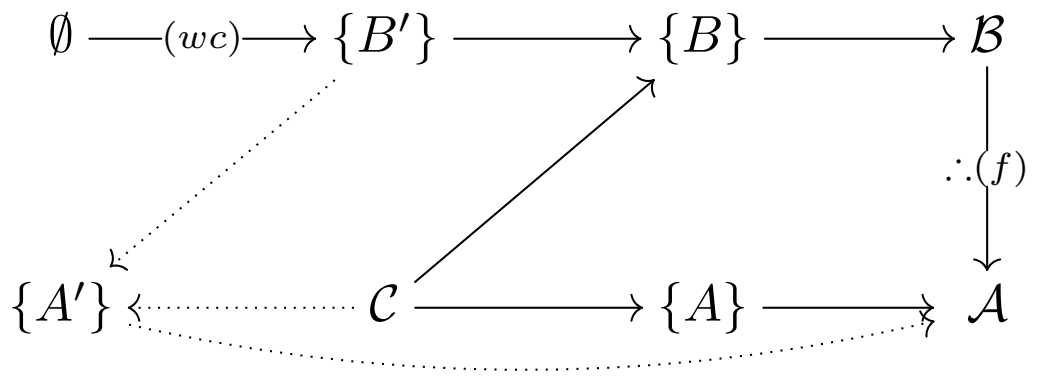

So the android has to verify that an arrow has the lifting property with respect to all $(wc)_0$-arrows if and only if it satisfies the above diagram, which is not hard even for the older release of our android to verify. Similar arguments could be applied to obtain all other classifications.

**Corollary 24.** *In StNaamen the following diagram is true:*

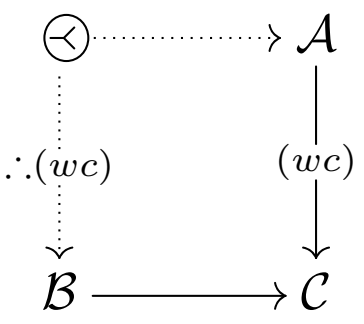

*Proof.* Recall that the inverse limit of $\mathcal{A}, \mathcal{B}$ is $\mathcal{A} \cap \mathcal{B}$. Now, if $\{B\} \longrightarrow \mathcal{B} \longrightarrow \mathcal{C}$ there exists $\{C\} \longrightarrow \mathcal{C}$ such that $\{B\} \longrightarrow \{C\}$ and $\{A\} \longrightarrow \mathcal{A}$ such that $C \setminus A$ is finite. Thus, $B$ is contained in $A$ up to finitely many elements. Since $\{B\} \longrightarrow \mathcal{B}$ was arbitrary, by Proposition 22 the corollary follows. □

3.3.2. *Axiom M5: 2-out-of-3.* At this stage our hero may have worked a little harder than he expected to, but he is quite content. He has only Axiom (M5) remaining to check in order to verify that his labelled category *StNaamen* is a model category.

**Remark 25** (Spoiler alert)**.** We could, already at this stage, give the counter-example showing that *StNaamen* does not satisfy Axiom (M5). But to reach this counter example without first checking all possibilities would be beyond our simple minded hero's abilities.

To verify Axiom (M5) our hero starts checking various special cases. He starts with arrows labelled (wc) and (wf). Of course, the composition of a (wc)-arrow with a (wf)-arrow is a (w)-arrow (by definition). He also recalls that he already knows that the composition of two (wc)-arrows is a (wc)-arrow (Claim 8), and, by duality, the composition of two (wf)-arrows is also a (wf)-arrow. He draws the following diagrams of cases that still have to be verified:

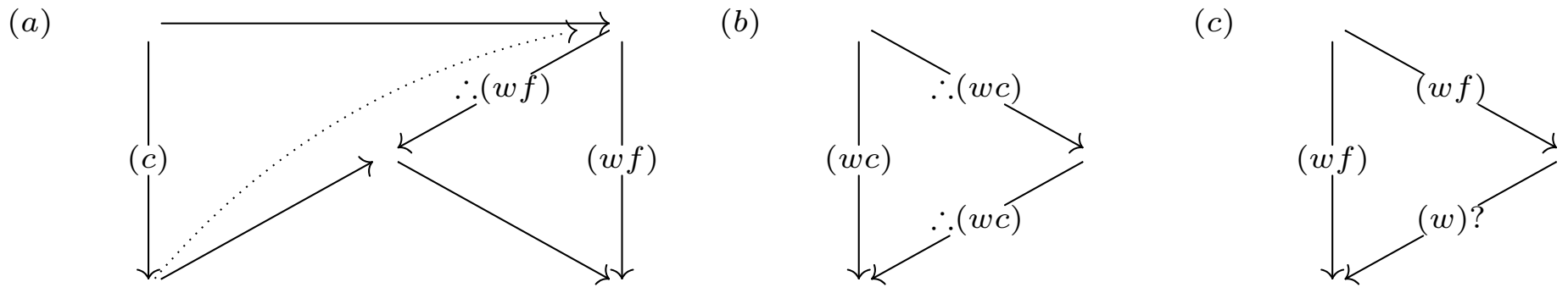

FIGURE 25. (a) If the composition of two arrows is labelled (wf) so is the first member of the composition. (b) If the composition of two arrows is labelled (wc) then, by the combinatorial classification of (wc)-arrows, so is every member of the composition (c) Here our hero got stuck.

Our android is trying to prove Diagram (c) above. He knows that, by Axiom (M2), any arrow decomposes into a (wc)-arrow followed by an (f)-arrow. He observes that this decomposition is, essentially, unique. Namely:

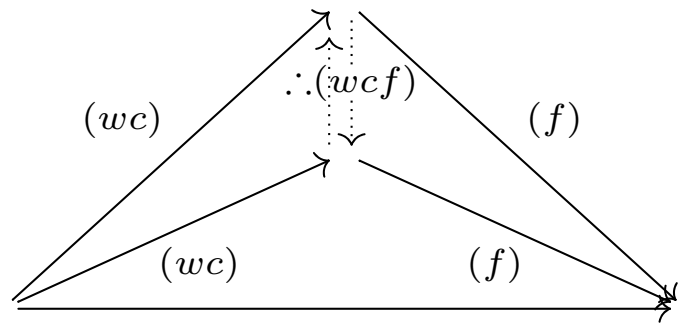

FIGURE 26. Since arrows are unique, the up and down lifting arrows are isomorphisms and therefore labelled (wcf).

He sees, therefore that if $\mathcal{X} \xrightarrow{(wc)} \mathcal{Y} \xrightarrow{(wf)} \mathcal{Z}$ then $\mathcal{X} \xrightarrow{(wc)} \mathcal{Y}' \xrightarrow{(f)} \mathcal{Z}$ implies that $\mathcal{Y}' \xrightarrow{(wf)} \mathcal{Z}$. Thus, in order to prove (c) of Figure 25 he has to show that

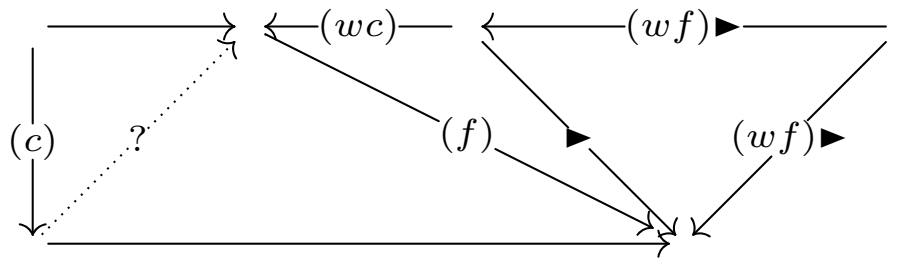

FIGURE 27. Given the diagram of ►-labelled arrows, for any system of solid arrows, our hero has to show that the lifting arrow exists.

In particular, this would be true in the special case where the bottom arrow is the identity. This gives:

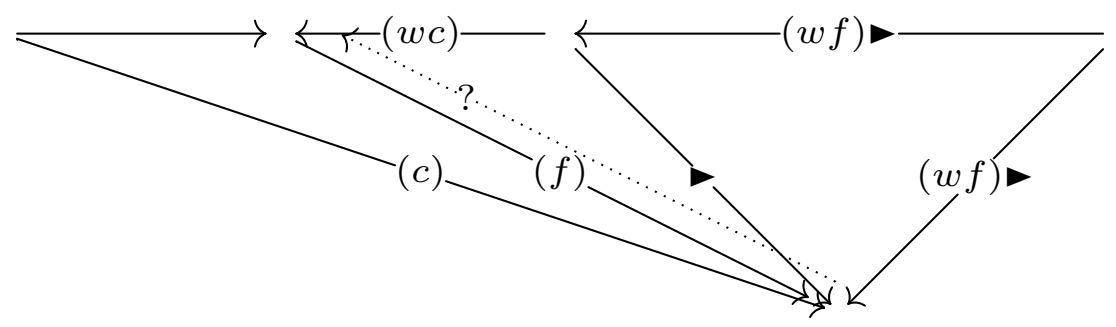

So that if the lifting arrow exists it is an isomorphism, and the diagram reduces to:

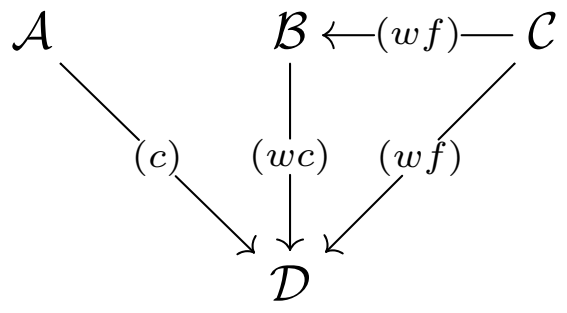

FIGURE 28

But, crestfallen, he already knows by the combinatorial classification of $StNaamen$ that for uncountable sets $A, B$ the following diagram can only be true if $A \setminus B$ is finite, which is - of course - not true in general:

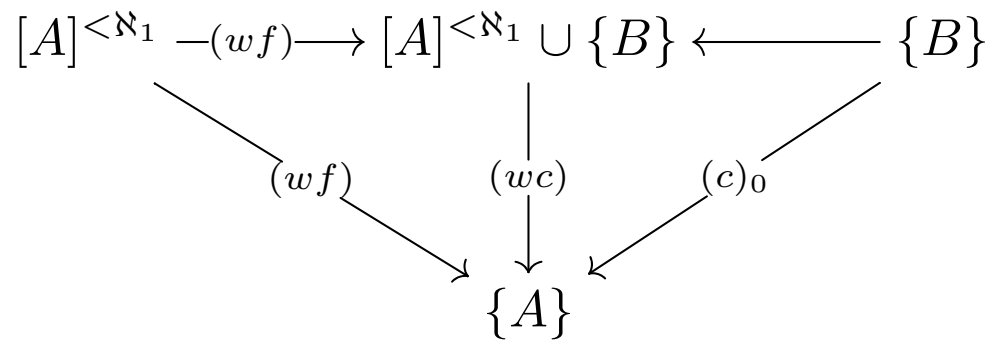

FIGURE 29. For a set $A$ we let $[A]^{<\aleph_1}$ be the set of all countable subsets of $A$. Here the arrow $\{B\} \longrightarrow [A]^{\aleph_1} \cup \{B\}$ appears in addition to the arrows required by Figure 28.

But, thankfully, our robot does not give up. He realises that some corrections have to be made. He should remove some objects (or arrows) to avoid such unfortunate counter examples. Our homotopy theorist has seen enough set theory already to realise that the objects $[A]^{<\aleph_1}$ and $\{B\}$ are quite natural in set theory, and he's reluctant to discard them. But the middle object $[A]^{<\aleph_1} \cup \{B\}$ does not seem that natural. Moreover, our hero notices that it is (w)-equivalent to $[A]^{<\aleph_1}$, so from a homotopy theorist's point of view is somewhat redundant. He is, therefore, happy to get rid of it. But how to do it?

3.4. **The model category QtNaamen.** Our hero ponders how to get rid of the problematic objects he ran into. He understands that the problem arises from some lack of "uniformity" in the classes he was considering. So he is looking for some uniformity condition:

**Definition 1.** (1) Let $\mathcal{X}$ be an object in *StNaamen*. Say that $\mathcal{X}$ *uniformly covers a set* $X$ (or, more precisely, $\mathbb{P}(X)$) if every countable subset of $X$ is contained in a member of $\mathcal{X}$, and for all $X' \subseteq X$, if $X' \subseteq C$ for some $C \in \mathcal{X}$ then every $X'' \in [X]^{\leqslant \operatorname{card} X'}$ is contained in some member of $\mathcal{X}$.

(2) Say that $\mathcal{X}$ as above is uniform, if it uniformly covers every set $X$ all of whose countable subsets are contained in members of $\mathcal{X}$.

By now our hero knows enough set theory to realise that non-uniform objects were at least one source of trouble in his previous attempt of constructing a model category. Of course, our hero is never happy with such set theoretic definitions, so he is looking for a homotopy theoretic characterisation of uniform objects. The proof of the following is obvious by Proposition 22:

**Claim 26.** An object $\mathcal{X}$ in *StNaamen* is uniform if and only if the following diagram is true:

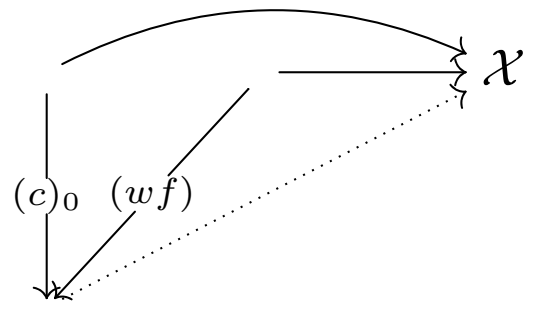

FIGURE 30

**Remark 27.** The concept of uniformity may be too complicated for our android to recognise. But he could have come up with a diagram chasing argument using the lifting property to obtain the category *QtNaamen*. He knows that he somehow has to avoid the diagram in Figure 28. He thinks of doing so by removing some of the arrows in *StNaamen*. Any arrow $\mathcal{A} \xrightarrow{(c)} \mathcal{D}$ and $\mathcal{A} \xrightarrow{(c)} \mathcal{B}$ easily fits in a diagram as in Figure 28, and dually for any arrow of the form $\mathcal{C} \xrightarrow{(wf)} \mathcal{D}$ and $\mathcal{C} \xrightarrow{(wf)} \mathcal{B}$. So the only class of arrows he can reasonably hope to eliminate is the class of those arrows which fit into a diagram as in Figure 28 as the middle arrow $\mathcal{B} \longrightarrow \mathcal{D}$.

Our automaton, trained to chase diagrams, has no use for definitions using negation. So he looks for a definition using the lifting property. He knows that given a class $K$ of morphisms in a category $\mathfrak{C}$ the class of morphisms $K^{\rightthreetimes} := \{h : K \rightthreetimes h\}$ defines naturally a subcategory of $\mathfrak{C}$. Moreover, if $h \in K$ then $h \in K^{\rightthreetimes}$ implies that $h$ is an isomorphism (because $h \rightthreetimes h$ if and only if $h$ is an isomorphism). Thus, letting $K$ be the class of all morphisms fitting into the diagram in Figure 28 as the middle arrow, we know that $K^{\rightthreetimes} \cap K$ contains only isomorphisms, but those do not lead to a counter example, so our andoid is happy to consider $K^{\rightthreetimes}$ more closely. As he is constructing a category that should have initial and terminal objects, he passes to the subcategory of $K^{\rightthreetimes}$ consisting of those objects, $\mathcal{X}$, such that $\emptyset \longrightarrow \mathcal{X}$ and $\mathcal{X} \longrightarrow \top$ are both in $K^{\rightthreetimes}$. Setting $\mathfrak{C} = StNaamen$ this is readily seen to be a *full* sub-category. Indeed, an arrow $\mathcal{X} \longrightarrow \mathcal{Y}$ is in $K^{\rightthreetimes}$ if and only if $\mathcal{X} \longrightarrow \top$ is in $K^{\rightthreetimes}$.

Thus, it is now a easy exercise to verify that $\mathcal{X} \longrightarrow \top$ and $\emptyset \longrightarrow \mathcal{X}$ are both in $K^{\prec}$ if and only if $\mathcal{X}$ is an object in $QtNaamen$.

Our hero is still not satisfied with the definition of Figure 30. On the set theoretic level $(c)_0$-arrows are simple enough to understand, but he does not feel comfortable enough manipulating them homotopy theoretically. So he verifies:

**Lemma 28.** *Let $\mathcal{X}$ be an object of $StNaamen$. Then $\mathcal{X}$ is uniform if and only if the following diagram is true for every choice of solid arrows.*

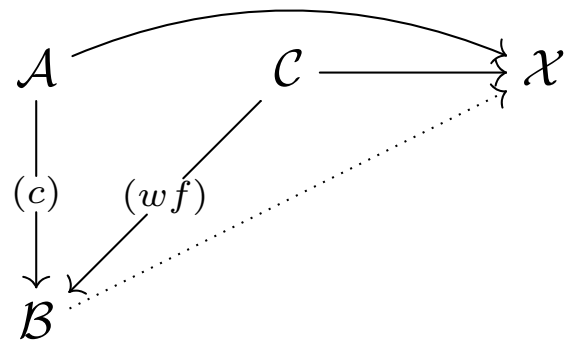

FIGURE 31

*Proof.* Even our hero knows that if $\mathcal{X}$ satisfies the diagram it must be uniform. So he proceeds to prove the other direction. He supposes that $\mathcal{X}$ is uniform, and verifies that it satisfies the above diagram. When given the above diagram he has to find a lifting arrow. He also knows that it will suffice to show that for every $B \in \mathcal{B}$ there exists an arrow $\{B\} \longrightarrow \mathcal{X}$. Fixing such a $B$ and expanding the set theoretic definition of $\mathcal{A} \xrightarrow{(c)} \mathcal{B}$ he chooses $A$ and $\{B_0, \ldots, B_n\}$ (with $B_n = B$) witnessing that $A \xrightarrow{\mathcal{B}} B$. He obtains the following diagram:

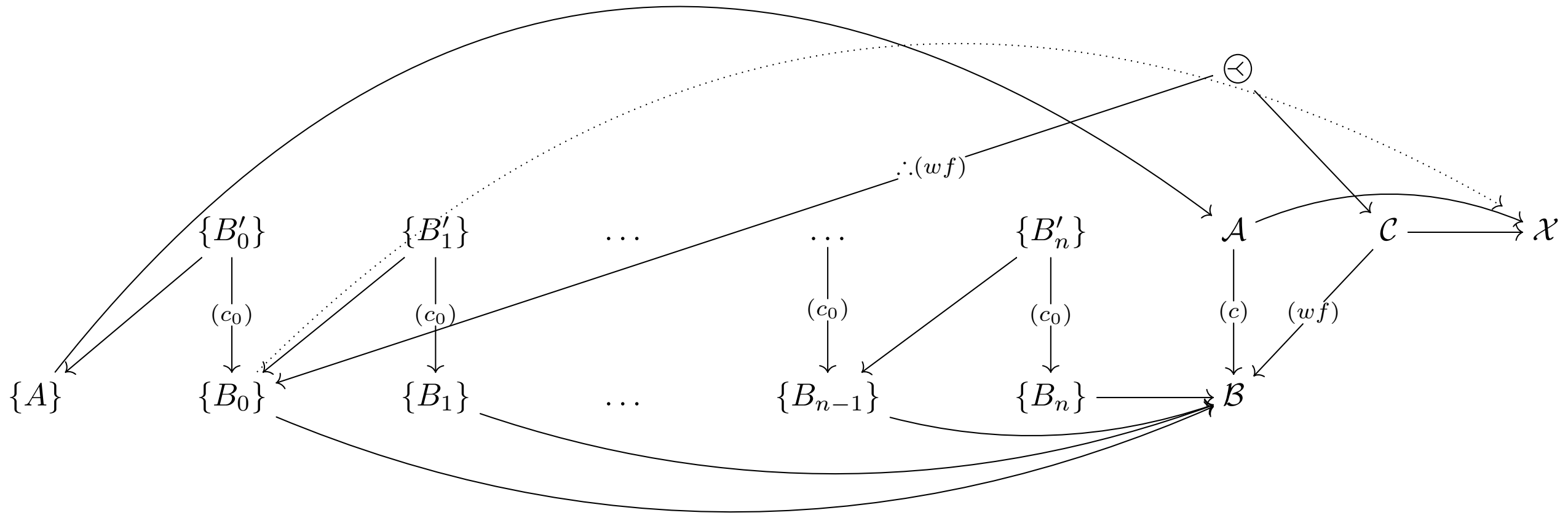

FIGURE 32. To obtain the labelling of the arrow from the limit of $\{B_0\}$ and $\mathcal{C}$ to $\{B_0\}$ use Axiom (M3) and Axiom (M4). To obtain the lifting arrow use the assumption with respect to $\{B'_0\} \xrightarrow{(c)_0} \{B_0\}$ and the previously discussed arrow. Observe that $\{B'_0\} \longrightarrow \mathcal{X}$ as the composition $\{B'_0\} \longrightarrow \{A\} \longrightarrow \mathcal{A} \longrightarrow \mathcal{X}$.

Now he notes that if $\mathcal{A} \xrightarrow{(c)} \mathcal{B}$ and $\{B\} \longrightarrow \mathcal{B}$ then $\mathcal{A}_B := \mathcal{A} \cup \{B\}$ also satisfies $\mathcal{A}_B \xrightarrow{(c)} \mathcal{B}$. So, using the lifting arrow he managed to find in Figure 32 he sees that he can replace $\mathcal{A}$ with $\mathcal{A}_B$, and that now $B$ has a chain of length at most $n-1$ witnessing that $B_0 \xrightarrow{\mathcal{B}} B$. So induction gives the desired result. □

Summing up the above discussion it may be worth pointing out :

**Lemma 29.** *Let $\mathcal{X} \in ObStNaamen$. Then the following conditions are equivalent:*

(1) $\mathcal{X} \in \mathcal{O}b$QtNaamen.

(2) *If* $\mathcal{C} \xrightarrow{(wf)} \mathcal{B}$ *then* $\mathcal{A} \xrightarrow{(c)} \mathcal{B} \measuredangle \mathcal{X} \longrightarrow T$ *and if* $\mathcal{A} \xrightarrow{(c)} \mathcal{B}$ *then* $\mathcal{C} \xrightarrow{(wf)} \mathcal{B} \measuredangle \mathcal{X} \longrightarrow T$ *for all objects* $\mathcal{A}, \mathcal{B}$ *and* $\mathcal{C}$ *of StNaamen such that* $\mathcal{A}, \mathcal{C} \longrightarrow \mathcal{X}$.
(3) $\bigcup\{X'' :$ *for some* $X_0, X_1, X'' \xleftarrow{(wf)} X_0 \longrightarrow \mathcal{X} \longleftarrow X_1 \xrightarrow{(c)} X''\} \dashrightarrow \mathcal{X}$
(4) *If* $\{a\} \longrightarrow \mathcal{X}$, $\mathcal{C} \longrightarrow X$ *and* $\{a\} \xrightarrow{(c)} \{b\}$ *and* $\mathcal{C} \xrightarrow{(wf)} \{b\}$, *then* $\{b\} \longrightarrow \mathcal{X}$.
(5) *For any set* $A$, *if* $A^{\leqslant \aleph_0} \longrightarrow \mathcal{X}$ *then* $A^{\leqslant \operatorname{card}(x \cap A)} \longrightarrow \mathcal{X}$ *for all* $x \in \mathcal{X}$ *(where* $A^{\leqslant \lambda} := \{L \subseteq A : \operatorname{card} L \leqslant \lambda\}$*)*

At last, our hero is satisfied with the definition of uniform objects. He now wonders whether there are really *enough* uniform objects in $StNaamen$:

**Lemma 30.** *Assume that* $\mathcal{C}$ *is a co-fibrant object of StNaamen. Then* $\mathcal{C}$ *is uniform.*

*Proof.* This follows immediately from Lemma 22 and the set theoretic characterisation of uniform objects. But our automaton is happier chasing diagrams, so he starts drawing them:

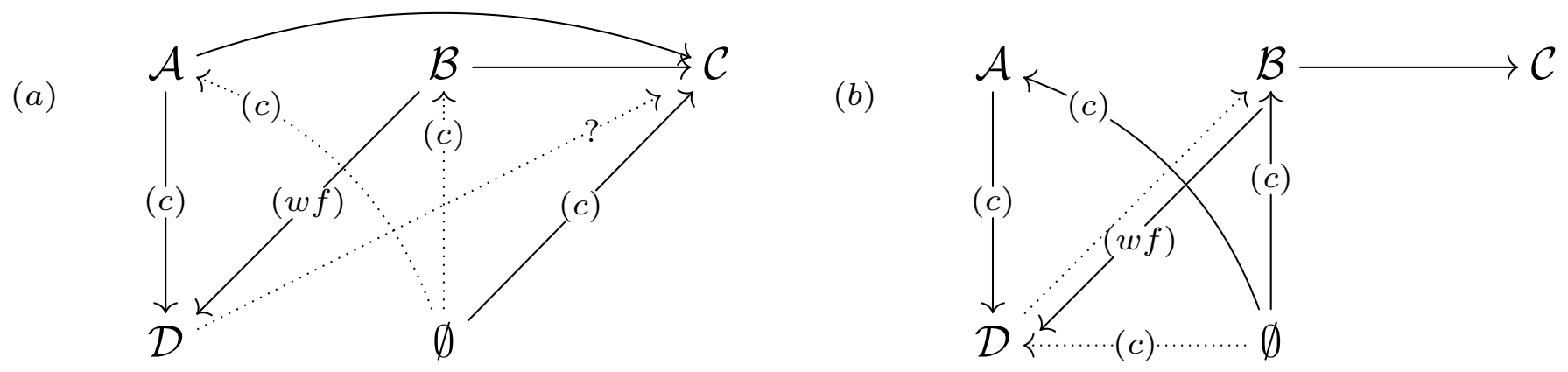

FIGURE 33. (a) The new (c)-arrows exists by Claim 8(4). (b) The arrow $\emptyset \longrightarrow \mathcal{D}$ is labelled (c) as a combination of $\emptyset \xrightarrow{(c)} \mathcal{A} \xrightarrow{(c)} \mathcal{D}$. The lifting arrow $\mathcal{D} \longrightarrow \mathcal{B}$ exists by $\emptyset \xrightarrow{(c)} \mathcal{D} \measuredangle \mathcal{B} \xrightarrow{(wf)} \mathcal{D}$.

□

Note that given any object $\mathcal{A}$ in $StNaamen$ by Axiom (M2) we know that $\emptyset \longrightarrow \mathcal{A} = \emptyset \xrightarrow{(c)} \mathcal{C} \xrightarrow{(wf)} \mathcal{A}$ for some $\mathcal{C}$. So any object in $StNaamen$ is weakly equivalent to a co-fibrant, hence a uniform object. This makes our robot quite happy, and he defines:

**Definition 2.** The category $QtNaamen$ is the full subcategory induced on the collection of uniform objects in $StNaamen$. The labelling of arrows in $QtNaamen$ is induced form the labelling in $StNaamen$.

**Remark 31.** It follows right from the definition that for every set $A$ the singleton $\{A\}$ is an object in $QtNaamen$. In particular, all $(c)_0$-arrows survive the passage to $QtNaamen$.

**Notation 32.** Let $\mathcal{X}$ be any object in $StNaamen$. We denote $\mathcal{X}_c$ the unique object such that $\emptyset \xrightarrow{(c)} \mathcal{X}_c \xrightarrow{(wf)} \mathcal{X}$.

How far could $QtNaamen$ be from being a model category? Clearly, Axioms (M3) and (M4) of model categories are preserved when restricting to an induced subcategory. Since the subcategory is full, the same is also true of Axiom (M1). But the remaining axioms have to be re-verified (since it is not clear, for example, that if $\mathcal{A} \longrightarrow \mathcal{B} = \mathcal{A} \xrightarrow{(wc)} \mathcal{C} \xrightarrow{(f)} \mathcal{B}$ and $\mathcal{A}, \mathcal{B}$ are objects in $QtNaamen$ so is $\mathcal{C}$. Moreover, though $QtNaamen$ may have limits (Axiom (M0)) it is clear from the example analysed in the end of the previous sub-section that those will not, in general, be the same limits as in $StNaamen$.

It is his first goal to prove that $QtNaamen$ has limits (otherwise the situation is hopeless already on the category theoretic level). His original hope is, given an object $\mathcal{X}$ in $StNaamen$ to find a *Qt-fication* of $\mathcal{X}$, namely, an object $\tilde{\mathcal{X}}$ such that:

(1) $\tilde{\mathcal{X}}$ is in $QtNaamen$.
(2) $\mathcal{X} \longrightarrow \tilde{\mathcal{X}}$.
(3) If $\mathcal{Y}$ is in $QtNaamen$ and $\mathcal{X} \longrightarrow \mathcal{Y}$ (in $StNaamen$) then $\tilde{\mathcal{X}} \longrightarrow \mathcal{Y}$.

So $\tilde{\mathcal{X}}$ (if it exists) is the object in $QtNaamen$ *closest* to $\mathcal{X}$. In particular, this would imply that, given any diagram. $\mathcal{D}$, in $QtNaamen$ if $\mathcal{L}$ is the direct limit of $\mathcal{D}$ (in $StNaamen$) then $\tilde{\mathcal{L}}$ must be its limit in $QtNaamen$.

**Remark 33.** The third point in the definition of the $Qt$-fication can be described as a lifting property of two classes of arrows that our android knows (and cares) about already: $\mathcal{X} \longrightarrow \tilde{\mathcal{X}} \measuredangle \mathcal{Y} \longrightarrow \mathcal{Z}$ for every arrow $\mathcal{Y} \longrightarrow \mathcal{Z}$ in $QtNaamen$. To see that this implies (3), take $\mathcal{Z} = \top$ the terminal object: if $\mathcal{Y}$ is in $QtNaamen$ and $\mathcal{X} \longrightarrow \mathcal{Y}$ (in $StNaamen$) then $\mathcal{Y} \longrightarrow \top$ is in $QtNaamen$, $\tilde{\mathcal{X}} \longrightarrow \top$ and the lifting property implies $\tilde{\mathcal{X}} \longrightarrow \mathcal{Y}$. The uniqueness of the $Qt$-fication implies that, indeed, the two conditions are equivalent.

Thus, our hero's first goal is to find for every $\mathcal{X}$ in *StNaamen* a Qt-fication $\tilde{\mathcal{X}}$. From the homotopy theoretic point of view, it is clear to our homotopy theory automaton what to do: an object $\mathcal{X}$ in *StNaamen* is not in *QtNaamen* if and only if there are $\mathcal{A}, \mathcal{B}$ and $\mathcal{C}$ as in the solid arrows of Figure 31, but such that the lifting arrow $\mathcal{B} \longrightarrow \mathcal{C}$ does not exist. The object in *StNaamen* "closest" to $\mathcal{X}$ for which such a lifting arrow exists is the direct limit of $\mathcal{X}$ and $\mathcal{B}$. Thus, his first candidate for the *Qt*-fication of $\mathcal{X}$ is the direct limit of all objects $\mathcal{B}$ for which there exist $\mathcal{A}, \mathcal{B}$ such that Figure 31 holds of $\mathcal{A}, \mathcal{B}, \mathcal{C}$ and $\mathcal{X}$.

**Notation 34.** Given $\mathcal{X}$, an object of *StNaamen* denote

$$\tilde{\mathcal{X}} := \bigcup \{\mathcal{X}'' : \mathcal{X} \longleftarrow \mathcal{X}_0 \xrightarrow{(c)} \mathcal{X}'' \xleftarrow{(wf)} \mathcal{X}_1 \longrightarrow \mathcal{X}\}.$$

In other words, given $\mathcal{X}$ the object $\tilde{\mathcal{X}}$ is the direct limit of all those objects $\mathcal{X}''$ for which the following diagram holds:

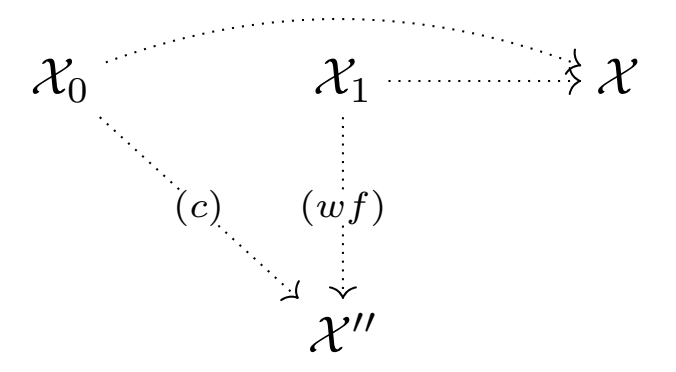

FIGURE 34. The object $\tilde{\mathcal{X}}$ is the direct limit of all $\mathcal{X}''$ satisfying the diagram.

**Remark 35.** A similar definition could be reached from the set theoretic angle. Let $\mathcal{X}$ be an object in *StNaamen*. If $\mathcal{X}$ has a Qt-fication, $\tilde{\mathcal{X}}$, then for every set $\{B\}$ such that every countable subset of $B$ is contained in some element of $\mathcal{X}$, if $B' \subseteq B$ is such that $B' \subseteq X$ for some $X \in \mathcal{X}$ then $[B]^{\leq \operatorname{card} B'} \longrightarrow \tilde{\mathcal{X}}$. Diagrammatically, this can be represented by

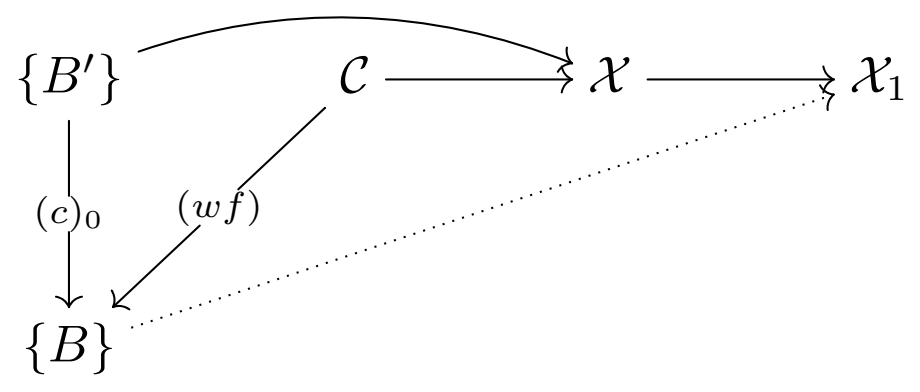

Set theoretically, $\mathcal{X}_1$ is simply the union of all $\{B\}$ as in the above diagram. So $\mathcal{X}_1$ is uniform - but only with respect to sets embedded in $\mathcal{X}$. The automaton denotes $\mathcal{X}_0 = \mathcal{X}$ and $\mathcal{X}_1$ the resulting set object. However, he cannot see any reason why $\mathcal{X}_1$ should be in *QtNaamen*. Indeed, with a little luck, he can even find a counter example. So - as he obtained $\mathcal{X}_1$ from $\mathcal{X}_0$ - goes on inductively to define $\mathcal{X}_2, \mathcal{X}_3, \ldots$ and $\tilde{\mathcal{X}} = \bigcup_{i \in \omega} \mathcal{X}_i$, and it is a straightforward argument to verify that $\tilde{\mathcal{X}}$ is the *Qt*-fication of $\mathcal{X}$. However, simple as it may seem, this set theoretic proof of the existence of *Qt*-fications is too complicated for our automaton.

It may be worth pointing out, thought that, using Proposition 22, the two definitions of $\tilde{\mathcal{X}}$ coincide. This would follow, of course, also from the uniqueness of the *Qt*-fication, once we know that both definitions are in fact *Qt*-fications.

At all events, our hero now wants to prove:

**Lemma 36.** *For every $\mathcal{X}$ in StNaamen the object $\tilde{\mathcal{X}}$ is the Qt-fication of $\mathcal{X}$.*

He has three properties to verify, of which the requirement that $\mathcal{X} \longrightarrow \tilde{\mathcal{X}}$ is obvious. So he immediately proceed to proving:

**Claim 37.** Let $\mathcal{V}$ be an object in *QtNaamen*. Then for all $\mathcal{X}$, if $\mathcal{X} \longrightarrow \mathcal{V}$ then $\tilde{\mathcal{X}} \longrightarrow \mathcal{V}$.

*Proof.* Let $\{B\} \in \tilde{\mathcal{X}}$. So the following diagram is true:

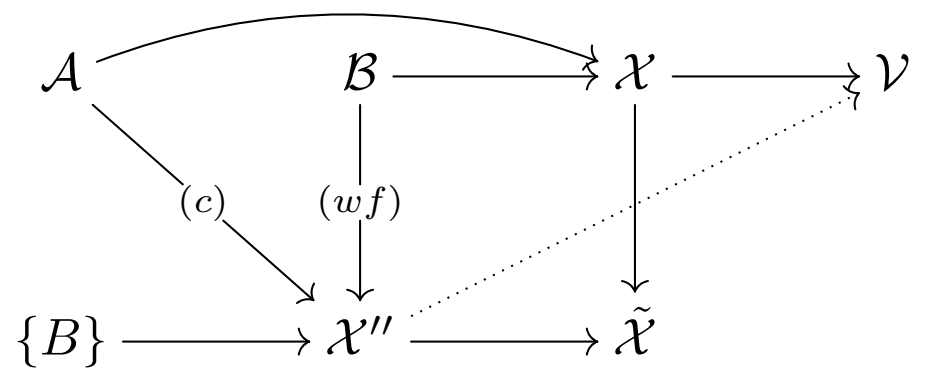

FIGURE 35. The objects $\mathcal{A}, \mathcal{B}, \mathcal{X}''$ (with their associated labelled arrows) exist by the definition of $\tilde{X}$. The arrow $\mathcal{X}'' \longrightarrow \mathcal{V}$ exists because $\mathcal{V}$ is in $QtNaamen$.

□

Anticipating the future, we isolate a simple observation that our android will run into shortly:

**Claim 38.** For every object $\mathcal{A}$ in $StNaamen$ we have $\emptyset \xrightarrow{(c)} \mathcal{A} \not\!\!\times \mathcal{X} \longrightarrow \tilde{\mathcal{X}}$.

*Proof.* Assume first that $\mathcal{A}$ is the singleton $\{A\}$ and $\emptyset \xrightarrow{(c)_0} \{A\}$. By definition of $\tilde{\mathcal{X}}$ the fact that $\{A\} \longrightarrow \tilde{\mathcal{X}}$ implies that there exists $\mathcal{X}''$ such that $\{A\} \longrightarrow \mathcal{X}''$ and such that $\mathcal{X}''$ is as in Figure 34. In particular there exists some $\mathcal{X}_1$ such that $\mathcal{X}_1 \xrightarrow{(wf)} \mathcal{X}''$. Therefore, we get the following diagram:

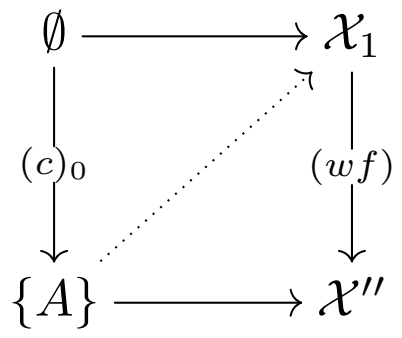

So that $\{A\} \longrightarrow \mathcal{X}_1 \longrightarrow \mathcal{X}$ is the required arrow. Finally, observe that by the combinatorial classification of co-fibrant objects, we know that $\mathcal{A}$ is co-fibrant if and only if $\emptyset \xrightarrow{(c)_0} \{A\}$ for all $A \in \mathcal{A}$. So a co-fibarant object is a direct limit of sets $\emptyset \xrightarrow{(c)_0} A$, and since the lifting arrow exists for every such $A$ it also exists for the limit. □

Now our android can prove:

**Lemma 39.** *For every object $\mathcal{X}$ in StNaamen the object $\tilde{\mathcal{X}}$ is in QtNaamen.*

*Proof.* Our hero has to show that for any choice of solid arrows the following diagram is true:

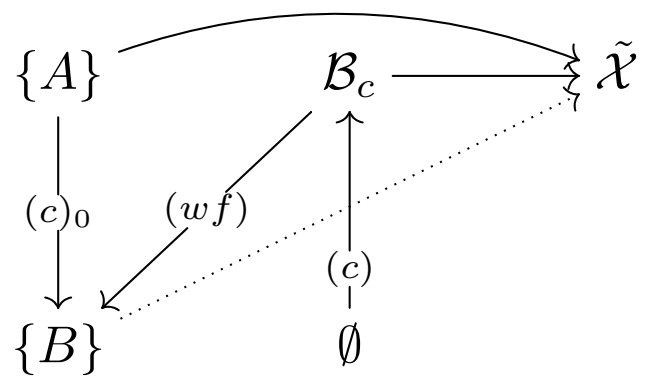

Since $\{A\}$ is a singleton, $\{A\} \longrightarrow \tilde{\mathcal{X}}$ implies that there exists $\mathcal{X}'' \longrightarrow \tilde{\mathcal{X}}$ such that for some $\mathcal{C}, \mathcal{D}$

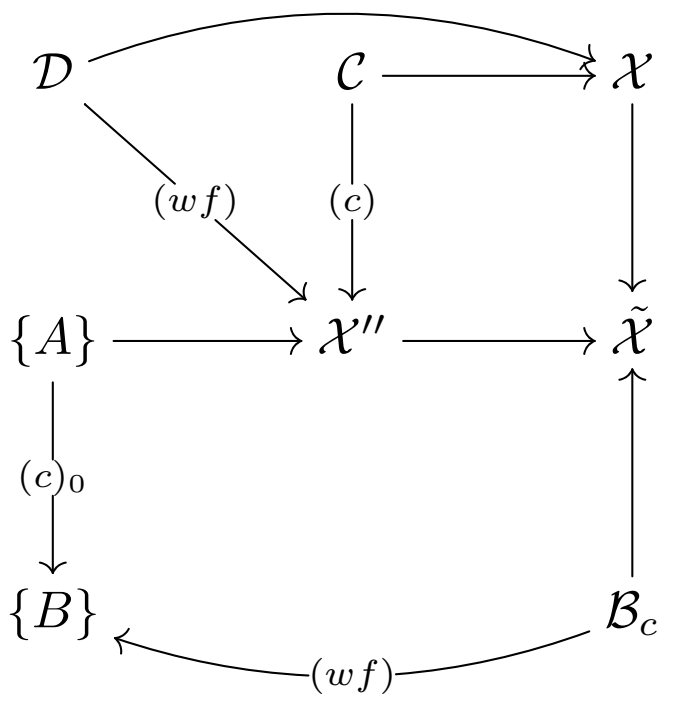

So he now considers the limit of $\mathcal{X}''$ and $\{B\}$, to obtain the following:

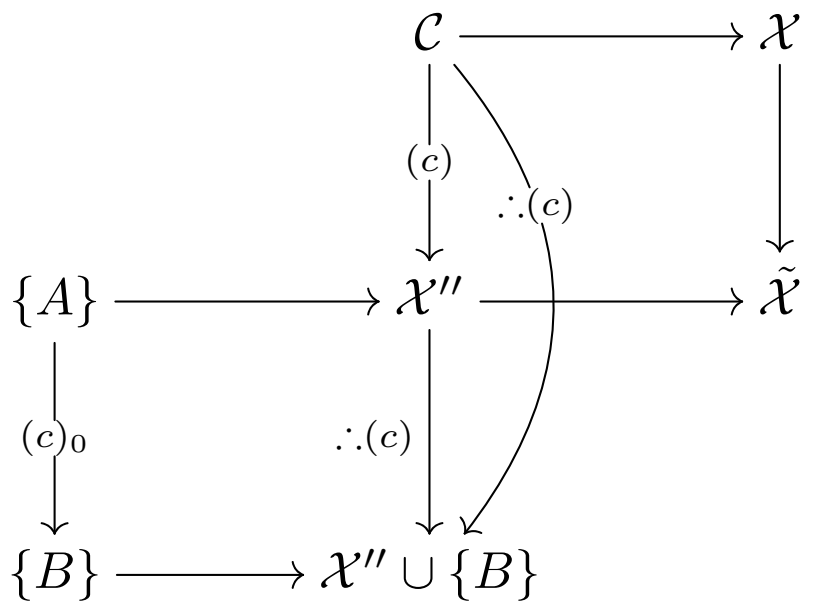

FIGURE 36. The arrow $\mathcal{X}'' \longrightarrow \mathcal{X}'' \cup \{B\}$ is labelled (c) by Axiom (M3)(c). The arrow $\mathcal{C} \longrightarrow \mathcal{X}'' \cup \{B\}$ is labelled (c) as the composition of two (c)-arrows.

Our android wants to show that if $\emptyset \xrightarrow{(c)} \mathcal{Y} \xrightarrow{(wf)} \mathcal{X}'' \cup \{B\}$ then $\mathcal{Y} \longrightarrow \mathcal{X}$, as - by definition - this would imply that $\mathcal{X}'' \cup \{B\} \longrightarrow \tilde{\mathcal{X}}$, which is what he wants.

What he can readily show is that, keeping the above notation, $\mathcal{Y} \longrightarrow \tilde{\mathcal{X}}$. To check this it will suffice that he shows that $\mathcal{Y} \longrightarrow \mathcal{D} \cup \mathcal{B}_c$. Since $\mathcal{D} \xrightarrow{(wf)} \mathcal{X}''$, if we let $\emptyset \xrightarrow{(c)} \mathcal{X}''_c \xrightarrow{(wf)} \mathcal{X}''$ then $\emptyset \xrightarrow{(c)} \mathcal{X}''_c \wedge \mathcal{D} \xrightarrow{(wf)} \mathcal{X}''$ implies that $\mathcal{X}''_c \longrightarrow \mathcal{D}$. Thus, in fact, it will suffice to show that $\mathcal{Y} \longrightarrow \mathcal{X}''_c \cup \mathcal{B}_c$, which is immediate from Proposition 22 (from which it follows that $(\mathcal{X}'' \cup \{B\})_c = \mathcal{X}''_c \cup \mathcal{B}_c$). By Claim 38 $\emptyset \xrightarrow{(c)} \mathcal{Y} \wedge \mathcal{X} \longrightarrow \tilde{\mathcal{X}}$ so $\mathcal{X}'' \cup \{B\} \xleftarrow{(c)} \mathcal{C} \longrightarrow \mathcal{X} \longleftarrow \mathcal{Y} \xrightarrow{(wf)} \mathcal{X}'' \cup \{B\}$ witness that $\mathcal{X}'' \cup \{B\} \longrightarrow \tilde{\mathcal{X}}$, and in particular $\{B\} \longrightarrow \tilde{\mathcal{X}}$, as required.

$\square$

Now, at last, Mr. A. can verify that the category $QtNaamen$ is closed under (finite) limits.

**Claim 40.** The category $QtNaamen$ satisfies Axiom (M0).

*Proof.* Assume that $\mathcal{X}_1, \ldots, \mathcal{X}_n$ are objects in $QtNaamen$ then

$$\bigcap_{i=1}^{n} \mathcal{X}_i := \{\bigcap_{i=1}^{n} X_i : X_i \in \mathcal{X}_i\}$$

is the inverse limit of the $\mathcal{X}_i$ in $QtNaamen$ and $\tilde{\mathcal{D}}$ is their direct limit in $QtNaamen$, where

$$\mathcal{D} := \bigcup_{i=1}^{n} \mathcal{X}_i.$$

Since $\mathcal{I} := \bigcap_{i=1}^{n} \mathcal{X}_i$ is the inverse limit of the $\mathcal{X}_i$ in $StNaamen$ it will suffice to show, in order to prove the first part of the claim, that it is in $QtNaamen$. But this is obvious, due to the following diagram:

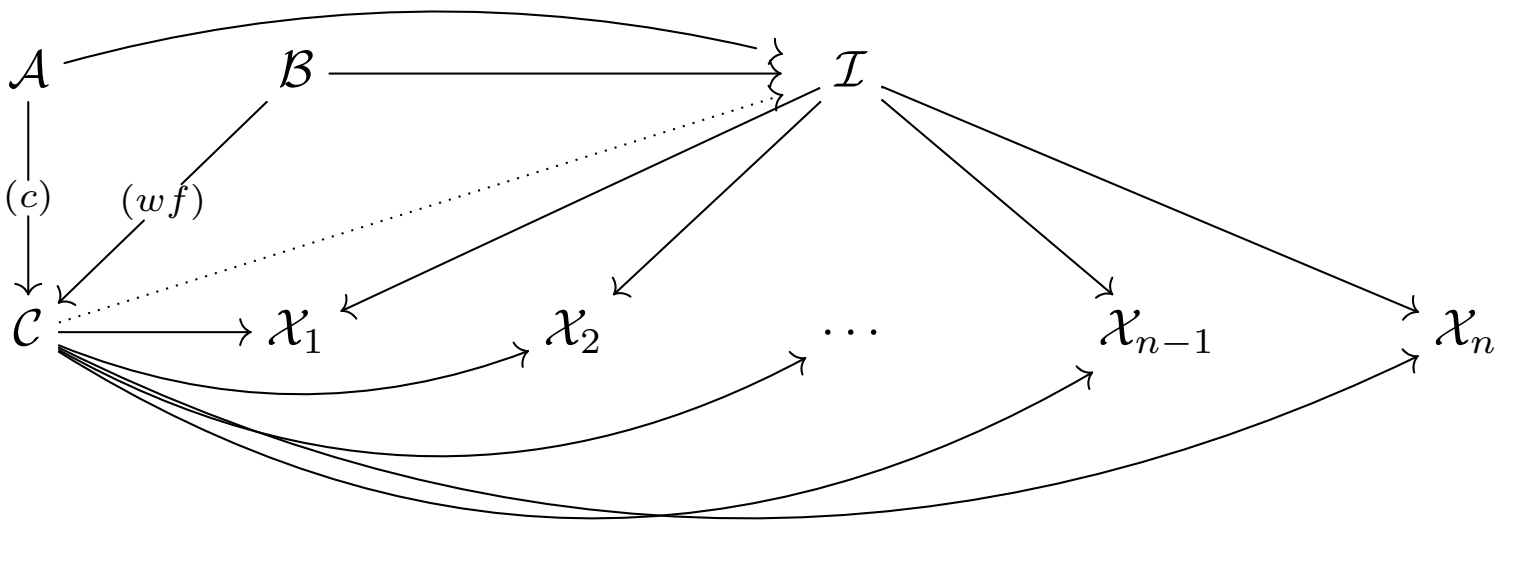

FIGURE 37. The arrows from $\mathcal{C}$ to all the $\mathcal{X}_i$ exists by the definition of $QtNaamen$. The lifting arrow exists by virtue of $\mathcal{I}$ being the inverse limit of the $\mathcal{X}_i$ in $StNaamen$.

For the second part, it will suffice to show that $\tilde{\mathcal{D}}$ is a direct limit (since Mr. A. already knows it to be in $QtNaamen$). Our hero has to to show that whenever $\mathcal{X}$ is an object of $StNaamen$ such that $\mathcal{X}_i \longrightarrow \mathcal{X}$ for all $1 \le i \le n$, there exists an arrow $\tilde{\mathcal{D}} \longrightarrow \mathcal{X}$. But, in *StNaamen*, if $\mathcal{X}_i \longrightarrow \mathcal{X}$ for all $i$ then $\mathcal{D} \longrightarrow \mathcal{X}$. So by Claim 37 the result follows. □

Our hero is rather satisfied. He obtained Axiom (M0), he knows that axioms (M1), (M3) and (M4) are inherited from *StNaamen*. So it only remains to verify Axiom (M2) and Axiom (M5). He is encouraged by the fact that Axiom (M2) seems easy:

**Claim 41.** If $\mathcal{X} \longrightarrow \mathcal{Y}$ in *QtNaamen* and $\mathcal{X} \xrightarrow{(c)} \mathcal{X}_c \xrightarrow{(wf)} \mathcal{Y}$ then $\mathcal{X}_c$ is an object in *QtNaamen*.

*Proof.* He has to check that

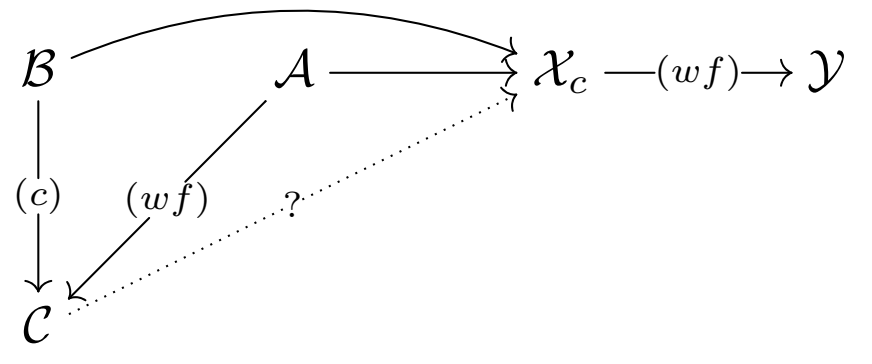

But by the definition of *QtNaamen*, and Axiom (M1) he gets this:

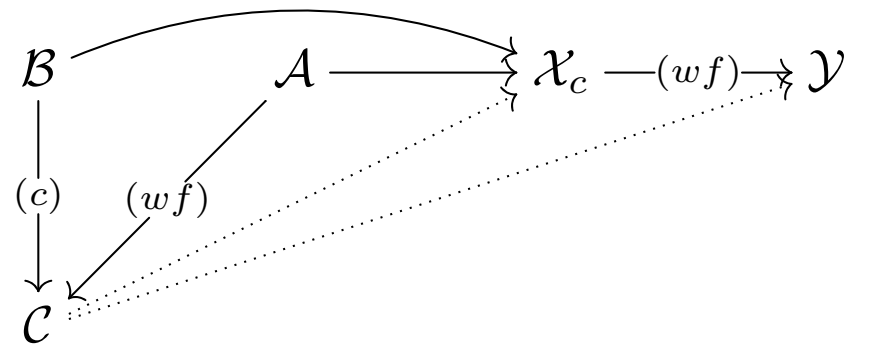

FIGURE 38. The lifting arrow to $\mathcal{Y}$ exists by definition of *QtNaamen*. The arrow to $\mathcal{X}_c$ is obtained by $\mathcal{B} \xrightarrow{(c)} \mathcal{C} \swarrow \mathcal{X}_c \xrightarrow{(wf)} \mathcal{Y}$.

□

The proof of the second part of Axiom (M2) is harder for our hero. In fact, we are not certain it is he who proved it. The key to the proof is:

**Lemma 42** (A continuity fixed-point argument)**.** *The following diagrams are true:*

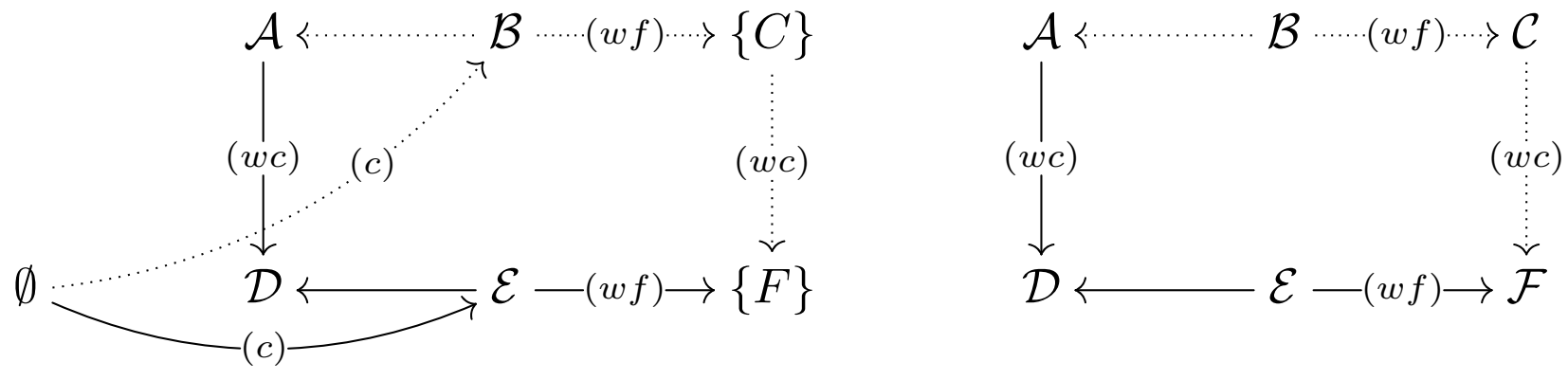

*Proof.* First, we show that the two diagrams are equivalent. Of course, we only have to show that the diagram on the left implies the one on the right. So assume the diagram on the left is true, and we are given $\mathcal{A}$, $\mathcal{D}$, $\mathcal{E}$ and $\mathcal{F}$ as in the figure on the right. First, observe that we may assume without loss that $\mathcal{E} \cong \mathcal{F}_c$, so without loss $\emptyset \xrightarrow{(c)} \mathcal{E}$.

We can now apply the diagram on the left as follows: for every $F \in \mathcal{F}$ observe that $\mathcal{A}$, $\mathcal{D}$, $\mathcal{F}_c \cap \{F\}$ and $\{F\}$ satisy the assumptions of the left hand side diagram. So we can find $\mathcal{B}_F$ and $\{C_F\}$ as provided there. Taking $\mathcal{B} := \bigcup\{B_F : F \in \mathcal{F}\}$ and $\mathcal{C} := \bigcup\{C_F : F \in \mathcal{F}\}$ we get $\emptyset \xrightarrow{(c)} \mathcal{B}$ and by Proposition 22, $\mathcal{B} \xrightarrow{(wf)} \mathcal{C}$, and $\mathcal{C} \xrightarrow{(wc)} \mathcal{F}$ so the diagram on the right is satisfied.

To prove the diagram on the left, first observe that using Proposition 22 it is clear that $\mathcal{A} \xrightarrow{(wc)} \mathcal{D}$ implies that $\mathcal{A}_c \xrightarrow{(wc)} \mathcal{D}_c$. Since $\mathcal{A}_c \longrightarrow \mathcal{A}$ it will suffice to prove the diagram with $\mathcal{A}_c$ instead of $\mathcal{A}$ and $\mathcal{D}_c$ instead of $\mathcal{D}$. So we may assume that $\mathcal{A}$ is a family of countable sets. Thus, combinatorially, the lemma asserts that given a family $\mathcal{A}$ of countable sets and a set $F$ such that for every countable $F_c \subseteq F$ there exists some $A \in \mathcal{A}$ such that $\mathrm{card}\,(F_c \setminus A) < \aleph_0$ then there exists a sets $C \subseteq F$ such that $\mathrm{card}\,(F \setminus C) < \aleph_0$ and every countable subset of $C$ is a subset of an element of $\mathcal{A}$.

Assume by way of contradiction that this is not the case. I.e., we assume that for every finite $\bar{b} \subseteq F$ there exists a countable set $C_{\bar{b}} \subseteq F$ such that $C_b$ is not contained in any element of $\mathcal{A}$. Let $C_0 = C_\emptyset$. Define inductively for $i > 0$

$$C_{i+1} = C_i \cup \{C_{\bar{b}} : \bar{b} \in [C_i]^{<\omega}, C_i \setminus \bar{b} \subseteq A_{i,\bar{b}},\ \text{some } A_{i,\bar{b}} \in \mathcal{A}\}$$

Let $C_\omega := \bigcup_{i\in\omega} C_i$. Then $C_\omega$ is countable (as the countable union of countable sets). But $C_\omega \setminus \bar{b} \not\subseteq A$ for all finite $\bar{b} \subseteq C_\omega$ and $A \in \mathcal{A}$, a contradiction. indeed, note that for all $\bar{b}$ as above there exists $n \in \omega$ such that $\bar{b} \subseteq C_n$. So $C_{\bar{b}} \subseteq C_{n+1} \subseteq C_\omega$, with the desired conclusion. □

For future reference, we point out the following generalization of the fixed point argument, appearing in the last part of the above proof:

**Remark 43.** Let $\kappa$ be a regular cardinal, $\mathcal{A}$ a class of sets of cardinality smaller than $\kappa$, and $F$ any set. Assume that for all $\bar{b} \subseteq F$ with $\mathrm{card}\,(\bar{b}) < \kappa$ there exists $C_{\bar{b}} \subseteq F$ with $\mathrm{card}\,(C_{\bar{b}}) \leq \kappa$ such that $C_{\bar{b}} \not\subseteq A$ for all $A \in \mathcal{A}$. Then there exists $C \subseteq F$ with $\mathrm{card}\, C = \kappa$ such that $\mathrm{card}\,(C \setminus A) = \kappa$ for all $A \in \mathcal{A}$.

The proof of the above remark goes through precisely as in the above lemma, with a single exception. In the fixed point argument above we define $C_0 = C_\emptyset$ and

$$C_{i+1} = C_i \cup \{C_{\bar{b}} : \bar{b} \in [C_i]^{<\omega}, C_i \setminus \bar{b} \subseteq A_{i,\bar{b}},\ \text{some } A_{i,\bar{b}} \in \mathcal{A}\}.$$

If we tried to repeat the same thing replacing systematically $\omega$ with $\kappa$ we might run into trouble, since $\mathrm{card}\, C_0 = \kappa$, but in general, it will no longer be true that $\mathrm{card}\,([C_0]^{<\kappa})$, and so already $C_1$ could be too large. Observe, however, that in order to make the argument go through we need not go over all $\bar{b} \in [C_0]^{<\kappa}$. Indeed, all we need is to go over a *covering* family of such subsets. Namely, the argument would go through unaltered if at stage $i+1$ we fixed $\mathcal{C}_i \subseteq [C_i]^{<\kappa}$ such that for $\bar{b} \in [C_i]^{<\kappa}$ there exists $\bar{c} \in \mathcal{C}_i$ with $\bar{b} \subseteq \bar{c}$. Thus, it suffices to show that we can choose $\mathcal{C}_i$ to be of cardinality at most $\kappa$. The minimum cardinality of such a family $\mathcal{C}_i$ is precisely $\mathrm{cov}(\kappa, \kappa, \kappa, 2)$ (see [She94, §5] for the definition). And the fact that this cardinality is precisely $\kappa$ is an immediate consequence of [She94] Observation 5.2(2), 5.2(5) and Observation 5.3(2).

We conclude the discussion of this last remark by noting that the relation between Shelah's covering numbers $\mathrm{cov}(\lambda, \kappa, \theta, \sigma)$ and various constructs arising form the analysis of the model category *QtNaamen* and some variants is the main theme of the second part of this paper. Interstinly, a somewhat degenerate version of these covering numbers seems to be necessary for the proof of the above remark.

We can now return to our hero and the proof that *QtNaamen* is a model category.

**Lemma 44.** *If* $\mathcal{X} \longrightarrow \mathcal{Y}$ *in QtNaamen and* $\mathcal{X} \xrightarrow{(wc)} \mathcal{X}_f \xrightarrow{(f)} \mathcal{Y}$ *then* $\mathcal{X}_f$ *is an object in QtNaamen.*

*Proof.* Our robot draws a picture:

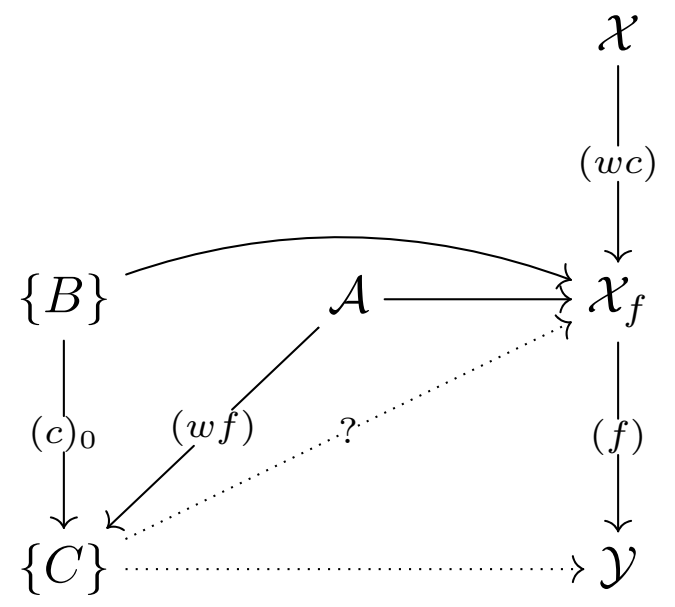

FIGURE 39. The arrow $\{C\} \longrightarrow \mathcal{Y}$ exists because $\mathcal{Y}$ is an object of *QtNaamen*.

With the arrow $\{C\} \longrightarrow \mathcal{Y}$ our hero sees only one way to find the lifting arrow $\{C\} \longrightarrow \mathcal{X}_f$. He has to find an object $\mathcal{B}'$ such that $\mathcal{B}' \xrightarrow{(wc)} \{C\}$ and $\mathcal{B}' \longrightarrow \mathcal{X}_f$. If he can find such an object his work will be done, because $\mathcal{B}' \xrightarrow{(wc)} \{C\} \rightthreetimes \mathcal{X}_f \xrightarrow{(f)} \mathcal{Y}$. Moreover, he sees that it will suffice to make sure that $\mathcal{B}' \longrightarrow \mathcal{X}$, for which he can use the fact that $\mathcal{X}$ is in *QtNaamen*.

In the notation of Figure 39, the system

$$\mathcal{X} \xrightarrow{(wc)} \mathcal{X}_f \longleftarrow \mathcal{A} \xrightarrow{(wf)} \{C\}$$

satisfies the assumption of Lemma 42. So he can add some arrows to Figure 39:

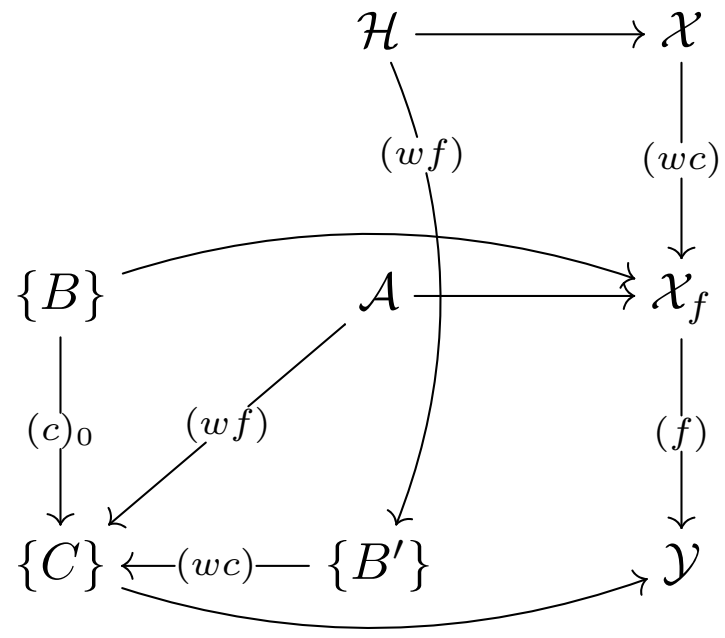

Now, if our robot can find some $\mathcal{I} \longrightarrow \mathcal{X}$ such that $\mathcal{I} \xrightarrow{(c)} \{B'\}$, he will be done. For then, since $\mathcal{X}$ is in $QtNaamen$, he could obtain the lifting arrow $\{B'\} \longrightarrow \mathcal{X}$ he wanted to get. So it remains to find $\mathcal{I}$. But this is now easy: Let $\{B''\}$ be the inverse limit of $\mathcal{X}$ and $\{B\}$. Then by Corollary 24 $\{B''\} \xrightarrow{(wc)} \{B\}$. By the same argument, setting $\mathcal{I}$ to be the inverse limit of $\{B''\}$ & $\{B'\}$ we have $\mathcal{I} \xrightarrow{(c)} \{B'\}$ and $\mathcal{I} \longrightarrow \{B''\} \longrightarrow \mathcal{X}$, which is what we wanted.

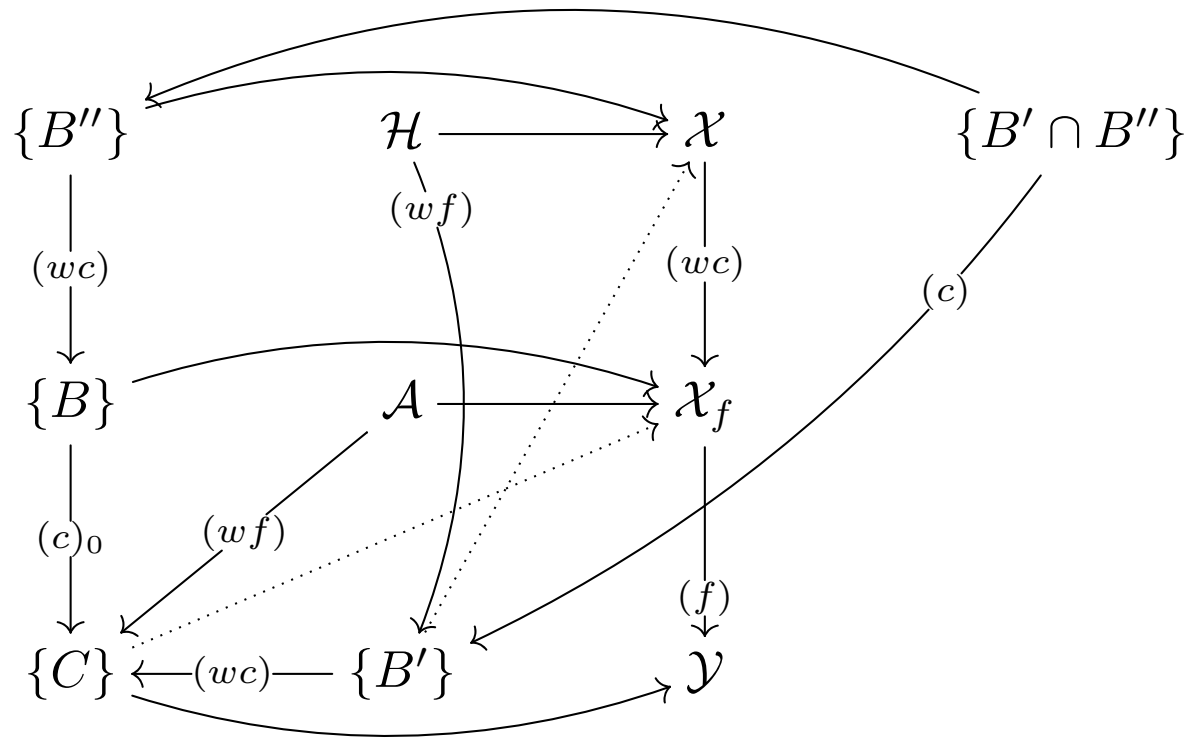

FIGURE 40. The arrow $\{B'\} \longrightarrow \mathcal{X}$ exists because $\mathcal{X}$ is in $QtNaamen$. The arrow $\{C\} \longrightarrow \mathcal{X}_f$ comes from $\{B'\} \xrightarrow{(wc)} \{C\} \nearrow \mathcal{X}_f \xrightarrow{(f)} \mathcal{Y}$

□

Now it is easy for our hero to verify that $QtNaamen$ satisfies Axiom (M6). Though he does not particularly care about this axiom, he knows that it will be considerably simpler for him to show that $QtNaamen$ is a model category if he has a homological characterisation of all the labelling.

**Lemma 45.** *The labelled category QtNaamen satisfies Axiom (M6).*

*Proof.* We have to show the labelling of arrows in $QtNaamen$ is determined by the appropriate lifting properties:

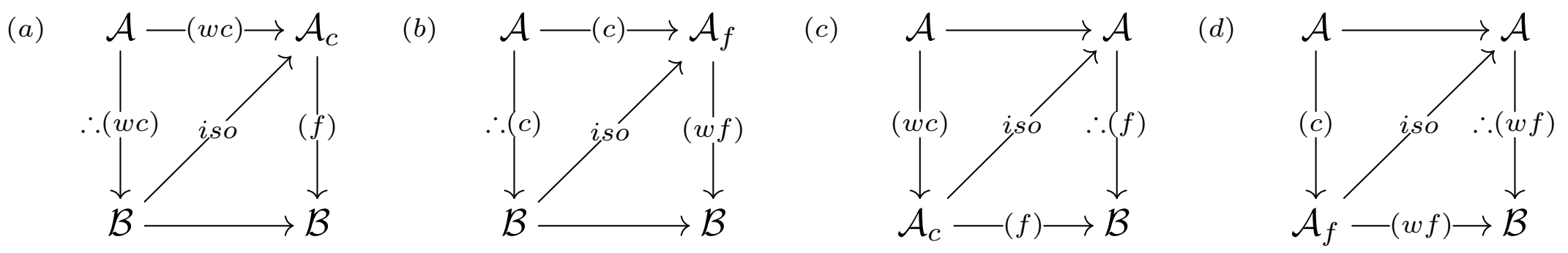

FIGURE 41. (a) Decompose $\mathcal{A} \longrightarrow \mathcal{B} = \mathcal{A} \xrightarrow{(wc)} \mathcal{A}_c \xrightarrow{(f)} \mathcal{B}$. If $\mathcal{A} \longrightarrow \mathcal{B} \rightthreetimes (f)$ it lifts with respect to $\mathcal{A}_c \xrightarrow{(f)} \mathcal{B}$ (because $\mathcal{A}_c$ is in $QtNaamen$). So $B \longrightarrow \mathcal{A}_c$ is an isomorphism. Whence $\mathcal{A} \xrightarrow{(wc)} \mathcal{B}$. (b)-(d) Similar argument.

□

So our hero has only the task of proving Axiom (M5) for $QtNaamen$. He ponders how he should approach the problem. He knows Axiom (M5) asserts, in some sense, that weak equivalence is an equivalence relation. So he looks for an equivalence relation that might capture the notion of weak equivalence. He remembers that by Proposition 22 $\mathcal{X} \xrightarrow{(w)} \mathcal{Y}$ if and only if, for every $X \in \mathcal{X}$ and every $Y' \subseteq Y \in \mathcal{Y}$, if $\operatorname{card} Y' \leq \operatorname{card}(X \cap Y) + \aleph_0$ then $Y'$ is contained, up to finitely many elements in some $X' \in \mathcal{X}$.

As he has already seen, this is not an equivalence relation in $StNaamen$. So he looks for another characterisation which, in $QtNaamen$, may indeed turn out to be an equivalence relation. He notes that the following, countable version of weak equivalence, is a special case:

**Definition 5.** Let $\mathcal{X}, \mathcal{Y}$ be objects in $StNaamen$. Denote $\mathcal{X} \sim^w \mathcal{Y}$ if for every countable set $L$, there exists $X \in \mathcal{X}$ such that $L \setminus X$ is finite if and only if there exists $Y \in \mathcal{Y}$ such that $L \setminus Y$ is finite.

The following is obvious:

**Claim 46.** The relation $\sim^w$ is an equivalence relation. It can be expressed by the following diagram:

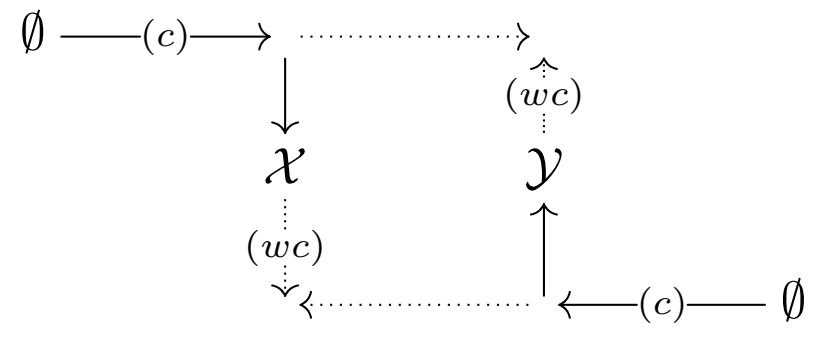

*Proof.* That $\sim^w$ is an equivalence relation is obvious from the definition. We will show that the figure above is equivalent to $\mathcal{X} \sim^w \mathcal{Y}$. First we show that if $\mathcal{X}, \mathcal{Y}$ satisfy the above figure then they are $\sim^w$ equivalent. So assume that $L$ is a countable set such that $L \subseteq X \cup \ell$ for some $X \in \mathcal{X}$ and finite $\ell$. Let $L' := L \setminus \ell$. Then $\emptyset \xrightarrow{(c)} \{L'\} \longrightarrow \mathcal{X}$. Then there exists $\mathcal{Z}$ such that $\mathcal{Y} \xrightarrow{(wc)} \mathcal{Z} \longleftarrow \{L'\}$. So there exists $Z \in \mathcal{Z}$ such that $L' \subseteq Z$. By the combinatorial classification of (cw)-arrows there is $Y \in \mathcal{Y}$ such that $Z \subseteq Y \cup \bar{z}$ for some finite $\bar{z}$. Thus $L \subseteq Y \cup (\ell \cap \bar{z})$, as required. Since the diagram is symmetric the same goes for $\mathcal{Y}$. The other direction is similar. □

So $\mathcal{X} \sim^w \mathcal{Y}$ is an equivalence relation, and if $\mathcal{X} \longrightarrow \mathcal{Y}$ it is a coarsening of $\mathcal{X} \xrightarrow{(w)} \mathcal{Y}$. Our hero tries to see how much coarser this notion could be. To his surprise he gets:

**Lemma 47.** *Assume that $\mathcal{X} \longrightarrow \mathcal{Y}$ and $\mathcal{X} \sim^w \mathcal{Y}$. Then $\mathcal{X} \xrightarrow{(w)} \mathcal{Y}$.*

*Proof.* The best he could hope for is, of course, that $\mathcal{X} \xrightarrow{(w)} \mathcal{Y}$. So he tries to prove this, to see where he will fail. By Axiom (M2), he knows that $\mathcal{X} \xrightarrow{(wc)} \mathcal{V} \xrightarrow{(f)} \mathcal{Y}$. In order to show that $\mathcal{X} \xrightarrow{(w)} \mathcal{Y}$ he has, by Axiom (M6), to show that $\mathcal{V} \xrightarrow{(wf)} \mathcal{Y}$, namely that $(c) \rightthreetimes \mathcal{V} \longrightarrow \mathcal{Y}$. In fact, it will suffice to show that $(c)_0 \rightthreetimes \mathcal{V} \longrightarrow \mathcal{Y}$. So he has to show that

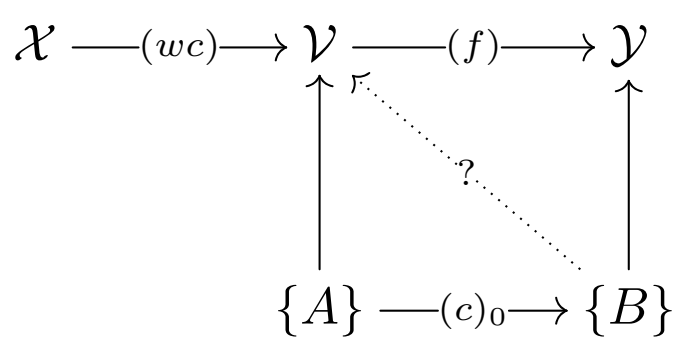

He is, by now, quite used to the category *QtNaamen*, so his first guess is that in order to find the lifting arrow he will have to use the fact that $\mathcal{V}$ is in *QtNaamen*. So he tries to find a copy of Figure 30 with $\mathcal{V}$ in the role of $\mathcal{X}$ there. In order to use the assumption that $\mathcal{X} \sim^w \mathcal{Y}$ he tries to work with countable sets. So he gets the following figure:

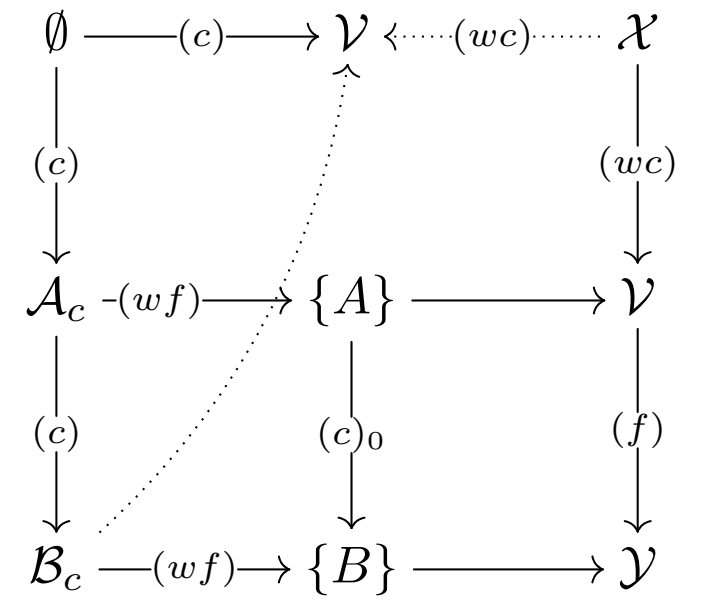

FIGURE 42. Decompose $\emptyset \longrightarrow \{A\}$ and $\emptyset \longrightarrow \{B\}$ by Axiom (M2). Since the composition of two (c)-arrows is a (c)-arrow, and the decomposition provided by Axiom (M2) is unique $\emptyset \xrightarrow{(c)} \mathcal{B}_c \xrightarrow{(wf)} \{B\}$ decomposes into $\emptyset \xrightarrow{(c)} \mathcal{A}_c \xrightarrow{(c)} \mathcal{B}_c \xrightarrow{(wf)} \{B\}$. The dashed arrows exist by the definition of $X \sim^w Y$.

Now he observes that the configuration $\mathcal{X} \xrightarrow{(wc)} \mathcal{V} \longleftarrow \mathcal{B}_c \xrightarrow{(wf)} \{B\}$ fits into the assumptions of Lemma 42. So he concludes that there is a set $B'$ such that:

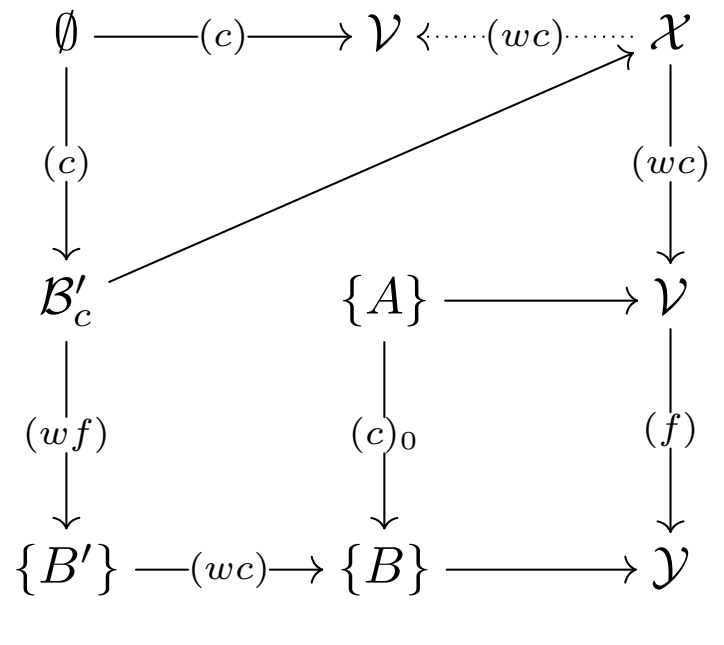

Taking the inverse limit of $\{A\}$ and $\{B'\}$ he gets:

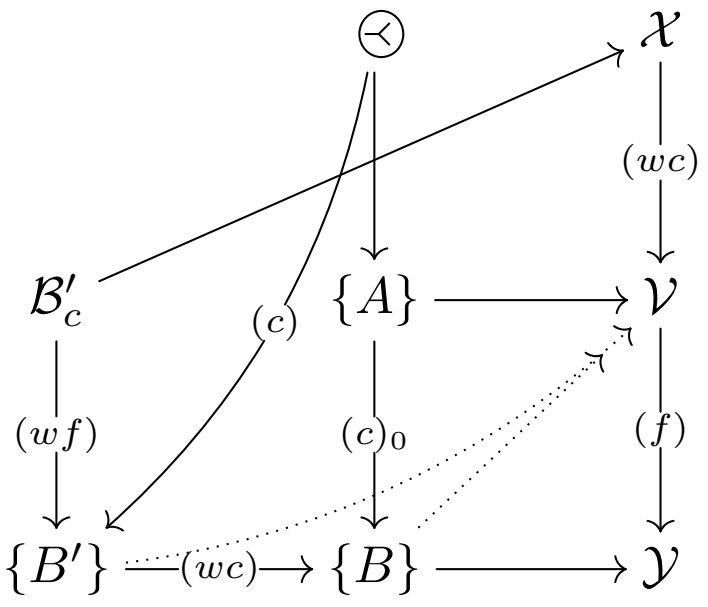

FIGURE 43. The arrow $\{B\}' \longrightarrow \mathcal{V}$ exists, because $\mathcal{V}$ is in *QtNaamen*. The arrow $\{B\} \longrightarrow \mathcal{V}$ exists by $\{B'\} \xrightarrow{(c)} \{B\} \rightthreetimes \mathcal{V} \xrightarrow{(f)} \mathcal{Y}$.

Surprised, our hero discovers that, indeed, he has discovered an equivalence relation equivalent to weak equivalence. □

So our automaton has proved:

**Theorem 48.** *QtNaamen is a closed model category. More precisely, let StNaamen be the category whose objects are classes of sets, and such that for* $\mathcal{C}, \mathcal{D} \in \mathfrak{Ob}StNaamen$ *there is an arrow* $\mathcal{C} \longrightarrow \mathcal{D}$ *if and only if for all* $C \in \mathcal{C}$ *there exists* $D \in \mathcal{D}$ *such that* $C \subseteq D$.

*Let* $\mathfrak{C}_0$ *be the following colouring of the arrows of StNaamen. Label an arrow* $\{A\} \xrightarrow{(c_0)} \{B\}$ *if* $\operatorname{card} A + \aleph_0 = \operatorname{card} B + \aleph_0$ *and* $\{A\} \xrightarrow{(wc)} \{B\}$ *if* $B \setminus A$ *is finite. Let* $\mathfrak{C}_0 \subseteq \mathfrak{C}$ *be the minimal (with respect to inclusion) colouring of StNaamen satisfying Quillen's axioms of model categories (M0),(M1),(M2) and (M6).*

*Let QtNaamen be the full subcategory of StNaamen, whose objects are all those* $\mathcal{X} \in \mathfrak{Ob}StNaamen$ *such that*

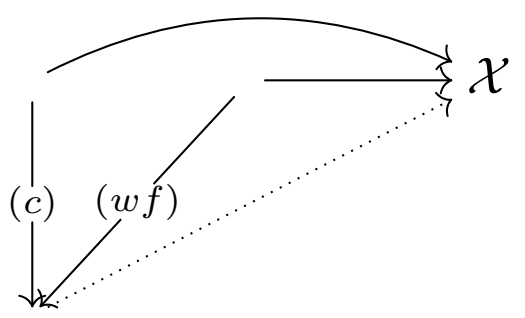

*Then:*

(1) *QtNaamen is a model category.*
(2) *Every co-co-fibrant object of StNaamen is in QtNaamen.*

## References

[Bay09] Martin Bays. *Categoricity Results for Exponential Maps of 1-Dimensional Algebraic Groups & Schanuel Conjectures for Powers and the CIT*. PhD thesis, Oxford University, 2009. Available at http://people.maths.ox.ac.uk/bays/dist/thesis/.

[BZ11] Martin Bays and Boris Zilber. Covers of multiplicative groups of algebraically closed fields of arbitrary characteristic. *Bull. Lond. Math. Soc.*, 43(4):689–702, 2011.

[Dur07] Nikolai Durov. New approach to arakelov geometry. Available on arXiv http://arxiv.org/abs/0704.2030, 2007.

[Gav06] M. Gavrilovich. *Model Theory of the Universal Covering Spaces of Complex Algebraic Varieties Misha Gavrilovich*. PhD thesis, Oxford University, 2006. Available at http://people.maths.ox.ac.uk/bays.

[GH10] M. Gavrilovich and Assaf Hasson. Exercises de style: a homotopy theory for set theory II. Available on `http://corrigenda.ru/by:gavrilovich-and-hasson/what:a-homotopy-theory-for-set-theory/a-homotopy-theory-for-set-theory.pdf`, 2010.

[Gro09] M. Gromov. Structures, learning and ergosystems. Available at http://www.ihes.fr/ gromov/PDF/ergobrain.pdf, 2009.

[Qui67] Daniel G. Quillen. *Homotopical algebra*. Lecture Notes in Mathematics, No. 43. Springer-Verlag, Berlin, 1967.

[She94] Saharon Shelah. *Cardinal arithmetic*, volume 29 of *Oxford Logic Guides*. The Clarendon Press Oxford University Press, New York, 1994. Oxford Science Publications.

[Zil06] Boris Zilber. Covers of the multiplicative group of an algebraically closed field of characteristic zero. *J. London Math. Soc. (2)*, 74(1):41–58, 2006.